%% file: Iso.tex
\documentclass[10pt]{smfart}
\usepackage{a4wide}
\usepackage{amsmath,amssymb,smfthm}
\usepackage[francais]{babel}
\usepackage[latin1]{inputenc}
\usepackage[all]{xy}
\usepackage{graphicx}
\usepackage{epic,eepic}
\usepackage{epsfig}
\usepackage{color}

\def\R{\mathbb{R}}

\def\C{\mathbb{C}}

\def\Cp{\mathbb{C}_p}

\def\Fq{\mathbb{F}_q}
\def\Fqb{\overline{\mathbb{F}}_q}
\def\L{\mathcal{L}}
\def\La{\Lambda}

\def\O{\mathcal{O}} 

\def\I{\mathcal{I}}

\def\*{^\times }
\def\dpt{\displaystyle}

\def\Gm{\mathbb{G}_m}         
\def\l{\lambda}

\def\a{\alpha}
\def\b{\beta}

\def\s{\sigma} 
\def\ph{\varphi}

\def\lssi{\Longleftrightarrow}

\def\limpl{\Longrightarrow}
\def\drt{\rightarrow}
\def\ldrt{\longrightarrow}
\def\Q{\mathbb{Q}}

\def\Qp{\mathbb{Q}_p}
\def\Qpb{\overline{\mathbb{Q}}_p}
\def\Zp{\mathbb{Z}_p}
\def\Z{\mathbb{Z}}
\def\N{\mathbb{N}}

\def\Hom{\text{Hom}}

\def\Gal{\text{Gal}}
\def\={\! = \!}

\def\spec{\text{Spec}}
\def\spf{\text{Spf}}

\def\E{\mathcal{E}}

\def\limp{\underset{\longleftarrow}{\text{ lim }}\;}
\def\limi{\underset{\longrightarrow}{\text{ lim }}\;}

\def\iso{\xrightarrow{\;\sim\;}}

\def\End{\text{End}}
\def\Aut{\text{Aut}}
\def\GL{\hbox{GL}}

\def\xrig{\xrightarrow}

\def\X{\mathfrak{X}}
\def\GG{\Gamma}

\def\bc{\backslash}
\def\AA{\mathcal{A}}

\def\<{<\hspace{-1mm}}
\def\>{\hspace{-1mm}>}

\def\Lie{\text{Lie}}

\def\dem{{\it Démonstration. }}

\def\Fil{\text{Fil}}

\def\unp{[ \frac{1}{p}]}
\def\LT{\mathcal{L}\mathcal{T}}

\def\Hb{\mathbb{H}}
\def\Gb{\mathbb{G}}

\def\DD{\mathbb{D}}
\def\unpi{[\frac{1}{\pi}]}

\author{  Laurent Fargues} 
\address{CNRS-université Paris-Sud-IHES} 
\email{laurent.fargues@math.u-psud.fr} 
\date{}

\begin{document}

\title[L'isomorphisme entre les tours de Lubin-Tate et de Drinfeld]{L'isomorphisme entre les tours de Lubin-Tate et de Drinfeld
  : démonstration du résultat principal}
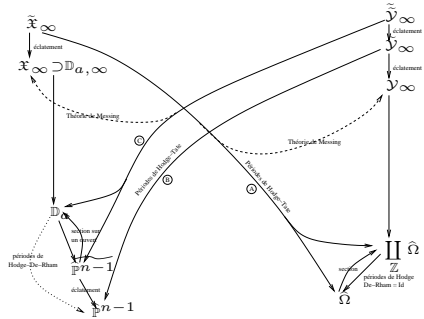
\maketitle

\vspace{-1cm}
\begin{figure}[h]
   \begin{center}
      \input{desol.pstex_t}
   \end{center}
\end{figure}

\begin{abstract}
Il s'agit du quatrième article consacré à l'isomorphisme
entre les tours de Lubin-Tate et de Drinfeld. Nous y démontrons le
résultat final concernant l'isomorphisme c'est à dire l'existence d'un isomorphisme
équivariant entre un
certain éclaté du schéma formel associé à l'espace de Lubin-Tate en
niveau infini construit dans \cite{Cellulaire} et un éclaté d'un autre
schéma formel construit à partir de l'espace de Drinfeld. Nous suivons pour cela
 la stratégie fournie par Faltings dans \cite{Faltings8}.  
\end{abstract}

\begin{altabstract}
This is the fourth article about the isomorphism between Lubin-Tate
and Drinfeld towers. We prove the final result concerning the
isomorphism that is to say the existence of an equivariant isomorphism
between some blow-up of the formal scheme associated to Lubin-Tate
space with infinite level constructed in \cite{Cellulaire} and a
blow-up of another formal scheme associated to Drinfeld tower. We
follow the strategy provided by Faltings in \cite{Faltings8}.  
\end{altabstract}

\section*{Introduction}

Dans l'article \cite{Cellulaire} on a construit un schéma formel
$\pi$-adique lié à l'espace de Lubin-Tate ``en niveau infini''. 
Dans la première section de cet article on construit un schéma formel
$\pi$-adique analogue du coté de l'espace de Drinfeld. Sa construction
est beaucoup plus simple que celle de \cite{Cellulaire}. Il s'agit du
normalisé du schéma formel de Deligne-Drinfeld $\widehat{\Omega}$ dans
les revêtements de sa fibre générique rigide donnés par les structures
de niveau sur le $\O_D$-module formel spécial universel sur
$\widehat{\Omega}$ (en fait, étant donné que nous adoptons le point de
vue des espace de Rapoport-Zink (\cite{RZ}) il s'agit d'une union disjointe
indexée par $\Z$ de
tels espaces).   

Le but est maintenant de construire un isomorphisme
équivariant entre des éclatements formels admissibles de ces deux
schémas formels. Pour cela on construit des morphismes dans les deux
directions (Lubin-Tate vers Drinfeld et réciproquement) puis on vérifie
qu'ils sont inverses l'un de l'autre. 

Avant d'expliquer le principe de la construction introduisons un point
clef : les deux types de périodes vivant sur ces espaces. 

\subsection*{Périodes de Hodge De-Rham}

Il s'agit des morphismes de périodes étudiés dans le cas de l'espace
de Lubin-Tate dans \cite{HopkinsGross} (cf. également \cite{Cellulaire}) et
définis en toutes généralités dans \cite{RZ}. Il s'agit de l'analogue
du plongement $X\hookrightarrow \check{X}$ où $X$ est un espace
symétrique hermitien de dual compact $\check{X}$. 

Le principe est le suivant : si $R$ est une 
$W(\overline{\mathbb{F}}_p )$-algèbre $p$-adique  et $H$
est un groupe $p$-divisible sur $\spf (R)$ il y a une filtration de
De-Rham sur l'homologie de De-Rham de $H$
$$
0 \ldrt V(H) \ldrt \Lie\, E(H) \ldrt \Lie\, H \ldrt 0
$$
où $E(H)$ est l'extension vectorielle universelle de $H$ et $V(H)$ sa
partie vectorielle. Supposons nous donné un groupe $p$-divisible $H_0$
sur $\overline{\mathbb{F}}_p$ et une quasi-isogénie 
$$
\rho : H_0\otimes_{\overline{\mathbb{F}}_p} R/p R \ldrt H\otimes_{R} R/pR
$$
Celle ci fournit grâce à la nature cristalline de l'extension
vectorielle universelle (\cite{Messing1})
couplée aux puissances divisées sur l'idéal $pR$  une rigidification de
l'homologie de De-Rham 
$$
\rho_* : \DD (H_0) \otimes_{W(\overline{\mathbb{F}}_p)} R\unp \iso \Lie\, E(H) \unp
$$
où $\DD (-)$ désigne le module de Dieudonné covariant. 
Et donc la suite précédente donne après inversion de $p$ un morphisme
de $\spec (R\unp)$ 
vers une Grassmanienne associée à l'espace $\DD (H_0)\unp$. 

\subsection*{Périodes de Hodge-Tate}

Reprenons les notations précédentes. Il y a un morphisme 
$$
\a_H : \Hom ( \Qp/\Zp, H) \ldrt \omega_{H^D}
$$
où $H^D$ désigne le dual de Cartier de $H$. \`A $\chi:\Qp/\Zp \ldrt H$
on associe $(\chi^D)^* \frac{dT}{T}$ où $\chi^D : H^D \ldrt
\mu_{p^\infty}$ et $(\chi^D)^* : R.\frac{dT}{T}
=\omega_{\mu_{p^\infty}} \ldrt \omega_{H^D}$. Notons $T_p (H) =\Hom
(\Qp/\Zp ,H)$ qui si $R$ est sans $p$-torsion 
intégralement fermé dans $R\unp$
s'identifie à $\Hom (\Qp/\Zp, H\otimes_R R\unp)$ et ne dépend donc
que du groupe étale fibre générique $H\otimes_R R\unp$. Supposons de
plus que $\mu_{p^\infty} (\Qpb) \subset R$
 Il y a alors un
accouplement parfait 
$$
T_p (H) \times T_p (H^D)\ldrt \Zp (1)
$$
et l'on a donc une suite 
$$
0\ldrt \omega_{H}^* (1) \xrig{\;\,^t \a_{H^D} (1)\;} T_p
(H)\otimes_{\Zp} R \xrig{\; \a_H\;} \omega_{H^D}  \ldrt 0
$$
Sous certaines hypothèses (et c'est là une des difficultés majeures
de cet article) on peut rendre cette suite exacte après inversion de
$p$. Si maintenant on se donne une rigidification du module de Tate 
$\Zp^n \iso T_p (H)$ on obtient une suite exacte de Hodge-Tate
$$
0\ldrt  \omega_{H}^* \unp (1) \xrig{\;\,^t \a_{H^D} (1)\;} R\unp^n \xrig{\; \a_H\;} \omega_{H^D}\unp  \ldrt 0
$$
qui fournit donc un morphisme de $spec (R\unp)$ vers une Grassmanienne 
associée à l'espace vectoriel $\Qp^n$.

\subsection*{Plan de la démonstration}

\begin{figure}[h]
   \begin{center}
      \input{diagramme_final.pstex_t}
   \end{center}
\caption{\footnotesize Le schéma de la démonstration}
\label{lsfjatff}
\end{figure}
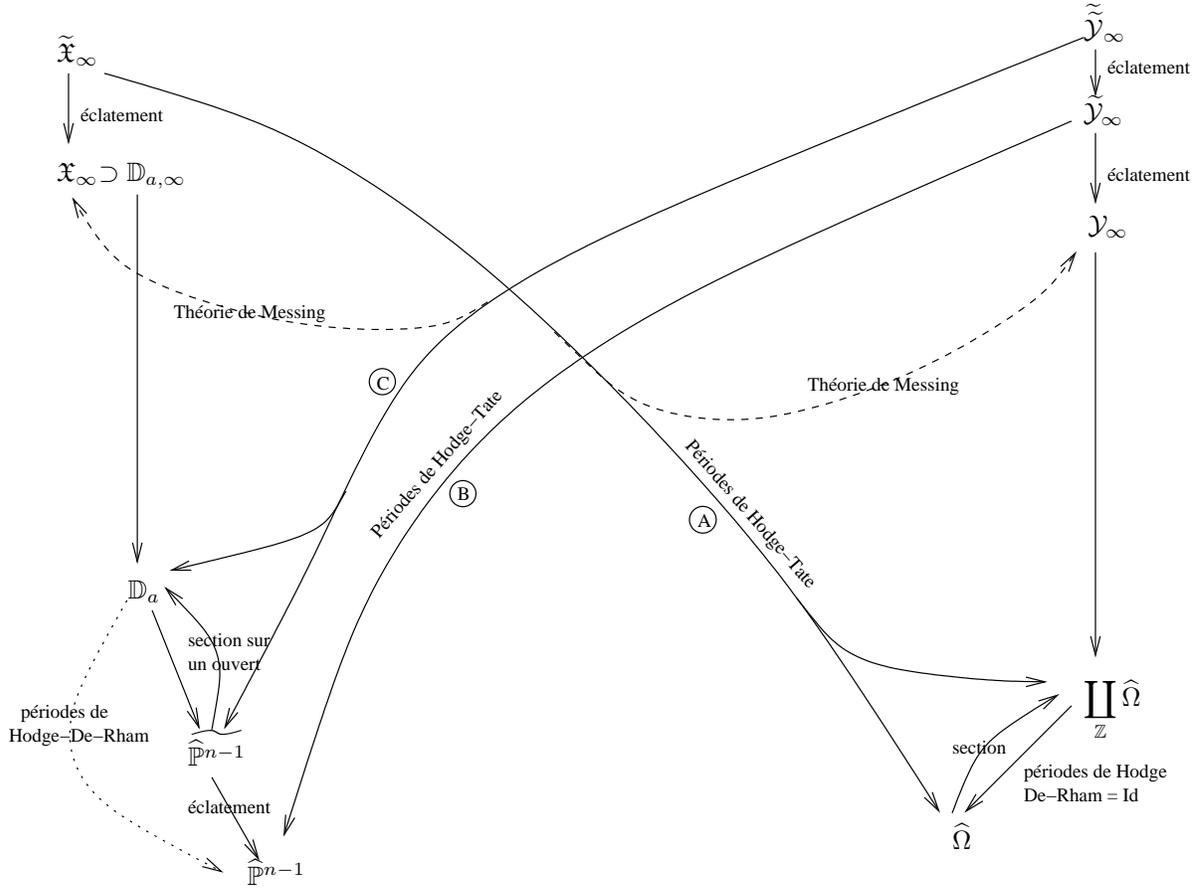

Le principe général est le même dans les deux sens. On note
$\X_\infty$ le schéma formel associé à l'espace de Lubin-Tate
construit dans \cite{Cellulaire}.  
\begin{itemize}
\item On construit un morphisme de périodes de Hodge-Tate 
 d'un espace en niveau infini vers
  l'espace des périodes de Hodge-De-Rham de l'autre espace, le
  principe étant que via l'isomorphisme entre les deux tours les
  périodes de Hodge-De-Rham et de Hodge-Tate sont permutées. 

L'espace des périodes de Hodge-De-Rham associé à l'espace de Drinfeld
$\widehat{\Omega}$ 
est l'espace de Drinfeld lui-même : l'application des périodes de
Hodge-De-Rham est un isomorphisme sur son image pour cet espace (et possède donc un
modèle entier (i.e. sans inverser $p$) donné par l'identité). 

 Ainsi dans la section 2 on
  construit un morphisme de périodes de Hodge-Tate
(noté $(A)$ dans la figure \ref{lsfjatff})
 de l'espace de Lubin-Tate ``en niveau infini''
  éclaté $\widetilde{\X}_\infty$ vers l'espace de Drinfeld
  $\widehat{\Omega}$. Pour définir ce morphisme nous  avons besoin d'éclater (et normaliser) le
  schéma formel $\X_\infty$ construit dans \cite{Cellulaire} pour deux raisons :
\begin{itemize}
\item La première est que l'on doit platifier par éclatements l'image
  d'une certaine application de Hodge-Tate afin de définir une version
  entière de celle-ci (la suite de Hodge-Tate précédente ne peut être exacte
  qu'après inversion de $p$)
\item La seconde vient du fait qu'une fois qu'on a platifié la suite
  de Hodge-Tate on a un morphisme d'un certain éclaté de $\X_\infty$
  vers un espace projectif, mais $\widehat{\Omega}$ est lui-même obtenu
  par éclatement de l'espace projectif. Il faut donc tirer en arrière
  de tels éclatements vers $\X_\infty$ afin de relever le morphisme
  défini vers l'espace projectif à $\widehat{\Omega}$. Cela est rendu
  nécessaire par le fait que la décomposition cellulaire de
  $\X_\infty$ construite dans \cite{Cellulaire} et indexée par un
  immeuble de Bruhat-Tits ne correspond pas à la décomposition
  cellulaire usuelle de $\widehat{\Omega}$ associée à l'immeuble du
  groupe linéaire. En effet, on a vu dans \cite{Rami} que l'image du
  domaine fondamental de Gross-Hopkins (qui vit dans l'espace de
  Lubin-Tate) dans $\Omega$ est un domaine polyèdral qui n'apparaît
  pas dans la structure simpliciale de l'immeuble, et donc n'apparaît
  pas dans la fibre spéciale de $\widehat{\Omega}$.
\end{itemize}

L'espace des périodes de Hodge-De-Rham associé à l'espace de
Lubin-Tate est le schéma formel $\widehat{\mathbb{P}}^{n-1}$. 
De la même façon que précédemment on définit 
dans la section \ref{kfjeg264} un morphisme de périodes
de Hodge-Tate 
(noté $(B)$) d'un éclaté $\widetilde{\mathcal{Y}}_\infty$  
du schéma formel associé à l'espace de Drinfeld en niveau infini noté
$\mathcal{Y}_\infty$ vers l'espace des périodes de Hodge-De-Rham de
l'espace 
Lubin-Tate $\widehat{\mathbb{P}}^{n-1}$. L'éclatement sert à 
platifier la suite de Hodge-Tate. 

\item Après avoir construit ces deux morphismes de périodes de
  Hodge-Tate on les relève aux espaces de modules sans structure de
  niveau. 
\begin{itemize}
\item
Dans le cas de l'espace de Drinfeld c'est immédiat puisque
  l'application des périodes est l'identité (ou plutôt la projection 
$\coprod_\Z \widehat{\Omega}\ldrt \widehat{\Omega}$ puisque l'on
considère les espaces de Rapoport-Zink). 
On peut donc relever le
morphisme noté précédemment $(A)$ en un morphisme
$\widetilde{\X}_\infty \ldrt \coprod_\Z \widehat{\Omega}$. 
\item
Dans le cas de l'espace de Lubin-Tate cela s'avère plus complexe. 
Commençons par rappeler que le schéma formel $\X_\infty$ a été
construit comme recollement de cellules (des schémas formels
$\pi$-adiques affines) ``en niveau infini''
 $\DD_{a,\infty}$ où l'indice $a$ varie dans les
sommets d'un immeuble de Bruhat-Tits. La cellule $\DD_{a,\infty}$ est
au dessus d'une cellule sans niveau $\DD_a$.
On utilise le fait (\cite{HopkinsGross}, \cite{Cellulaire}) que l'application
des périodes rigide $\DD_a^{rig} \ldrt (\mathbb{P}^{n-1})^{rig}$ est
un isomorphisme sur son image, un ouvert admissible dans
$(\mathbb{P}^{n-1})^{rig}$. Quitte à former un éclatement formel
admissible $\widetilde{\widehat{\mathbb{P}}^{n-1}}\ldrt
\widehat{\mathbb{P}}^{n-1}$ on peut faire apparaître cet ouvert dans
la fibre spéciale de $\widetilde{\widehat{\mathbb{P}}^{n-1}}$. Alors
l'application des périodes de Hodge-De-Rham devient entière sur
$\DD_a$ et induit un isomorphisme entre $\DD_a$ et cet ouvert. 
Tirant en arrière l'éclatement $\widetilde{\widehat{\mathbb{P}}^{n-1}}\ldrt
\widehat{\mathbb{P}}^{n-1}$ sur $\widetilde{\mathcal{Y}}_\infty$ on
obtient un éclatement $\widetilde{\widetilde{\mathcal{Y}}}_\infty
\ldrt \widetilde{\mathcal{Y}}_\infty$ tel que sur cet éclaté le morphisme
de périodes de Hodge-Tate $(B)$ s'étende en un morphisme noté $(C)$
vers $\widetilde{\widehat{\mathbb{P}}^{n-1}}$. Grâce au fait que le
morphisme des périodes est un isomorphisme entre $\DD_a$ et un ouvert
de    $\widetilde{\widehat{\mathbb{P}}^{n-1}}$ on a une section sur
cet ouvert qui permet de relever le morphisme $(C)$ (c'est un peu plus
compliqué que cela car la section n'existe que sur un ouvert...) 
de $\widetilde{\widetilde{\mathcal{Y}}}_\infty$ vers la cellule sans
niveau $\DD_a$. Tout cela est fait dans la section \ref{lqvout159}. 
\end{itemize}
\item Il reste maintenant à relever nos deux morphismes dont le but
  est un des schémas formels sans niveau vers le schéma formel en
  niveau infini. Pour cela il faut construire des éléments dans le
  module de Tate de l'image réciproque par notre morphisme du  groupe $p$-divisibles
  universel sur l'espace au but du morphisme.
On utilise pour cela la théorie de Messing. Celle-ci permet de
transférer des éléments du module de Tate d'un groupe $p$-divisible
vers un autre en les transférant modulo $p$ puis en utilisant le
critère de relèvement de Messing (le fait qu'un morphisme induit au
niveau des cristaux est compatible aux filtrations de Hodge). 
Cela est fait dans la section \ref{jvhmzohr3} pour le morphisme de Lubin-Tate
vers Drinfeld et dans la section \ref{dsjdmi184} dans l'autre sens.
\item Enfin dans la section \ref{consteug27} on éclate de nouveaux
  idéaux (et normalise) afin de construire les deux schémas formels
  finaux éclatés isomorphes. Pour cela on développe dans l'appendice 
une théorie des
  éclatements formels admissibles et normalisation dans la fibre
  générique 
 pour des schémas formels $\pi$-adiques ne satisfaisant aucune
 condition de finitude.  Les résultats de l'appendice permettent
 également de voir que les éclatements/normalisation effectués durant la
 démonstration peuvent se faire directement en niveau infini. 
\end{itemize}

\hspace{8mm}

{\it Prérequis : La lecture des articles
  \cite{Cellulaire} et \cite{Points} est indispensable puisque cet
  article leur fait suite. Par contre nous n'utiliserons  quasiment
  pas l'article
  \cite{Rami} (nous l'utilisons dans la section \ref{dkluiyar16},
  cependant la section 2.2 de \cite{Points} serait suffisante mais
  donnerait lieu à des éclatements supplémentaires). Nous
  utilisons la théorie de \cite{Faltings7} des ``strict $\O$-actions''
  afin de ``simplifier'' les notations. Le lecteur ne connaissant pas
  \cite{Faltings7} pourra supposer $F=\Qp$.
  Outre la
  connaissance des diverses références déjà citées dans
  les prérequis de \cite{Cellulaire} et \cite{Points} le lecteur doit
  être familier avec la théorie de Raynaud des éclatements formels admissibles telle
  qu'exposée dans \cite{BLI}. Bien qu'elle n'aparaisse pas dans les démonstrations,
  la théorie de Fontaine est au coeur de cet isomorphisme et a toujours
  été présente à l'esprit de l'auteur (cf. \cite{Points}). 
}
\\

{\it Enfin, il apparaitra comme clair au lecteur que l'auteur de cet article s'est 
largement inspiré des travaux de Faltings \cite{Faltings8}}.

\section*{Notations}

On fixe $F$ une extension de degré fini de $\Qp$. On note $\O$ son anneau des
entiers et $\Fq$ son corps résiduel. Enfin on fixe une clôture algébrique $\Fqb$ de $\Fq$
et l'on note  $\breve{\O} = W_\O (\Fqb)$. 

Fixons $\LT$ un groupe de Lubin-Tate de hauteur $1$ 
relativement à l'extension $F|\Qp$
sur $\breve{\O}$.
 On note alors  
$F/\O_F (1) = \LT [\pi^\infty] (\O_{\overline{F}})$ et $\O_F(1)
 = T_p (\LT)$
  son module de Tate. Enfin on fixe un générateur
 $\b_{\LT}$ du $\breve{\O}$-module de rang $1$ $\omega_{\LT}$ des
 formes différentielles invariantes sur $\LT$.  
\\

Soit $\Hb$ un $\O$-module $\pi$-divisible formel de dimension $1$ et hauteur $n$  sur $\Fqb$.
Soit $D$ une algèbre à division d'invariant $1/n$ sur $F$. On identifie alors $\O_D$ à
$\End (\Hb)$. 
On fixe un $\O_D$-module $\pi$-divisible
 formel spécial $\mathbb{G}$ de hauteur $n^2$  sur $\Fqb$ muni d'une isogénie $\O_D$-équivariante 
$$
\Delta : \Hb^n \ldrt \Gb 
$$
(on renvoie à \cite{Points}). On note $\DD (\Hb)$, resp. $\DD (\Gb)$,
les modules de Dieudonné covariants de $\Hb$, resp. $\Gb$, relativement à l'extension $F|\Qp$ (on renvoit à l'appendice B de \cite{Cellulaire}). Ce sont des $\breve{\O}$-modules de rang $n$, resp. $n^2$, munis d'un endomorphisme $\s^{-1}$-linéaire 
$V$, le Verschiebung, où $\s$ désigne le Frobenius de $W_\O (\Fqb)$. 

Soit $\Pi$ une uniformisante de $\O_D$. Soit $F_n$ une extension non-ramifiée de degré $n$ de $F$ dans $D$. La bijection $\O_D \iso \End (\Hb)$ induit un morphisme $\O_{F_n}\ldrt \Fqb$ par action 
de $\O_{F_n}$ sur l'algèbre de Lie de $\Hb$. Cela permet d'identifier le corps résiduel de $F_n$ à $\mathbb{F}_{q^n}$ et fournit un plongement $F_n \hookrightarrow \breve{F}$. 

Fixons un isomorphisme 
$$
F\iso \mathbb{D} (\Hb)_{\Q,0}^{V=\Pi}
$$
où $\DD (\Hb)_{\Q,0}$ désigne le facteur direct de $\DD (\Hb)_\Q$ où
$\O_{F_n} \subset \O_D=\End (\Hb)$ agit via le plongement
$F_n\hookrightarrow \breve{F}$. 

Via $\Delta$ cela induit un isomorphisme 
$$
F^n \iso \DD (\Gb)_{\Q,0}^{V=\Pi}
$$
Cela nous permet d'identifier
$$
\End_{\O_D} (\Gb)_\Q = \End_{\O_D} (\DD (\Gb)_\Q,V )= \End ( \DD ( \Gb)_{\Q,0}^{V=\Pi}) = M_n (F) 
$$

Pour un groupe $M$ on note $\underline{M}$ le faisceau constant
associé.

\section{La bestiole du coté Drinfeld}

Dans l'article \cite{Cellulaire} on a construit un schéma formel
$p$-adique $\X_\infty$ sur $\spf (\breve{\O})$ muni d'une action de
$\GL_n (F)\times D^\times$.
Nous construisons maintenant son équivalent du coté Drinfeld.
\\

Rappelons qu'un schéma formel localement de type fini plat sur $\spf (\breve{\O})$
est un schéma formel localement de la forme $\spf (A)$ où $A$ est une
$\breve{\O}$-algèbre $\pi$-adiquement complète topologiquement de type
fini sans $\pi$-torsion c'est à dire un quotient sans $\pi$-torsion
d'une algèbre du type $\breve{\O}<T_1,\dots,T_n>$, les séries formels
$\pi$-adiquement convergentes. Rappelons également que pour ce type de
schémas formels on a une bonne notion de schéma formel normal ainsi que de normalisé
(on renvoie à l'appendice A de \cite{Cellulaire}).

Rappelons également qu'on appelle schéma formel $\pi$-adique
(sous entendu sur $\spf (\breve{\O})$) un schéma formel $\mathfrak{Z}$ pour lequel 
$\pi\O_{\mathfrak{Z}}$ est un idéal de définition. La catégorie des schémas formels $\pi$-adiques est donc équivalente à la 2-limite projective de la catégorie des schémas sur 
 $(\spec (\breve{\O}/\pi^k\breve{\O}))_{k\geq 1}$, c'est à dire les
 familles $(Z_k)_{k\geq 1}$ munies de données de `réduction 
 où $Z_k$ est un schéma sur $\spec (\breve{\O}/\pi^k\breve{\O})$ et la donnée de
réduction est un ensemble d'isomorphismes $Z_{k+1}\otimes \breve{\O}/\pi^k\breve{\O} \iso Z_k$ satisfaisant des conditions de cocyles évidentes. 

Enfin tous les schémas formels considérés seront quasi-séparés. 

\begin{lemm}
Soit $K$ un sous-groupe compact ouvert de $\O_D^\times$. Le
foncteur défini sur la catégorie des schémas formels localement de
type fini sur $\spf (\breve{\O})$ plats et normaux qui à
$\mathfrak{Z}$ associe l'ensemble des classes d'isomorphismes de
triplets $(G,\rho,\eta)$ où
\begin{itemize}
\item $G$ est un $\O_D$-module formel spécial sur $\mathfrak{Z}$
\item $\rho : \mathbb{G}\times_{\Fqb} (\mathfrak{Z} \text{ mod } \pi )
\ldrt G\times_{\mathfrak{Z}} (\mathfrak{Z} \text{ mod } \pi)$ est une
quasi-isogénie $\O_D$-équivariante 
\item $\eta : \O_D \iso T_p (G^{rig})\; [K]$ est une structure de
  niveau $K$ sur $\mathfrak{Z}^{rig}$
\end{itemize}
est représentable. 
\end{lemm}
\dem Si $K=\O_D^\times$ le foncteur classifie les couples $(G,\rho)$
modulo isomorphisme. On sait alors qu'il est représentable (non-canoniquement) par le
schéma formel 
$$
\coprod_{\Z} \widehat{\Omega}
$$
où $\widehat{\Omega}$ est le schéma formel de Deligne-Drinfeld sur $\spf (
\breve{\O})$ et la composante indexée par l'entier $h$ dans l'union
disjointe sur $\Z$ classifie les couples $(G,\rho)$ avec $\text{ht}\,
\rho =nh$ (par hauteur on entend hauteur au sens des $\O$-modules $\pi$-divisibles).  
\\
Pour $K$ général l'espace classifiant des structures de niveau $K$ 
sur l'espace rigide fibre générique de l'espace précédent 
est un espace rigide 
$$
\coprod_{\Z} \Omega_K
$$
étale fini au dessus de 
$$
\coprod_{\Z} \widehat{\Omega}^{rig}
$$
D'après l'appendice A.3 de \cite{Cellulaire} 
 le foncteur de l'énoncé est donc représentable par le
normalisé de $\coprod_{\Z} \widehat{\Omega}$ dans $\coprod_{\Z}
\Omega_K$.
\qed

\begin{rema}
Avec les notations de la démonstration précédente rappelons que
si $$\Pi_K : \Omega_K \ldrt \widehat{\Omega}^{rig}$$ est le morphisme
d'oubli de la structure de niveau, $$\text{sp} :
\widehat{\Omega}^{rig} \ldrt \widehat{\Omega}$$
le morphisme de topos annelés de ``spécialisation'' et
$\O^0_{\Omega_K}\subset \O_{\Omega_K}$ le sous-faisceau des fonctions
rigides 
$$
\O^0_{\Omega_K} = \{ f\in \O_{\Omega_K}\;|\; \|f\|_\infty \leq 1 \}
$$
alors $ \text{sp}_* \Pi_{K*}\O^0_{\Omega_K}$ est une
$\O_{\widehat{\Omega}}$-algèbre cohérente et si 
$$
\widehat{\Omega}_K =\spf \left ( \text{sp}_* \Pi_{K*}\O^0_{\Omega_K} \right ) 
$$
alors 
$\widehat{\Omega}_K$ est un schéma formel localement de type fini sur
$\spf ( \breve{\O})$ et $\coprod_{\Z} \widehat{\Omega}_K$ représente
le schéma formel du lemme précédent.
\end{rema}

\begin{defi}\label{zer3476}
On note $\mathcal{Y}_K$ le schéma formel précédent. 
On obtient ainsi une tour de schéma formels $\pi$-adiques dont les
morphismes de transition sont finis et qui 
 est munie d'une
action de $\GL_n (F)\times \O_D^\times$ via 
$$
\forall (g,d)\in \GL_n (F) \times \O_D^\times \;\; (g,d) :
\mathcal{Y}_K\ldrt \mathcal{Y}_{dKd^{-1}}
$$
où $\GL_n (F)$ agit à gauche et $\O_D^\times$ à droite. Cette action
est définie en posant 
$$
(g,d). ( G,\rho,\eta) = (G,\rho\circ g^{-1}, \eta \circ d^{-1})
$$
où 
\begin{eqnarray*}
{d^{-1}} : \O_D &\ldrt & \O_D \\
 d' &\longmapsto & d'd^{-1}
\end{eqnarray*}
et on utilise l'identification de $\GL_n (F)$ à $\End_{\O_D} ( \mathbb{G})_\Q^\times$.  
\end{defi}

\begin{lemm}
L'action de $\GL_n (F)\times \O_D^\times$ s'étend en une action de
$\GL_n (F) \times D^\times$. 
\end{lemm}
\dem
Pour $(G,\rho,\eta)\in\mathcal{Y}_K$ et $(g,d)\in \GL_n (F)\times  D^\times$ posons 
$
(g,d). ( G,\rho,\eta) = (G',\rho',\eta ')
$
où
\begin{itemize}
\item si $d\in \O_D, d\in \Pi^a \O_D$, alors 
$$
G'= G/G[\Pi^a]
$$ 
et si $\ph$ est l'isogénie $\ph: G\twoheadrightarrow G/G[\Pi^a]$ alors 
$$
\rho'= \ph \circ \rho \circ g^{-1}
$$
et $\eta'$ fait commuter le diagramme suivant 
$$
\xymatrix@C=14mm{
\O_D \ar[d]^{\bullet d} \ar[r]^{\eta}_\simeq & T_p (G^{rig})  \ar@{^(->}[d]^{\ph_*}
& \hspace{-15mm} [K] \\
\O_D \ar[r]^{\eta'}_\simeq & T_p ({G'}^{rig}) & \hspace{-15mm} [dKd^{-1}] 
}
$$
\item pour $d$ général on pose 
$$
(g,d). (G,\rho,\eta) = (\pi^b g , \pi^b d ). (G,\rho,\eta)
$$
avec $b>>0$.
\end{itemize}
On vérifie que cela définit bien une action pour laquelle $\forall
x\in F^\times$ $(x,x)\in \GL_n (F)\times D^\times$ agit trivialement.
\qed

\begin{rema}
Le fait que l'action de $\O_D^\times$ s'étende en une action  de
$D^\times$ au normalisé résulte en fait de ce que le sous-groupe
compact maximal $\O_D^\times$ est distingué dans $D^\times$. Cela est
faux en général pour un espace de Rapoport-Zink quelconque.
\end{rema}

\begin{defi}
On note $\mathcal{Y}_\infty$ le schéma formel $\pi$-adique $\GL_n (F)\times D^\times$-équivariant 
sur $\spf (\breve{\O})$ 
limite projective dans la catégorie des schémas formels $\pi$-adiques des $(\mathcal{Y}_K)_{K\subset \O_D^\times}$. 
\end{defi}

\begin{rema}
Dans \cite{Cellulaire} on a construit le schéma formel $\X_\infty$ du coté Lubin-Tate par recollement de cellules (des schéma formels) indexées par les sommets d'un immeuble de Bruhat-Tits. 
Le schéma formel $\mathcal{Y_\infty}$ admet une telle décomposition
cellulaire puisque c'est le cas de $\widehat{\Omega}$. Cependant pour $\mathcal{Y}_\infty$ les cellules sont indexées par les simplexes de l'immeuble. La raison en est que via l'isomorphisme entre les deux tours le domaine fondamental de Gross-Hopkins dans l'espace de Lubin-Tate n'est pas envoyé sur un simplexe de l'immeuble du coté Drinfeld mais sur un sous-ensemble plus petit (il s'agit d'un domaine fondamental pour toutes les correspondances de Hecke sphériques alors qu'un simplexe maximal est un domaine fondamental pour 
les correspondances de degré un multiple de $n$, cf. \cite{Rami}). Les deux décompositions cellulaires de $\X_\infty$ et $\mathcal{Y}_\infty$ ne se correspondent donc pas directement. C'est une des raisons pour lesquelles on devra éclater le schéma formel $\mathcal{Y}_\infty$ pour faire apparaître cet ensemble plus petit que le simplexe qui n'apparaît pas dans la fibre de spéciale de $\widehat{\Omega}$. 
\end{rema}

\section{Construction du morphisme $\widetilde{\X}_\infty\ldrt \widehat{\Omega}$}\label{decosdza}

Le but de ce chapitre est de construire un morphisme équivariant
d'un certain
éclaté de $\X_\infty$ vers le schéma formel de Deligne-Drinfeld
$\widehat{\Omega}$. 

\subsection{Définition de l'application de Hodge-Tate sur une cellule
  de $\X_\infty$} \label{chaverztu}

On reprend les notations de \cite{Cellulaire}. Soit $a=[\La,M]$ un
sommet de l'immeuble $\mathcal{I}$ de $(\GL_{n/F}\times
D^\times)/\Gm$. Soit 
$$
\mathbb{D}_{a,\infty}  = \underset{K\subset GL (\La)}{\limp} \mathbb{D}_{a,K}
$$
le schéma formel affine $\pi$-adique sur $\spf ( \breve{\O})$ construit dans
\cite{Cellulaire} (la cellule associée au sommet $a$). Afin d'alléger
les notations nous noterons momentanément 
\begin{eqnarray*}
\DD& :=& \DD_{a,GL (\La)} \\
\forall k\geq 1 \; \; \DD_k &:= &\DD_{a,Id+\pi^k \End (\La)} \\
\text{ et  } \;\DD_\infty &:=& \underset{k}{\limp} \DD_k
\end{eqnarray*}

{\it Fait admis  : On admettra que $\forall k\geq 1 \; \O_F/\pi^k
\O_F (1)$ est trivial sur $\DD_k$, c'est à dire $\O_{\DD_k}$ contient
les points de $\pi^k$-torsion d'un groupe de Lubin-Tate de hauteur $1$
(si $F=\Qp$ $\; \DD_k\ldrt \spf (\Zp[\zeta_{p^k}])$). 
L'auteur prévoit d'en donner la dmonstration  dans \cite{Periodes},
comme corollaire de la construction d'une application déterminant}
\\
 
Soit $H$ le $\O$-module $\pi$-divisible universel sur $\DD$. On note
$H^\vee$ son dual strict au sens de Faltings
(\cite{Faltings7}). D'après la remarque 8.2 de \cite{Cellulaire} pour
tout entier $k\geq 1$ le groupe $H$ est muni d'une structure de niveau
de Drinfeld de niveau $k$ sur $\DD_k$ 
$$
\eta : \pi^{-k}\La/\La \ldrt H[\pi^k] (\DD_k)
$$
Rappelons qu'étant donnée que $\DD_k$ est normal l'application fibre
générique induit une bijection entre ensembles finis 
$$
\pi_0 (\DD_k)\iso \pi_0 (\DD_k^{rig})
$$
 Si $\mathfrak{Z}$ est une composante connexe de $\DD_k$ alors
 $\eta$ induit une bijection 
$$
\eta : \pi^{-k}\La/\La \iso H[\pi^k] (\mathfrak{Z})
$$
puisque le schéma formel $\mathfrak{Z}$ étant normal et $H[\pi^k]$ un groupe
plat fini on a 
$$
H[\pi^k] ( \mathfrak{Z}) = H[\pi^k]^{rig}  ( \mathfrak{Z}^{rig}) 
$$
et que par définition d'une structure de niveau sur la fibre  fibre générique un
tel $\eta$ doit induire un isomorphisme en fibre générique sur chaque
composante connexe. 
De plus pour une telle composante connexe l'accouplement de dualité 
$$
H[\pi^k] (\mathfrak{Z}) \times H[\pi^k]^\vee (\mathfrak{Z}) \ldrt F/\O_F (1)
$$
est parfait puisque c'est le cas en fibre générique. 
 L'isomorphisme $\eta$ induit donc pour chaque composante connexe
 $\mathfrak{Z}$ de $\DD_k$ un isomorphisme 
$$
\eta^\vee : \La^*/\pi^k \La^* (1) \iso H[\pi^k]^\vee (\mathfrak{Z})
$$
Rappelons que l'inclusion $H[\pi^k]\hookrightarrow H$ induit un
isomorphisme
$$
\omega_H/\pi^k \omega_H \iso \omega_{H[\pi^k]}
$$
 Pour
tout $\DD$-schéma formel $\mathfrak{Z}$ il y a  une application de Hodge-Tate modulo $\pi^k$ 
\begin{eqnarray*}
H[\pi^k]^\vee  (\mathfrak{Z}) & \ldrt & \omega_{H[\pi^k]}\otimes_{\O_\DD}
\O_{\mathfrak{Z}} = \omega_H\otimes  \O_{\mathfrak{Z}}/\pi^k
\O_{\mathfrak{Z}} \\
x & \longmapsto & (x^\vee)^* \beta_{\LT}
\end{eqnarray*}
où si $x\in H[\pi^k]^\vee (\mathfrak{Z})$, $x$ correspond à un
morphisme strict 
$$
x^\vee : H[\pi^k] \times_{\DD}\mathfrak{Z} \ldrt \LT [\pi^k]
$$
et $(x^\vee)^*$ est le morphisme induit au niveau des formes
différentielles invariantes. 
\\

En particulier appliquant cela pour $\mathfrak{Z}$ variant parmi les composantes
connexes de $\DD_k$ et 
utilisant la
rigidification $\eta^\vee$, on obtient une application de Hodge-Tate
$$
\a_{H^\vee [\pi^k]} : \La^*\otimes \O_{\DD_k}/\pi^k \O_{\DD_k} (1)
\ldrt \omega_H \otimes \O_{\DD_k}/\pi^k \O_{\DD_k} 
$$
Ces différentes applications sont compatibles lorsque $k$ varie :
$$
\a_{H^\vee[\pi^{k+1}]}\text{ mod }\pi^k
\equiv \a_{H^\vee[\pi^k]} \otimes_{\O_{\DD_k}} \O_{\DD_{k+1}}
$$
 On en
déduit un morphisme de Hodge-Tate 
$$
\boxed{
\a_{H^\vee} : \La^* \otimes_{\O_F} \O_{\DD_\infty} (1) \ldrt
\omega_H\otimes \O_{\DD_\infty}}
$$
vérifiant 
$$
\forall k\;\; \a_{H^\vee}  \text{ mod } \pi^k \equiv \a_{H^\vee
  [\pi^k]}\otimes_{\O_{\DD_k}} \O_{\DD_\infty}
$$

\subsection{Rappels de quelques résultats de \cite{Rami} sur l'application de Hodge-Tate dans le cas d'un point}
\label{dkluiyar16}

\begin{prop}\label{rugytok}
Soit $K|F$ un corps valué complet pour une valuation $v$ à valeurs
dans $\R$ étendant celle de $F$. Soit $H$ un $\O$-module
$\pi$-divisible formel de dimension $1$ sur $\O_K$ et 
$$
\a_{H^\vee} : T_p (H^\vee)\ldrt \omega_H\otimes \O_{\widehat{\overline{K}}}
$$
l'application de Hodge-Tate de $H^\vee$.
Supposons que le polygone de Newton de la multiplication par $\pi$ associé à $H$ soit dans le domaine fondamental de Gross-Hopkins. 
 Alors
$$
\forall w\in T_p (H^\vee)\setminus \pi T_p (H^\vee)\;\;\; \a_{H^\vee} (w) \notin \pi^2 . \omega_H\otimes  \O_{\widehat{\overline{K}}}
$$
\end{prop}
\dem
Cela résulte de la formule donnée dans la démonstration du théorème
6.1 de  \cite{Rami} qui exprime la
valuation de $\a_{H^\vee} (w)$ en fonction du polygone de Newton $\mathcal{P}$ lorsque celui-ci est dans le ``simplexe fondamental'' : 
$$
v(\a_{H^\vee} (w)) = \frac{q}{q-1} \left (  1 -  (\mathcal{P} (q^{i-1})- \mathcal{P} (q^i))\right )
$$
\qed

Citons également la proposition suivante qui résulte immédiatement de \cite{Rami}.

\begin{prop}\label{ujifxd}
Reprenons les hypothèses de la proposition précédente. Soit $\Omega  (\widehat{\overline{K}}  ) \subset \mathbb{P} ( V_p (H^\vee)) ( \widehat{\overline{K}} )$ l'espace de Drinfeld et $\mathcal{I}$ l'immeuble de $\text{PGL} (V_p (H^\vee))$. Soit $\l : \Omega (\widehat{\overline{K}}) \ldrt |\mathcal{I}|$ la rétraction sur la réalisation géométrique de l'immeuble. Soit $E$ l'étoile du sommet $[T_p (H^\vee)]$ dans $\mathcal{I}$ et $|E|$ sa réalisation géométrique, 
$$
|E|=\bigcup_{\s\text{ simplexe de }\mathcal{I}\atop [T_p (H^\vee)]\in \s} |\s|
$$
Alors si $x\in  \Omega (\widehat{\overline{K}})$ est le point correspondant à $\a_{H^\vee}$
$$
\l (x) \in \text{Intérieur} ( |E|)
$$
\end{prop}

\subsection{Quelques rappels sur le schéma formel de
  Deligne-Drinfeld}\label{ramlqs}

\subsubsection{L'ouvert associé à un simplexe}

Soit $\s$ un simplexe de l'immeuble de $\text{PGL}_{n/F}$ représenté
par une suite de réseaux
$$
\pi\La_1\subsetneq \La_r\subsetneq \dots \subsetneq \La_1
$$
On note alors $\widehat{\Omega}_\s$ le schéma formel classifiant les
diagrammes 
$$
\xymatrix@C=5mm{
\La_{r+1} =\pi \La_1 \ar@{^(->}[r] \ar[d]^{u_1 \pi^{-1}} & \La_{r} \ar@{^(->}[r]\ar[d]^{u_r} & \dots
\ar@{^(->}[r] & \La_i \ar@{^(->}[r] \ar[d]^{u_i} & \dots \ar@{^(->}[r]
& \La_1\ar[d]^{u_1} \\
\L_1 \ar[r] \ar@(rd,ld)[rrrrr]_{\times \pi}
 & \L_r \ar[r] & \dots \ar[r] & \L_i \ar[r] & \dots \ar[r]
& \L_1
}
$$
sur un $\O_F$-schéma $S$ sur lequel $\pi$ est nilpotent 
où les $\L_i$, $1\leq i\leq r$, sont des fibrés en droites, les
applications $u_i$ sont $\O_F$-linéaires 
 et $\forall
i, \;1\leq i\leq r, \;\forall x\in S$
$$
\ker ( \La_i/\pi\La_i \ldrt \L_i\otimes k(x)) \subset \La_{i+1}/\pi\La_i
$$
On remarquera en particulier que $\forall i$ les sections  $u_i (\La_i)$ engendrent
le fibré $\L_i$. 
\\

Soit maintenant $\widehat{\mathbb{P}}(\La_1)$ le schéma formel complété $\pi$-adique
de l'espace projectif $\mathbb{P} (\La_1)$ sur $\spec (\O_F)$. Le
simplexe définit un drapeau de variétés linéaires dans la fibre
spéciale de $\widehat{\mathbb{P}}(\La_1)$
$$
\mathbb{P} (\La_1/\La_2) \subsetneq \mathbb{P}(\La_1/\La_3) \subsetneq
\dots \subsetneq \mathbb{P}(\La_1/\La_r) \subsetneq
\mathbb{P}(\La_1/\pi\La_1 ) =\widehat{\mathbb{P}}(\La_1)\times_{\spf
  (\O_F)} \spec (\Fq)
$$
Notons $\widehat{W}_\s$ l'éclatement formel admissible de toutes ces
variétés de la fibre spéciale. Si
$$
u:\La_1\otimes_{\O_F} \O_{\widehat{\mathbb{P}}(\La_1)}
\twoheadrightarrow \O_{\widehat{\mathbb{P}}(\La_1)} (1)
$$
désigne le morphisme universel sur $\widehat{\mathbb{P}}(\La_1)$ et si
$$
\forall 1\leq i\leq r\;\;\; \mathcal{I}_{\La_i} = u( \La_i\otimes
\O_{\widehat{\mathbb{P}}(\La_1)}) (-1) \subset \O_{\widehat{\mathbb{P}}(\La_1)}
$$
est l'idéal cohérent définissant $\mathbb{P} (
\La_1/\La_i)\hookrightarrow \widehat{\mathbb{P}}(\La_1)$ alors $W_\s$
est l'éclatement formel admissible de l'idéal
$$
\prod_{1\leq i\leq r} \mathcal{I}_{\La_i}
$$
Ainsi sur $\widehat{W}_\s$ l'image par $u$ des $\La_i$ devient
localement libre de rang $1$ et $\widehat{W}_\s$ classifie les
diagrammes 
$$
\xymatrix@C=5mm{
\La_{r+1} =\pi \La_1 \ar@{^(->}[r] \ar[d]^{u_1 \pi^{-1}} & \La_{r} \ar@{^(->}[r]\ar[d]^{u_r} & \dots
\ar@{^(->}[r] & \La_i \ar@{^(->}[r] \ar[d]^{u_i} & \dots \ar@{^(->}[r]
& \La_1\ar[d]^{u_1} \\
\L_1 \ar[r] \ar@(rd,ld)[rrrrr]_{\times \pi}
 & \L_r \ar[r] & \dots \ar[r] & \L_i \ar[r] & \dots \ar[r]
& \L_1}
$$
sur un $\O_F$-schéma $S$ sur lequel $\pi$ est nilpotent 
où les $\L_i$ sont des fibrés en droites et $\forall i\;\; u_i
(\La_i)$ engendre $\L_i$, ou encore de façon équivalente $\forall x\in
S$ l'application $\La_i/\pi\La_i \ldrt \L_i\otimes k(x)$ est
non-nulle. 

Il résulte de cette description que 
$$
\widehat{\Omega}_\s \subset \widehat{W}_\s
$$
est un sous-schéma formel ouvert. 

\begin{exem}
Lorsque $n=2$ $\;\widehat{W}_\s$ est le schéma formel stable dont la
fibre spéciale est l'union d'une droite projective $D$ sur $\Fq$ 
et de droites projectives $(L_x)_{x\in D (\Fq)}$ telles que  la droite $L_x$
intersecte $D$ en $x$.  Le schéma formel $\widehat{\Omega}_\s$ est
l'ouvert obtenu en retirant les points sur $\Fq$ des droites
projectives $(L_x)_x$, excepté le point $\{ x \} =L_x\cap D$. 
\end{exem}

\begin{rema}
Une autre définition de $\widehat{W}_\s$ consisterait à éclater
$\mathbb{P} ( \La_1)$ le long de $\mathbb{P} ( \La_1/\La_2)$ puis
éclater le transformé strict de $\mathbb{P}( \La_1/\La_3)$, puis le
transformé strict de $\mathbb{P} ( \La_1/\La_4)$ et ainsi de suite,
puis de prendre le complété $p$-adique du schéma obtenu. On peut vérifier
que ces deux définitions coïncident (utiliser le fait que les
sous-variétés de la fibre spéciale que l'on éclate sont définies pas
des suites régulières). 
\end{rema}

Si $\s'\subset \s$ est un sous-simplexe il y a alors
une immersion ouverte naturelle $\widehat{\Omega}_{\s'}
\subset \widehat{\Omega}_\s$ et le schéma formel de Deligne-Drinfeld
$\widehat{\Omega}$ est ainsi obtenu par recollement des
$\widehat{\Omega}_\s$ en utilisant les relations de faces données par
l'immeuble.

\subsubsection{L'ouvert associé à l'étoile d'un sommet} 

Fixons maintenant un réseau $\La$. Soit $E(\La)$ l'étoile de $\La$
dans l'immeuble c'est à dire 
$$
E(\La)= \bigcup_{\s\text{ simplexe} \atop [\La]\in \s} \s
$$
Intéressons-nous à
$$
\widehat{\Omega}_{E(\La)}= \bigcup_{\s\atop [\La]\in \s}
\widehat{\Omega}_\s \subset \widehat{\Omega}
$$
Soit $\widehat{W}_{E(\La)}$ le schéma formel suivant. Partons de
$\widehat{\mathbb{P}}(\La)$ est considérons la famille de
sous-variétés linéaires de sa fibre spéciale
$$
\{ \mathbb{P}(\La/\La')\;|\; \pi\La \subsetneq \La'\subsetneq \La \}
$$
Définissons $\widehat{W}_{E(\La)}$ comme l'éclaté formel admissible de
cette famille de sous-variétés. Plus précisément, si 
$$
u:\La\otimes \O_{\widehat{\mathbb{P}}(\La)} \twoheadrightarrow \O_{\widehat{\mathbb{P}}(\La)}(1)
$$
désigne le morphisme universel sur $\widehat{\mathbb{P}}(\La)$ alors
$\widehat{W}_{E(\La)}$ est l'éclaté formel admissible du produit
des idéaux $u(\La'\otimes  \O_{\widehat{\mathbb{P}}(\La)}) (-1)$ pour $ \pi\La
\subsetneq \La'\subsetneq \La$. 

Alors 
$$
\widehat{\Omega}_{E(\La)} \subset \widehat{W}_{E(\La)}
$$
est un sous-schéma formel ouvert.
Pour tout schéma formel $\pi$-adique sans $\pi$-torsion
$\mathfrak{Z}$ sur $\spf (\O_F)$ l'ensemble $\widehat{W}_{E(\La)}
(\mathfrak{Z})$ est (fonctoriellement en $\mathfrak{Z}$) l'ensemble des classes d'isomorphismes 
$$
\left [ u: \La\otimes \O_{\mathfrak{Z}}\twoheadrightarrow \L \right ]
\in \widehat{\mathbb{P}}(\La)(\mathfrak{Z})
$$
telles que 
$\forall \La',\; \pi\La \subsetneq \La' \subsetneq \La\;$
$u(\La'\otimes \O_{\mathfrak{Z}})$ soit localement libre de rang $1$ 
 (si $\mathcal{Z}$ a de la $\pi$-torsion c'est à dire n'est pas plat
 sur $\spf (\O_F)$ alors la description de $\widehat{W}_{E(\La)}
(\mathfrak{Z})$ est plus complexe). 

Soit maintenant le domaine analytique fermé
associé dans l'espace de Berkovich fibre générique 
$$
\widehat{\Omega}_{E(\La)}^{an} \subset  \widehat{W}_{E(\La)}^{an} =
\mathbb{P} (\La)^{an} = ( \mathbb{P}^{n-1})^{an}
$$
Si $ \l : \Omega \ldrt |\mathcal{I}|$ désigne
l'application vers la réalisation géométrique de
l'immeuble $\mathcal{I}$  alors 
$$
\widehat{\Omega}_{E(\La)}^{an} =\l^{-1} \left ( |E(\La)| \right )
$$

Nous aurons besoin du lemme suivant qui dit que sur l'ouvert
$\widehat{\Omega}_{E(\La)}$ de $\widehat{W}_{E(\La)}$ il est inutile
d'éclater de nouveaux réseaux. 

\begin{lemm}\label{klupitag}
Soit $\mathfrak{Z}$ un schéma formel $\pi$-adique sans $\pi$-torsion
sur $\spf (O)$. Soit $\L$ un fibré en droites sur $\mathfrak{Z}$ et 
$u:\La\otimes\O_{\mathfrak{Z}}\twoheadrightarrow \L$ un morphisme
donnant naissance à un morphisme $\mathfrak{Z}\ldrt
\widehat{\Omega}_{E(\La)}$. Alors, pour tout réseau $\La'\subset \La$
$\;\; u(\La'\otimes \O_{\mathfrak{Z}}) \subset \L$ est localement
libre de rang $1$.
\end{lemm}
\dem
Il suffit de le montrer pour $\mathfrak{Z} =\widehat{\Omega}_\s$ où
$\s$ est un simplexe maximal possédant $[\La]$ comme sommet. Soit donc
$\s$ un simplexe maximal associé à une chaîne $\La_{n+1}=\pi\La
\subsetneq \La_n\subsetneq \dots \subsetneq \La_2\subsetneq
\La_1=\La$.
Soit 
$$
\xymatrix@C=5mm{
\La_{n+1} =\pi \La_1 \ar@{^(->}[r] \ar[d]^{u_1 \pi^{-1}} & \La_{n} \ar@{^(->}[r]\ar[d]^{u_n} & \dots
\ar@{^(->}[r] & \La_i \ar@{^(->}[r] \ar[d]^{u_i} & \dots \ar@{^(->}[r]
& \La_1\ar[d]^{u_1} \\
\L_1 \ar[r] \ar@(rd,ld)[rrrrr]_{\times \pi}
 & \L_n \ar[r] & \dots \ar[r] & \L_i \ar[r] & \dots \ar[r]
& \L_1}
$$
le diagramme universel sur $\widehat{\Omega}_{\s}$. Deux simplexes
de l'immeuble sont contenus dans un même appartement. Il existe donc
une base $(e_1,\dots,e_n)$ de $\La$ telle que 
$$\forall i\geq 2\; \La_i=<\pi e_1,\dots,\pi e_{i-1},e_i,\dots,e_n>$$ et
des entiers $(a_1,\dots,a_n)\in\N^n$ tels que
$\La'=<\pi^{a_1}e_1,\dots,\pi^{a_n}e_n>$. 
\\
Pour tout $i$, par définition de $\widehat{\Omega}_\s$,
$$
\O_{\widehat{\Omega}_\s}.u(e_i)=\L_i
$$
Soit donc $\a=\inf\{ a_i\;|\; 1\leq i\leq n\}$ et $i_0=\inf \{ i\;|\;
a_i=\a\}$. Alors 
$$
u(\La'\otimes \O_{\widehat{\Omega}_\s}) =\pi^{\a} \L_{i_0}
$$
\qed

\subsubsection{Recollement des morphismes vers $\widehat{\Omega}$}

L'idée de cette sous-section repose sur l'analogie suivante. Soit $S$ un schéma et $Y$ un $S$-schéma séparé. Soit $X$ un schéma muni d'un ouvert $U \subset X$ schématiquement dense (par exemple $U$ dense et $X$ réduit). Soient deux $S$-morphismes $\xymatrix{X\ar@<0.6ex>[r] \ar@<-0.6ex>[r] & Y}$ coïncidant sur l'ouvert $U$. Alors ces deux morphismes sont égaux. 
L'analogie  sera faite avec $X$ un schéma formel $\pi$-adique sans $\pi$-torsion, $U$ ``l'ouvert $\pi\neq 0$'' et $Y=\widehat{\Omega}$. 
\\

Sur l'espace annelé $(\widehat{\Omega},\O_{\widehat{\Omega}}[\frac{1}{\pi}])$ vit un $\O_{\widehat{\Omega}}[\frac{1}{\pi}]$-module localement libre de rang $1$ $\L$ ainsi qu'un morphisme injectif
$$
u: F^n\hookrightarrow \L
$$
induisant un épimorphisme $\O_{\widehat{\Omega}}[\frac{1}{\pi}]^n\twoheadrightarrow \L$. 
Cet objet est construit de la façon suivante. Soit $\s$ un simplexe de l'immeuble associé à une chaîne de réseaux $\pi\La_1\subsetneq \La_r\subsetneq \dots \subsetneq \La_1$. Soit $\widehat{\Omega}_\s\subset \widehat{\Omega}$ l'ouvert correspondant. Alors le diagramme 
$$
\xymatrix@C=5mm{
\La_{r+1} =\pi \La_1 \ar@{^(->}[r] \ar[d]^{u_1 \pi^{-1}} & \La_{r} \ar@{^(->}[r]\ar[d]^{u_r} & \dots
\ar@{^(->}[r] & \La_i \ar@{^(->}[r] \ar[d]^{u_i} & \dots \ar@{^(->}[r]
& \La_1\ar[d]^{u_1} \\
\L_1 \ar[r] \ar@(rd,ld)[rrrrr]_{\times \pi}
 & \L_r \ar[r] & \dots \ar[r] & \L_i \ar[r] & \dots \ar[r]
& \L_1}
$$
sur $\widehat{\Omega}_\s$ devient après inversion de $\pi$ c'est à dire sur 
l'espace annelé $(\widehat{\Omega}_\s,\O_{\widehat{\Omega}_\s}[\frac{1}{\pi}])$ 
$$
\xymatrix@C=5mm{
\La_{r+1}\unpi =\pi \La_1[\frac{1}{\pi}] \ar@{=}[r] \ar[d]^{u_1 \pi^{-1}} & \La_{r}[\frac{1}{\pi}] \ar@{=}[r] \ar[d]^{u_r} & \dots
\ar@{=}[r] & \La_i [\frac{1}{\pi}] \ar@{=}[r] \ar[d]^{u_i} & \dots \ar@{=}[r] 
& \La_1[\frac{1}{\pi}] \ar[d]^{u_1}  \ar@{=}[r]& F^n\\
\L_1 [ \frac{1}{\pi}] \ar[r]^\sim  & \L_r [\frac{1}{\pi}] \ar[r]^\sim & \dots \ar[r]^\sim & \L_i[\frac{1}{\pi}] \ar[r]^\sim  & \dots \ar[r]^\sim
& \L_1[\frac{1}{\pi}] }
$$
et fournit donc l'objet voulu sur $\widehat{\Omega}_\s$. On vérifie
aussitôt que ces différents objets se recollent sur
$\widehat{\Omega}$. En fait le
$\O_{\widehat{\Omega}}[\frac{1}{\pi}]$-module $\L$ est libre de rang
$1$ isomorphe à $\O_{\widehat{\Omega}}[\frac{1}{\pi}]$ mais nous ne
fixons pas d'isomorphisme entre $\L$ et $\O_{\widehat{\Omega}}[\frac{1}{\pi}]$.
Lorsque $\La$ varie parmi les réseaux dans $F^n$ les
sous-$\O_{\widehat{\Omega}}$-modules
$\O_{\widehat{\Omega}}.u(\La)\subset \L$ forment une chaîne de fibrés
en droites 
$$
(\L_i)_{i\in \Z},\;\;\; \L_i \subset \L,\;\; \L_{i+1}\subset
\L_i\;\text{ et }\; \pi \L_i = \L_{i+n}
$$
Si $\underline{F}^n$ désigne le faisceau constant alors si 
$$
\forall i\in\Z\;\;\eta_i=u^{-1} (\L_i) \subset \underline{F}^n
$$
la chaîne de sous-faisceaux $(\eta_i)_{i\in\Z}$ vérifie
$\eta_{i+1}\subset \eta_i$ et $\pi \eta_i=\eta_{i+1}$ (ce sont les
faisceaux constructibles associés à $\widehat{\Omega}$ par Drinfeld,
cf. \cite{DrinfeldOmega}, \cite{BoutotCarayol}). Et si $\s$ est un simplexe associé à une chaîne périodique de
réseaux $(\La_i)_{i\in \Z}$ alors 
$$
\widehat{\Omega}_\s = \{x\in \widehat{\Omega}\;|\; \forall i
\;\eta_{i,x}\in \{ \La_j\;|\; j\in \Z \}\}
$$

\begin{rema}
Si $\mathfrak{Z}$ est un schéma formel quasi-séparé 
$\pi$-adique sans $\pi$-torsion $\O_\mathfrak{Z} [\frac{1}{\pi}]$ est le faisceau associé au préfaisceau $\mathcal{U}\mapsto \GG(\mathcal{U},\O_{\mathfrak{Z}})[\frac{1}{\pi}]$. Ainsi sur un ouvert quasicompact $\mathcal{U}$ on a $\GG( \mathcal{U},\O_\mathfrak{Z} [\frac{1}{\pi}])= \GG(\mathcal{U},\O_{\mathfrak{Z}})[\frac{1}{\pi}]$ mais cela est faux en général. Par exemple sur $\widehat{\Omega}$ les $n$ sections de $\L$ associées à $F^n \ldrt \L$ ne sont pas associées à des sections entières provenant d'un fibré en droites sur $\widehat{\Omega}$;
 lorsque l'on sort de tout ouvert quasicompact de $\widehat{\Omega}$ les puissances de $\pi$ dans les dénominateurs de ces $n$ sections tendent vers l'infini. 
\end{rema}

Sur la fibre générique rigide, via le morphisme de topos annelés $sp:(\widehat{\Omega}^{rig},\O_{\widehat{\Omega}^{rig}})\ldrt (\widehat{\Omega},\O_{\widehat{\Omega}} [\frac{1}{\pi}])$, $sp^*[F^n\hookrightarrow \L]$ induit le plongement $\widehat{\Omega}^{rig}\hookrightarrow (\mathbb{P}^n)^{rig}$.
\\

\begin{lemm}\label{ikdsez}
Soit $\mathfrak{Z}$ un schéma formel $\pi$-adique sans $\pi$-torsion sur $\spf (\O_F)$. Soient deux morphismes $\xymatrix{\mathfrak{Z} \ar@<0.6ex>[r] \ar@<-.6ex>[r] & \widehat{\Omega}}$ et $F^n\ldrt \L$, resp. $F^n\ldrt \L'$, les objets associés sur $(\mathfrak{Z},\O_{\mathfrak{Z}}[\frac{1}{\pi}])$. Ces deux morphismes sont égaux ssi il existe un isomorphisme $\a$ faisant commuter le diagramme suivant 
$$
\xymatrix@C=12mm@R=4mm{
 &  \L\ar[dd]^\a_{\simeq} \\
F^n\ar[ru] \ar[rd] \\
 & \L'
}
$$
\end{lemm}
\dem
Il résulte des considérations précédentes que le morphisme
$F^n\ldrt \L$ permet de reconstruire les chaînes de faisceaux
$(\L_i)_{i\in \Z}$ ainsi que
$(\eta_i)_{i\in \Z}$. Le résultat s'en
déduit facilement.
\qed

\subsection{Sur le conoyau de l'application de Hodge-Tate}\label{maishegndgj}

\begin{prop}\label{cocokf}
Soit l'application de Hodge-Tate tordue
$$
\a_{H^\vee}(-1) : \La^*\otimes \O_{\DD_\infty}  \ldrt \omega_H\otimes \O_{\DD_\infty} (-1)
$$
Alors
$$
\forall w\in \La^*\setminus \pi \La^*\;\; \pi^2\omega_H\otimes \O_{\DD_\infty} (-1) \subset \O_{\DD_\infty}. \a_{H^\vee} (w\otimes 1)
$$
De même 
$$
\forall k\geq 3\;\; \pi^2.\omega_H\otimes \O_{\DD_k}/\pi^k \O_{\DD_k} (-1) \subset \O_{\DD_k}/\pi^k \O_{\DD_k}. \a_{H^\vee[\pi^k]} (w\otimes 1)
$$
\end{prop}
\dem
Afin de ne pas alourdir les notations oublions la torsion à la Tate $F(-1)$ dans cette démonstration. 
Choisissons un générateur $t\in \GG ( \DD,\omega_H)$ de $\omega_H$ 
$$
\omega_H =\O_\DD .t
$$
Soit $f\in \GG ( \DD_3,\O_{\DD_3})$ tel que 
$$
\a_{H^\vee} ( w\otimes 1) \equiv \a_{H^\vee [\pi^3]} (w\otimes 1) \equiv ft \;\text{ mod }
\pi^3 \O_{\DD_3}
$$
Pour tout point $x\in \DD_3^{rig} ( \overline{\breve{F}})$, $x: \spf ( \O_{K})\ldrt \DD_3$
avec $K|\breve{F}$ de degré fini, 
$$
x^* \a_{H^\vee} (w\otimes 1) = \a_{(x^* H)^\vee } ( w\otimes 1) \equiv 
f(x).x^*t\;\text{ mod } \pi^3 \O_{\Cp}
$$
où $x^*H$ est un groupe $p$-divisible sur $\O_{K}$ et 
$\omega_{x^* H} = \O_{K}. x^*t$. 
D'après la proposition \ref{rugytok}
$$
 \a_{(x^* H)^\vee } ( w\otimes 1) \notin \pi^2 \O_{\Cp}. x^*t
$$
Donc
$$
\left | \frac{\pi^2}{f(x)} \right | \leq 1 
$$
La fonction rigide $\dpt{\frac{\pi^2}{f}\in \GG ( \DD_3^{rig},\O_{\DD_3}^{rig})}$ vérifie donc
$$
\left \| \frac{\pi^2}{f} \right \|_\infty \leq 1 \limpl \frac{\pi^2}{f} \in \GG ( \DD_3,\O_{\DD_3})
$$
puisque $\DD_3$ est normal (cf. l'appendice A de \cite{Cellulaire}). De cela on déduit que 
$$
\pi^2 t \in \O_{\DD_3}. \a_{H^\vee } (w\otimes 1) + \pi^3\O_{\DD_\infty}.t
$$
et donc, $\DD_\infty$ étant $\pi$-adique 
$$
\pi^2 t \in  \O_{\DD_\infty}. \a_{H^\vee} (w\otimes 1)
$$
De même pour tout $k\geq 3$
$$
\pi^2 t \in \O_{\DD_k}/\pi^k \O_{\DD_k}. \a_{H^\vee [\pi^k]} (w\otimes 1) + \pi^3\O_{\DD_k}/\pi^k \O_{\DD_k} .t
$$
et on conclut comme précédemment.
\qed

\subsection{\'Eclatement de la cellule et construction du morphisme de
  la cellule éclatée vers $\widehat{\Omega}$}

On identifie $(F^n)^*$ et $F^n$ via la base duale de la base canonique
de $F^n$. Rappelons que $\La\subset F^n$ et donc $\La^*$ est un réseau
de $F^n$.
\\

Pour tout réseau $\La'$ vérifiant $\pi\La\subsetneq \La'\subset \La^*$
et tout entier $k\geq 3$ soit l'idéal cohérent 
$$
\mathcal{I}_{\La',k}\subset \O_{\DD_k}
$$
tel que $\pi^k \O_{\DD_k} \subset \mathcal{I}_{\La',k}$ et 
$$
\a_{H^\vee[\pi^k]} \left ( \La'\otimes \O_{\DD_k}/\pi^k \O_{\DD_k}
\right ) = \omega_H\otimes \mathcal{I}_{\La',k}/\pi^k \O_{\DD_k} 
$$
D'après la proposition \ref{cocokf} 
$$
\pi^2 \O_{\DD_k} \subset \I_{\La',k}
$$
et si $k\geq 3$ et $\Pi_{k,3}: \DD_k\ldrt \DD_3$ alors
$$
\O_{\DD_k}.\Pi_{k,3}^{-1}( \mathcal{I}_{\La',3}) =  \mathcal{I}_{\La',k}
$$

\begin{defi}
Pour tout entier $k\geq 3$ on note $\widetilde{\DD}_k$ le normalisé de
l'éclatement formel admissible de $\DD_k$ relativement aux idéaux
$\mathcal{I}_{\La',k}$ où $\pi\La\subsetneq \La'\subset \La^*$. 
\end{defi}

Il résulte de la propriété énoncée précédemment,
$\O_{\DD_k}.\Pi_{k,3}^{-1}( \mathcal{I}_{\La',3}) =
\mathcal{I}_{\La',k}$, que $\widetilde{\DD}_k$ est le normalisé du transformé
strict de $\DD_k\ldrt \DD_3$ relativement à l'éclatement
$\widetilde{\DD}_3\ldrt \DD_3$ (cf. section \ref{sdguiec23} de
l'appendice) 
 et donc les morphismes de transition
$\widetilde{\DD}_l \ldrt \widetilde{\DD}_k$ pour $l\geq k\geq 3$ 
sont finis.

\begin{defi}
On note $\widetilde{\DD}_\infty = \underset{k\geq 3}{\limp}
\widetilde{\DD}_k$ dans la catégorie des schémas formels
$\pi$-adiques. 
\end{defi}

\begin{rema}
D'après les résultats de l'appendice on peut construire
$\widetilde{\DD}_\infty$ directement en niveau infini comme le
normalisé dans sa fibre générique de l'éclatement formel admissible des
idéaux $\O_{\widetilde{\DD}_\infty}.\Pi_{\infty,3}^{-1}
\mathcal{I}_{\La',3}$, cf. corollaire \ref{sdmlkgiob35} de
l'appendice. 
\end{rema}

Avec les notations de la section \ref{ramlqs} et d'après les rappels
de cette même section l'application de Hodge-Tate tordue 
$$
\a_{H^\vee}(-1) : \La^* \otimes \O_{\widetilde{\DD}_\infty} \ldrt
\omega_H\otimes \O_{\widetilde{\DD}_\infty} (-1)
$$
induit un morphisme 
$$
\widetilde{\DD}_\infty \ldrt \widehat{W}_{E(\La^*)} 
$$

\begin{prop}
Le morphisme $\widetilde{\DD}_\infty \ldrt \widehat{W}_{E(\La^*)} $ se
factorise par l'ouvert $\widehat{\Omega}_{E(\La^*)} \subset
\widehat{W}_{E(\La^*)}$ et définit donc un morphisme $\widetilde{\DD}_\infty \ldrt \widehat{\Omega}$.
\end{prop}
\dem
Il suffit de vérifier qu'au niveau des fibres spéciales le morphisme 
$$
\widetilde{\DD}_\infty\otimes \Fq \ldrt \widehat{W}_{E(\La^*)} \otimes \Fq
$$
se factorise via l'ouvert $\widehat{\Omega}_{E(\La^*)}\otimes \Fq$. 

Le schéma $\widehat{W}_{E(\La^*)}\otimes \Fq$ étant de présentation finie sur 
$\spec (\Fq)$ et $\widetilde{\DD}_\infty =\underset{k}{\limp} \widetilde{\DD}_k\otimes \Fq$, $\exists k\geq 3$ et une factorisation
$$
\xymatrix{
\widetilde{\DD}_\infty \otimes \Fq \ar[r]\ar[d] & \widehat{W}_{E(\La^*)}\otimes \Fq \\
\widetilde{\DD}_k\otimes \Fq \ar[ru]
}
$$
Il suffit alors de montrer que l'image du morphisme $\widetilde{\DD}_k\otimes \Fq \ldrt \widehat{W}_{E(\La^*)}\otimes \Fq$ est contenue dans l'ouvert $\widehat{\Omega}_{E(\La^*)}\otimes \Fq$. 

Soit $|\widetilde{\DD}_k^{an}|=|\DD_k^{an}|$
 l'espace analytique de Berkovich fibre générique de $\widetilde{\DD}_k$. D'après la proposition 2.4.4. page 36 de \cite{BerkSpectral} et puisque $\widetilde{\DD}_k$ est normal le morphisme de spécialisation
$$
sp: |\widetilde{\DD}_k^{an}|\ldrt |\widetilde{\DD_k}\otimes \Fq |
$$
est surjectif. Soit donc $x\in \widetilde{\DD}_k\otimes \Fq$ et $y\in \widetilde{\DD}_k^{an}$ tel que $sp (y)=x$. Soit $K|\breve{F}$ une extension valuée complète telle que $y$ provienne du point $z\in \widetilde{\DD}_k^{an} (K) = \widetilde{\DD}_k (\O_K)$, $z:\spf (\O_K)
\ldrt \widetilde{\DD}_k$. On peut de plus supposer que $K=\widehat{\overline{K}}$. Dès lors 
$\exists z'\in \widetilde{\DD}_\infty (K)$ tel que $z'\mapsto z$ via 
$\widetilde{\DD}_\infty \ldrt \widetilde{\DD}_k$. D'après la proposition \ref{ujifxd} l'image de $z'$ dans $\widehat{W}_{E(\La^*)} (\O_K)$ est contenue dans 
$\widehat{\Omega}_{E(\La^*)}(\O_K)=\widehat{\Omega}_{E(\La^*)}^{an}
(K) \subset  \widehat{W}_{E(\La^*)}^{an} (K)=\mathbb{P}^n (K)$. Il en est donc de même de l'image de $z$. Le morphisme de spécialisation s'inscrit dans un diagramme 
$$
\xymatrix{
|\widetilde{\DD}_k^{an} |\ar[r]\ar[d]^{sp} & |\widehat{W}_{E(\La^*)}^{an}| \ar[d]^{sp}
\\
|\widetilde{\DD}_k\otimes \Fq |\ar[r] & |\widehat{W}_{E(\La^*)}\otimes \Fq|
}
$$
on en déduit que l'image de $x$ dans $\widehat{W}_{E(\La^*)}\otimes \Fq$ est dans l'ouvert 
$\widehat{\Omega}_{E(\La^*)}\otimes \Fq$.
\qed

\subsection{Recollement des morphismes sur les cellules}

\subsubsection{Recollement des cellules éclatées}

Soit $a=[\La,M]$ un sommet de l'immeuble et $\forall k\geq 3\;\widetilde{\DD}_{a,k}$ la cellule éclatée définie dans les sections
précédentes en niveau $Id+\pi^k \End (\La)\subset \GL (\La)$. Soit
$a\drt a'$ une arrête de l'immeuble (cf. \cite{Cellulaire}) où $a'=[\La',M']$ avec
$$
\La\subsetneq \La'\subsetneq \pi^{-1}\La \;\text{ et }\;
M'=\Pi^{-[\La:\La']} M
$$
On a défini dans \cite{Cellulaire} un ouvert Zariski $\DD_{a\drt a',k}\subset
\DD_{a,k}$. On note $\widetilde{\DD}_{a\drt a',k}$  son image réciproque à
$\widetilde{\DD}_{a,k}$.

Rappelons que $H$ est le groupe $p$-divisible universel sur
$\DD_a$. Sur $\DD_{a\drt a'}$ il y a un sous-groupe plat fini
$C\subset H[\pi]\times_{\DD_a} \DD_{a\drt a',k}$ tel que $\eta$ induise un diagramme commutatif
$$
\xymatrix{
\pi^{-k}\La/\La \ar[r]^(.4)\eta & H[\pi^k]\times_{\DD_a} \DD_{a\drt a',k} \\
\La'/\La \ar[r]  \ar@{^(->}[u] & C  \ar@{^(->}[u] 
}
$$
et $\eta$ induise une rigidification en fibre générique $\eta^{rig}:
\underline{\La'}/\underline{\La} \iso C^{rig}$.  Le groupe $C$ est un
sous-groupe canonique ``généralisé''.
Si $H'=H/C$ alors $H'$ couplé à la déformation universelle sur
$\DD_{a\drt a',k}$ et $\eta$ fournissent un morphisme 
$$
\DD_{a\drt a',k} \ldrt \DD_{a'\drt a,k-1}
$$

\begin{prop}
Supposons $k\geq 4$.
Alors l'application de recollement s'étend à la cellule éclatée.
Plus précisément, 
il y a alors un morphisme $\a$ faisant commuter le
diagramme suivant
$$
\xymatrix{
\widetilde{\DD}_{a\drt a',k}\ar[d] \ar[r]^\a & \widetilde{\DD}_{a'\drt a
  ,k-1} \ar[d] \\
\DD_{a\drt a',k} \ar[r] & \DD_{a'\drt a ,k-1}
}
$$
et induisant le morphisme de recollement de \cite{Cellulaire} en fibre générique.
\end{prop}
\dem
Il s'agit de voir que par l'application de Hodge-Tate associé à
$H'[\pi^{k-1}]$ sur $\widetilde{\DD}_{a\drt a',k}$ l'image des
réseaux $\La''$ vérifiant $\pi \La'^* \subsetneq \La''\subset \La'^*$
est localement libre de rang $1$. Pour cela on va utiliser le lemme
\ref{klupitag}. Dans cette démonstration on oublie une fois de plus
les torsions à la Tate afin de ne pas alourdir les notations.

Pour $G$ un groupe plat fini muni d'une action stricte de $\O$ le
morphisme
$$
\a_{G^\vee}:G^\vee \ldrt \omega_G
$$
est fonctorielle en $G$. Le morphisme $H[\pi^{k-1}]\ldrt
H'[\pi^{k-1}]$ sur $\DD_{a\drt a',k}$ 
induit donc un
diagramme 
$$
\xymatrix@C=3.4cm{
\La'^*\otimes \O_{\DD_{a\drt a',k}}/\pi^{k-1} \O_{\DD_{a\drt a',k}}
\ar[d] 
\ar[r]^{\a_{H'^\vee [\pi^{k-1}]}} & \omega_{H'} \otimes \O_{\DD_{a\drt a',k}}/\pi^{k-1} \O_{\DD_{a\drt a',k}}
\ar[d] \\
\La^*\otimes \O_{\DD_{a\drt a',k}}/\pi^{k-1} \O_{\DD_{a\drt a',k}}
\ar[r]^{\a_{H^\vee [\pi^{k-1}]}} & \omega_H \otimes \O_{\DD_{a\drt a',k}}/\pi^{k-1} \O_{\DD_{a\drt a',k}}
}
$$
où la flèche verticale de gauche est induite par l'inclusion
$\La'^*\subset \La^*$.  
Fixons des relèvements $\beta$, resp. $\beta'$, de
$\a_{H^\vee[\pi^{k-1}]}$, resp. $\a_{H'^\vee[\pi^{k-1}]}$ fournissant 
un diagramme 
$$
\xymatrix@C=2.4cm{
\La'^*\otimes \O_{\DD_{a\drt a',k}}
\ar[d] 
\ar[r]^{\beta'} & \omega_{H'} \otimes \O_{\DD_{a\drt a',k}}
\ar[d] \\
\La^*\otimes \O_{\DD_{a\drt a',k}}
\ar[r]^{\beta} & \omega_H \otimes \O_{\DD_{a\drt a',k}}
}
$$
qui commute modulo $\pi^{k-1}$. 
\'Etant donné que $k\geq 4$ pour tout réseau $\La''$ vérifiant
$\pi\La\subsetneq \La'' \subsetneq \La^*$ puisque $\pi^2
\O_{\DD_{a,k}}\subset \mathcal{I}_{\La'',k}$ et d'après les propriétés
de compatibilité du système des idéaux $ \mathcal{I}_{\La'',k}$
lorsque $k$ varie on a 
$$
\beta ( \La''\otimes\O_{\DD_{a\drt a',k}}) = \omega_H\otimes
(\mathcal{I}_{\La'',k})_{|\DD_{a\drt a',k}}
$$
Le morphisme $\beta \otimes_{\O_{\DD_{a\drt
      a',k}}}\O_{\widetilde{\DD}_{a\drt a',k}}$ définit donc un
morphisme 
$$
\widetilde{\DD}_{a\drt a',k} \ldrt \widehat{W}_{E(\La^*)}
$$
Ce morphisme est congru modulo $\pi$ au morphisme défini précédemment
en niveau infini $\widetilde{\DD}_{a\drt a',\infty}\ldrt
\widehat{\Omega}_{E(\La^*)}$. Il se factorise donc en un morphisme 
$$
\widetilde{\DD}_{a\drt a',k} \ldrt\widehat{\Omega}_{E(\La^*)}
$$
Il résulte alors du lemme \ref{klupitag} que $\forall \La''$ vérifiant 
$\pi\La'^*\subsetneq \La''\subset \La'^*$ on a 
$$
\beta(\La''\otimes \O_{\widetilde{\DD}_{a\drt a',k}}) =
\omega_H\otimes \mathcal{J}
$$
où $\mathcal{J}$ est un idéal localement libre de rang $1$ dans
$\O_{\widetilde{\DD}_{a\drt a',k}}$. Fixons un tel $\La''$. 

De la même façon que précédemment pour $\beta$, $\beta'$ vérifie 
$$
\beta' (\La''\otimes  \O_{\widetilde{\DD}_{a\drt a',k}}) =\omega_{H'}\otimes \mathcal{J}'
$$
où $\mathcal{J'}$ est un idéal vérifiant $\pi^2 \O_{\widetilde{\DD}_{a\drt a',k}} \subset \mathcal{J}'$ (utiliser la proposition \ref{cocokf} appliquée à $\DD_{a',k-1}$
et l'application de Hodge-Tate de $H'^\vee [\pi^{k-1}]$, puis restreint à $\DD_{a'\drt a,k-1}$ et enfin tirée en arrière via l'application $\DD_{a\drt a',k}\ldrt\DD_{a'\drt a,k-1}$). 

La congruence de $\beta$ et $\beta'$ modulo $\pi^{k-1}$ implique que dans $\omega_H\otimes \O_{\widetilde{\DD}_{a\drt a',k}}$ on a 
$$
\omega_{H'}\otimes \mathcal{J}' + \pi^{k-1}. \omega_H\otimes \O_{\widetilde{\DD}_{a\drt a',k}}
 = \omega_H\otimes \mathcal{J} + \pi^{k-1}. \omega_H\O_{\widetilde{\DD}_{a\drt a',k}}
$$
Mais étant donné que $(\pi^2)\subset \mathcal{I}, (\pi^2)\subset \mathcal{I}'$, que 
$\pi\omega_H\subset \omega_{H'}$ (car $C$ est inclus dans les points de $\pi$-torsion de $H$) 
 et que $k\geq 4$ on en déduit
$$
\mathcal{J}=\mathcal{J}'
$$
De tout cela on
déduit que l'application
$$
\a_{H'^\vee[\pi^{k-1}]}\otimes\O_{\widetilde{\DD}_{a\drt
    a',k}}  : \La'^*\otimes \O_{\widetilde{\DD}_{a\drt
    a',k}} /\pi^{k-1}\O_{\widetilde{\DD}_{a\drt
    a',k}} \ldrt \omega_{H'}\otimes \O_{\widetilde{\DD}_{a\drt
    a',k}}/\pi^{k-1} \O_{\widetilde{\DD}_{a\drt
    a',k}}
$$
vérifie que $\forall \La''$ tel que $\pi\La'^* \subsetneq
\La''\subset \La'^*$ 
$$
\a_{H'^\vee[\pi^{k-1}]} (\La''\otimes \O_{\widetilde{\DD}_{a\drt
    a',k}} /\pi^{k-1} \O_{\widetilde{\DD}_{a\drt
    a',k}}   ) = \omega_H\otimes\mathcal{K}/\pi^{k-1} \O_{\widetilde{\DD}_{a\drt
    a',k}} 
$$
où $\mathcal{K}$ est un idéal localement libre de rang $1$ contenant $\pi^2 \O_{\widetilde{\DD}_{a\drt
    a',k}}$. D'où l'existence de la factorisation annoncée. 
\qed
\\

De plus, $\forall k\geq 5 $ les applications de recollement composées 
$$
\widetilde{\DD}_{a\drt a',k}\ldrt \widetilde{\DD}_{a'\drt a,k-1} \ldrt
\widetilde{\DD}_{a\drt a',k-2}
$$
sont les morphismes de changement de niveau. Elles induisent donc des
isomorphismes en niveau infini
$$
\xymatrix{
\widetilde{\DD}_{a\ldrt a',\infty} \ar[r]^\sim \ar@(rd,ld)[rr]^{Id}
 &
\widetilde{\DD}_{a'\ldrt a,\infty} \ar[r]^\sim & \widetilde{\DD}_{a\ldrt a',\infty}
}
$$
On vérifie aisément que l'action de $\GL_n(F)\times D^\times$ est
naturellement compatible aux éclatements : $\forall (g,d)\in
\GL_n(F)\times D^\times$ et tout entier $k\geq 3$ l'isomorphisme
naturel
$$
(g,d): \DD_{a,Id+\pi^k End(\La)}\iso \DD_{(g,d).a,Id+\pi^k End (g^{-1}\La)}
$$
s'étend en un isomorphisme 
$$
(g,d): \widetilde{\DD}_{a,Id+\pi^k End(\La)}\iso \widetilde{\DD}_{(g,d).a,Id+\pi^k End(g^{-1}\La)}
$$
et induit donc un isomorphisme en niveau infini en passant à la limite
projective. 

\begin{defi}
On note $\widetilde{\X}_\infty$ le schéma formel $\GL_n(F)\times
D^\times$-équivariant  défini par le diagramme de recollement 
$$
\xymatrix{
\coprod_{a\drt a'} \widetilde{\DD}_{a\drt a',\infty} \ar@<0.8ex>[r] \ar@<-0.8ex>[r] &
\coprod_{a}
\widetilde{\DD}_{a,\infty} \ar[r] & \widetilde{\X}_\infty
}
$$
où $\forall a\drt a'$ les applications de recollement sont 
$\;\widetilde{\DD}_{a\drt a',\infty} \hookrightarrow \widetilde{\DD}_{a,\infty}$ l'inclusion naturelle et $\widetilde{\DD}_{a\drt a',\infty} \iso \widetilde{\DD}_{a'\drt a,\infty}\hookrightarrow \widetilde{\DD}_{a',\infty}$ l'application de recollement précédente composée avec l'inclusion naturelle.
\end{defi}

\subsubsection{Recollement des morphismes vers $\widehat{\Omega}$}

\begin{prop}
Les morphismes construits précédemment
$\widetilde{\DD}_{a,\infty}\ldrt \widehat{\Omega}$ pour des sommets
$a$ de l'immeuble se recollent en un morphisme 
$$
\widetilde{\X}_\infty \ldrt \widehat{\Omega}
$$ 
L'action à droite de $\GL_n (F)$ sur $\widetilde{\X}_\infty$ 
et celle à gauche sur $\widehat{\Omega}$ se correspondent via
$g\mapsto \,^t g$. 
\end{prop}
\dem
Il suffit de vérifier que pour toute arête $a\drt a'$ le
diagramme suivant est commutatif
$$
\xymatrix@R=7mm{
 & \widetilde{\DD}_{a,\infty} \ar[rd]  \\
\widetilde{\DD}_{a\drt a',\infty} \ar@{^(->}[ru]\ar@{_(->}[rd]  & & \widehat{\Omega}
\\
& \widetilde{\DD}_{a',\infty} \ar[ru]
}
$$
Cela découle de la fonctorialité de l'application de Hodge-Tate 
rappelée précédemment appliquée à l'isogénie ``quotient par un
sous-groupe canonique'' sur le bord des cellules.
 Plus précisément, avec les notations de la
démonstration précédente il y a un diagramme commutatif
$$
\xymatrix@C=18mm{
\La'^*\otimes \O_{\widetilde{\DD}_{a\drt a'},\infty} \ar@{^(->}[d]
\ar[r]^(.6){\a_{H'^\vee}} & \omega_{H'}  \ar@{^(->}[d] \\
\La^*\otimes \O_{\widetilde{\DD}_{a\drt a',\infty}}
  \ar[r]^(.6){\a_{H^\vee}} & \omega_H
}
$$
où $H'$ est le quotient de $H\times_{\DD_a}\times
\widetilde{\DD}_{a\drt a',\infty}$ par le sous-groupe canonique défini par
$\eta (\La'/\La)$. La proposition est donc une conséquence du lemme \ref{ikdsez}
\qed

\begin{rema}
Le morphisme construit est invariant par l'action de $D^\times$,
puisque construit à partir de l'application de Hodge-Tate du groupe de
Lubin-Tate universel qui ne dépend pas de la rigidification modulo $p$ (la
déformation notée $\rho$).  
\end{rema}

\section{Construction du morphisme $\widetilde{\X}_\infty \ldrt
  \mathcal{Y}_\infty$}
\label{jvhmzohr3}

Le but de ce chapitre est de relever le morphisme $\GL_n(F)$-équivariant défini précédemment $\widetilde{\X}_\infty \ldrt \widehat{\Omega}$ en un morphisme $\GL_n (F)\times D^\times$-équivariant $\widetilde{\X}_\infty \ldrt \mathcal{Y}_\infty$. 

\subsection{ \'Etude des normalisés de $\widehat{\Omega}$ dans la tour de Drinfeld}\label{rrkcd34}

Soit $G$ le $\O_D$-module formel spécial universel sur $\widehat{\Omega}$. On note pour tout sous-groupe compact-ouvert $K\subset \O_D^\times$ $$\Omega_K\ldrt \widehat{\Omega}^{rig}$$ le revêtement étale fini classifiant les structures de niveau $K$ sur $G^{rig} = \underset{k\geq 1}{\limi} G[\pi^k]^{rig}$. On note $\widehat{\Omega}_K$ le normalisé de $\widehat{\Omega}$ dans $\Omega_K$. Rappelons qu'alors 
$$
\mathcal{Y}_K= \coprod_\Z \widehat{\Omega}_K
$$

Rappelons les deux faits suivants :
\begin{itemize}
\item
 Soit $\mathfrak{Z}$ un schéma formel localement de type fini sans $\pi$-torsion sur $\spf (\breve{\O})$. Supposons $\mathfrak{Z}$ normal. Alors l'application qui à un ouvert $\mathcal{U}$ de $\mathfrak{Z}$ associe l'ouvert admissible $\mathcal{U}^{rig}$ de $\mathfrak{Z}^{rig}$ induit une bijection entre les parties ouvertes/fermées de $\mathfrak{Z}$ et celles de $\mathfrak{Z}^{rig}$.
\item Soit $X\xrig{\;\ph\;}Y$ un morphisme étale fini d'espaces rigides et $\xymatrix{Y\ar@<0.5ex>[r]^{s_1} \ar@<-0.5ex>[r]_{s_2} & X}$ deux sections de $\ph$. Alors l'espace rigide $\{s_1\neq s_2\}$ est un ouvert/fermé de $Y$.
\end{itemize}

\begin{lemm}
Soit $m\in \N^*$. Soit $\a\in G[\Pi^m]^{rig} ( G[\Pi^m]^{rig})$ la section identité et $e$ la section neutre. soit $U$ l'ouvert/fermé de $G[\pi^m]^{rig}$ où $\Pi^{m-1} (\a)\neq e$ c'est à dire $U= (\Pi^{m-1})^{-1} ( G[\Pi]^{rig}\setminus \{ e \})$ avec $\Pi^{m-1} : G[\Pi^m]^{rig}\ldrt G[\Pi]^{rig}$. Alors 
$\Omega_{1+\Pi^m\O_D} = U$.  
\end{lemm}
\dem 
L'espace $\Omega_{1+\Pi^m\O_D}$ est le foncteur au dessus de $\Omega$
défini par 
$$\forall Y\ldrt \Omega \;\;\;\Omega_{1+\Pi^m\O_D} (Y) =
\{\eta: \underline{\Pi^{-m}\O_D/\O_D} \iso
G[\Pi^m]^{rig}\times_{\Omega} Y \}$$
où $\eta$ est $\O_D$-linéaire et $ \underline{\Pi^{-m}\O_D/\O_D}$
désigne 
le groupe étale constant.  Se donner un morphisme  $\eta$ $\O_D$-équivariant
comme ci-dessus sur $Y$ est équivalent à
se donner la section $\eta (1)\in G[\Pi^{m}]^{rig} (Y)$. La condition
que $\eta$ soit un isomorphisme est équivalente à ce que pour toute
composante connexe  $W$ de $Y$ $\Pi^{m-1}
(\eta(1))_{|W}\neq 0$ soit encore que le morphisme $Y\ldrt
G[\Pi^m]^{rig}$ se factorise par l'ouvert/fermé $U$.
\qed

\begin{coro}
Soit $\widetilde{G[\Pi^m]}$ le normalisé de $G[\Pi^m]$ dans $G[\Pi^m]^{rig}$ (un $\breve{\O}$-schéma formel localement de type fini). Alors $\widehat{\Omega}_{1+\Pi^m\O_D}$ est l'ouvert/fermé de $\widetilde{G[\Pi^m]}$ induisant l'ouvert fermé $U$ du lemme précédent en fibre générique.\end{coro}

\begin{coro}
Soit $V=\spf (A)\subset \widehat{\Omega}$ un ouvert affine et $G[\Pi^m]_{|V}=\spf ( B_m)\ldrt \spf (A)$. Soit $W\subset \spec (B_m[\frac{1}{\pi}])$ l'ouvert fermé $\{\Pi^{m-1}\neq e \}$ et  $\widetilde{B_m}$ le normalisé de $B_m$. Alors si $\spec (A_m)$ est l'adhérence schématique de $W$ dans $\spec ( \widetilde{B_m})$ le diagramme suivant est cartésien
$$
\xymatrix{
\spf (A_m) \ar@{^(->}[r] \ar[d] & \widehat{\Omega}_{1+\Pi^m\O_D} \ar[d] \\
\spf (A) \ar@{^(->}[r] & \widehat{\Omega} 
}
$$
\end{coro}

\begin{exem}
Pour illustrer les constructions précédentes considérons l'espace de module obtenu en mettant des points de torsion sur $\mu_{p^\infty}$ sur $\spec (\Zp)$ (ce qui correspond au cas $n=1$ et $F=\Qp$). On a $\mu_{p^n}= \spec (\Zp [T]/(T^{p^n}-1))$.
De plus $\Zp [T]/(T^{p^n}-1)=\Zp[T]/(\prod_{i=0}^n \Phi_{p^i} (T))$ où $\Phi_k$ désigne le $k$-ième polynôme cyclotomique. Alors le normalisé de $\mu_{p^n}$ s'écrit 
$$
\widetilde{\mu_{p^n}} = \coprod_{i=0}^n \spec (\underbrace{\Zp[\zeta_{p^i}]}_{\Zp[T]/(\Phi_{p^i})})
$$
L'ouvert $T^{p^{n-1}}\neq 1$ est $\spec (\Zp [\zeta_{p^n}])$ car $\Zp[T]/(\Phi_{p^i}) \left [\frac{1}{T^{p^{n-1}}-1}\right ]=0$
 si $i<n$ puisqu'alors $\Phi_{p^i} | T^{p^{n-1}}-1$.
\end{exem}

\begin{prop}\label{truyyhg}
Soit $\mathfrak{Z}$ un $\breve{\O}$-schéma formel $\pi$-adique sans
$\pi$-torsion tel que $\O_{\mathfrak{Z}}$ soit intégralement fermé
dans $\O_{\mathfrak{Z}}\unpi$ (cela ne signifie rien d'autre que
pour tout ouvert affine $\mathcal{U}$ de
$\mathfrak{Z}$ l'anneau $\GG(\mathcal{U},\O_{\mathcal{Z}})$ est
intégralement clos dans
$\GG(\mathcal{U},\O_\mathfrak{Z})[\frac{1}{\pi}]$, cf. 
section \ref{krutofj248} de l'appendice). Se donner un
morphisme 
$$
\mathfrak{Z}\ldrt \widehat{\Omega}_{1+\Pi^m\O_D}
$$
est équivalent à se donner un triplet $(G',\rho',s)$ à isomorphisme près
où
\begin{itemize}
\item
 $(G',\rho')\in
\widehat{\Omega} (\mathfrak{Z})$ 
\item $s\in G'[\Pi^m] (\mathfrak{Z})$ est
une section telle qu'il existe un recouvrement affine
$(\spf(A_i))_i$ de $\mathfrak{Z}$ tel que 
si $G'_i$ désigne le groupe $p$-divisible sur $\spec (A_i)$ associé au
groupe $p$-divisible $G'\times_{\mathfrak{Z}} \spf (A_i)$ alors 
sur $\spec
(A_i[\frac{1}{\pi}])$ la section $s$ induit un isomorphisme 
$$
\underline{\Pi^{-m}\O_D/\O_D} \iso  G'_i[\Pi^m]\times_{\spec (A_i)} \spec
(A_i [\frac{1}{\pi}])
$$  
\end{itemize}
\end{prop}
\dem
Rappelons que $G$ désigne le $\O_D$-module formel universel sur
$\widehat{\Omega}$. Le triplet $(G',\rho')$ fournit un morphisme 
$$
\ph_1: \mathfrak{Z}\ldrt \widehat{\Omega}
$$
par lequel l'image réciproque de $G$ est isomorphe à $G'$. La section
$s$ relève ce morphisme en un morphisme $\ph_2$ 
$$
\xymatrix@R=7mm{
& G[\Pi^m] \ar[d] \\
\mathfrak{Z}\ar[r]^{\ph_1} \ar[ru]^{\ph_2} &  \widehat{\Omega}
}
$$
Celui-ci s'étend de façon unique en un morphisme $\ph_3$
$$
\xymatrix@R=7mm{
 & \widetilde{G[\Pi^m]} \ar[d] \\
\mathfrak{Z} \ar[r]^{\ph_2} \ar[ru]^{\ph_3} & G[\Pi^m]
}
$$
En effet, $\widetilde{G[\Pi^m]} = \spf ( \mathcal{A})$ où
$\mathcal{A}$ est une $\O_{G[\Pi^m]}$-algèbre cohérente vérifiant 
$\mathcal{A}\subset \O_{G[\Pi^m]}[\frac{1}{\pi}]$ et toute section de
$\mathcal{A}$ sur un ouvert quasicompact est entière sur
$\O_{G[\Pi^m]}$. Il y a donc d'après l''hypothèse faite sur
$\mathfrak{Z}$ 
une unique factorisation
$$
\xymatrix{
\ph_2^{-1} \O_{G[\Pi^m]} \ar[r]^(.6){\ph_2^*} \ar@{^(->}[d] & \O_\mathfrak{Z} \\
\ph_2^{-1} \mathcal{A} \ar[ru]
}
$$
qui fournit le morphisme cherché $\mathfrak{Z}\ldrt \spf
(\mathcal{A})$. 
\\
On vérifie alors localement sur $\mathfrak{Z}$ que ce morphisme se
factorise à travers l'ouvert/fermé de $\widetilde{G[\Pi^m]}$ égal à $\widehat{\Omega}_{1+\Pi^m\O_D}$.
\qed

\subsection{Définition modulaire de $\mathcal{Y}_\infty$}

\begin{lemm}[Pureté du $\pi_0$]\label{purifitazm}
Soit $R$ une $\breve{\O}$-algèbre $\pi$-adique sans $\pi$-torsion telle
que $R$ soit intégralement fermé dans $R[\frac{1}{\pi}]$. Alors
$$
\spec (R/\pi R) \text{ connexe }\lssi \spec (R) \text{ connexe } \lssi
\spec (R[\frac{1}{\pi}] ) \text{ connexe }
$$ 
\end{lemm}

Rappelons également que pour $\mathfrak{Z}$ un
$\breve{\O}$-schéma formel localement de type fini sans $\pi$-torsion
et normal on a $\mathfrak{Z}$ connexe $\lssi$ $\mathfrak{Z}^{rig}$
l'est.

\begin{prop}\label{crutozp}
Soit $\mathfrak{Z}$ un schéma formel $\pi$-adique
sans $\pi$-torsion sur $\spf (\breve{\O})$ tel que $\O_{\mathfrak{Z}}$
soit intégralement fermé dans $\O_{\mathfrak{Z}}\unpi$. 
Supposons que pour tout ouvert quasicompact $\mathcal{U}$ de
$\mathfrak{Z}$
$\;\pi_0 (\mathcal{U})$ est fini, c'est à dire les composantes
connexes de $\mathcal{U}$ sont ouvertes. 
 Il y a alors une bijection
$D^\times$-équivariante entre l'ensemble des morphismes 
$\mathfrak{Z}\ldrt \mathcal{Y}_\infty$ et
les triplets $(G',\rho',\eta)$ à isomorphisme près 
 où $(G',\rho')\in \widehat{\Omega} (\mathfrak{Z})$
et $\eta \in 
\GG ( \mathfrak{Z},\underline{\Hom} ( F/\O_F,G')[\frac{1}{\pi}])$
est tel que pour tout ouvert $\mathcal{U}$ de $\mathfrak{Z}$ $\;\;
\eta_{|\mathcal{U}}\neq 0$; l'action
de $D^\times$ se faisant sur les triplets via
$\forall d\in D^\times \; d.(G',\rho',\eta)=(G',\rho',\eta  \circ
(\bullet d^{-1}))$ où l'on note $\bullet d^{-1}$ la quasi-isogénie
$\O_D$-équivariante de
$D/\O_D$ dans lui même 
définie par $d'\O_D\mapsto d'd^{-1}\O_D$. 
\end{prop}
\dem 
Construisons un morphisme $\mathfrak{Z}\ldrt \mathcal{Y}_\infty$ à
partir d'un triplet $(G',\rho',\eta)$. 

Soit $\spf (R)\subset \mathfrak{Z}$ un ouvert affine connexe. Notons 
$G''$ le groupe $p$-divisible sur $\spec (R)$ associé au groupe
$p$-divisible 
$G'\times_{\mathcal{Z}}\spf (R)$ sur $\spf (R)$. L'élément $\eta$ induit un élément
$$
\a \in \Hom_{\O_D} ( D/\O_D, G'')\unpi\setminus \{ 0\} =
\Hom_{\O_D} (D/\O_D, G''\times_{\spec (R)} \spec
(R\unpi))\unpi\setminus \{ 0 \}
$$
l'égalité résultant de ce que $R$ est intégralement fermé dans
$R\unpi$. Le groupe  $\Hom_{\O_D} ( D/\O_D, G''\otimes_{R}
R [\frac{1}{\pi}])[\frac{1}{\pi}]$ est muni d'une structure de
$D$-module 
 via l'action de $D$ à droite sur $D/\O_D$. D'après le lemme
 \ref{purifitazm} $\spec (R\unpi)$ est connexe et donc ce $D$-module
 est de rang $1$. 

Soit $N\in \Z$ tel que 
$$
\Pi^N \a\in \Hom_{\O_D} ( D/\O_D, G'')\setminus
\Pi.\Hom_{\O_D} ( D/\O_D, G'')
$$
Toujours parce que $\spec (R [\frac{1}{\pi}])$ est
connexe $\Pi^N\a$ induit un isomorphisme 
$$
\Pi^N\a : D/\O_D \iso G\times_{\spec (R)}\spec (
R [\frac{1}{\pi}])
$$
Donc d'après la proposition \ref{truyyhg} le système compatible des
$$
((\Pi^N\a) (\bar{\Pi}^{-m}))_{m\geq 1}\in (G[\Pi^m]
(R))_{m\geq 1}
$$
induit un système compatible de morphismes 
$$
\xymatrix@R=5mm{
 &  \ar@{..>}[d] \\
 & \mathcal{Y}_{1+\Pi^{m+1}\O_D} \ar[d] \\
\spf (R) \ar[r] \ar[ru]  & \mathcal{Y}_{1+\Pi^{m}\O_D} \ar@{..>}[d] \\
 & 
}
$$
D'où, puisque $\mathcal{Y}_\infty = \underset{m}{\limp}
\mathcal{Y}_{1+\Pi^m\O_D}$ dans la catégorie des schémas formels
$\pi$-adiques, un morphisme $\ph: \spf ( R)\ldrt
\mathcal{Y}_\infty$.  On associe  alors à $(G',\rho',x)$ le morphisme
composé 
$$
\spf (R) \xrig{\;\ph\;} \mathcal{Y}_{\infty} \xrig{\;
  \Pi^{-N}\;} \mathcal{Y}_\infty
$$
Soient maintenant $\spf (R')$ et $\spf (R'')\neq \emptyset$ 
deux autres ouverts affines connexes
 tels que $\spf (R'') \subset \spf (R)\cap \spf (R')$. Il y a un
 diagramme de morphismes
$$
\xymatrix@R=6mm@C=8mm{
& \spec (R) \\
\spec (R'') \ar[ru]\ar[rd] \\
& \spec (R')
}
$$
Le groupe $p$-divisible sur $ \spec ( R')$ associé au groupe
$p$-divisible $G'\times_{\mathfrak{Z}}\spf (R'')$ est égal à
$G''\times_{\spec (R)} \spec (R'')$. 
On en déduit que l'entier $N\in \Z$ associé précédemment à $R$ est
le même pour $R,R'$ et $R''$. On conclu aisément que les différents
morphismes sur les ouverts affines connexes se recollent en un
morphisme
$$
\mathfrak{Z}\ldrt \mathcal{Y}_\infty
$$
Le fait que cela donne une bijection entre $\mathcal{Y}_\infty
(\mathfrak{Z})$ est les classes d'isomorphismes de triples
$(G',\rho',\eta)$ est laissé au lecteur.
\qed

\subsection{Sur la suite de Hodge-Tate en niveau infini}

Reprenons les notations de la section \ref{decosdza}. Soit $a=[\La,M]$ un sommet de l'immeuble paramétrant les cellules de $\X_\infty$. Soit $\DD_{a,\infty} = \underset{k}{\limp} \DD_{a,k}$ la cellule associée en niveau infini où $\DD_{a,k} := \DD_{a,Id+\pi^k End (\La)}$. 
Dans la section  \ref{decosdza} on a défini une application de Hodge-Tate
$$
\a_{H^\vee}:\La^*\otimes \O_{\DD_{a,\infty}} (1) \ldrt \omega_{H}\otimes \O_{\DD_{a,\infty}}
$$
On définit de même 
$$
\a_H : \La\otimes \O_{\DD_{a,\infty}} \ldrt \omega_{H^\vee}\otimes \O_{\DD_{a,\infty}}
$$
(le module $\omega_{H^\vee}$ est celui noté $\widetilde{\omega}_{H^D}$
dans l'appendice B de \cite{Cellulaire}). 

\begin{prop}\label{exhodateg}
Dans la suite de Hodge-Tate 
$$
\omega_{H^\vee}^*\otimes \O_{\DD_{a,\infty}} (1)\xrig{\; \,^t\a_{H}(1) \;} 
\La^*\otimes \O_{\DD_{a,\infty}} (1) \xrig{\; \a_{H^\vee}\;} \omega_H \otimes \O_{\DD_{a,\infty}}
$$
on a $\a_{H^\vee}\circ \,^t \a_{H} (1) =0$. 
\end{prop}
\dem
On utilisera le résultat suivant (cf. \cite{Points}). Si $H_0$ est un $\O$-module $\pi$-divisible sur $\O_K$
où $K|\breve{F}$ est une extension de degré fini alors dans la suite de Hodge-Tate usuelle
$$
\omega_{H_0}^*\otimes \O_{\widehat{\overline{K}}} (1) \xrig{\;\,^t\a_{H_0^\vee} (1)\;} 
T_p (H_0)\otimes \O_{\widehat{\overline{K}}}  \xrig{\; \a_{H_0}\;} \omega_{H_0^\vee} \otimes \O_{\widehat{\overline{K}}}
$$
on a $\a_{H_0}\circ \,^t \a_{H_0^\vee}(1)=0$. 

Revenons à l'énoncé. Il suffit de montrer que 
$$
\forall k\geq 1 \;\; \a_{H^\vee}\circ \,^t\a_{H}(1) \equiv 0 \text{ mod }\pi^k
$$
Mais $\forall k\geq 1$ la suite de Hodge-Tate modulo $\pi^k$ provient du niveau fini $\DD_{a,k}$
pour le groupe plat fini $H^\vee [\pi^k]$ 
\begin{eqnarray*} &&
[\omega_{H^\vee}^*\otimes \O_{\DD_{a,\infty}} (1)\xrig{\; \,^t\a_{H}(1) \;} 
\La^*\otimes \O_{\DD_{a,\infty}} (1) \xrig{\; \a_{H^\vee}\;} \omega_H \otimes \O_{\DD_{a,\infty}}
] \text{ mod } \pi^k \\
&=& 
[ \omega_{H^\vee}^*\otimes \O_{\DD_{a,k}}/\pi^k (1)\xrig{\; \,^t\a_{H[\pi^k]}(1) \;} 
\La^*\otimes \O_{\DD_{a,k}}/\pi^k   (1)
 \xrig{\; \a_{H^\vee[\pi^k]}\;} \omega_H \otimes \O_{\DD_{a,k}}/\pi^k  ]
\otimes \O_{\DD_{a,\infty}}/\pi^k \O_{\DD_{a,\infty}}
\end{eqnarray*}
où pour abréger on a noté $\O_{\DD_{a,\infty}}/\pi^k := \O_{\DD_{a,\infty}}/\pi^k\O_{\DD_{a,\infty}}$. De plus d'après le rappel précédent $\forall x\in \DD_{a,k}^{rig} ( \overline{\breve{F}})$, 
$x:\spf (\O_K)\ldrt \DD_k$ où $K|\breve{F}$ est de degré fini, 
$x^* \,^t\a_{H[\pi^k]}(1)$ et $x^* \a_{H^\vee[\pi^k]}$ sont congrus modulo $\pi^k$, après extension des scalaires de $\O_K$ à  $\O_{\widehat{\overline{K}}}$, à 
$\,^t\a_{x^*H}(1)$ et $\a_{x^* H^\vee}$ où $x^*H$ est un $\O$-module $\pi$-divisible sur $\O_K$.
Donc d'après le rappel du début de la démonstration 
$$
x^* ( \a_{H^\vee[\pi^k]}\circ \,^t\a_{H[\pi^k]}(1) ) =0 
$$
Soit $t$ un générateur de $\omega_H$, $\omega_H= \O_{\DD_{a}}.t$. 
Cela étant vrai pour tout $x$ on en déduit que pour $v\in \omega_{H^\vee}\otimes \O_{\DD_{a,k}} (1)$
  si $f\in \O_{\DD_{a,k}}$ est tel que $ (\a_{H^\vee[\pi^k]}\circ \,^t\a_{H[\pi^k]}(1)) (v) \equiv f.t$ alors 
$$
\forall x\in \DD_k^{rig} (\overline{\breve{F}})\;\;\;\left \| \frac{f}{\pi^k} \right \|_\infty \leq 1
$$
et que donc $\DD_k$ étant normal on a 
$$
f\in \pi^k\O_{\DD_k}
$$
\qed

\subsection{Construction d'éléments dans le module de Tate du
  $\O_D$-module formel spécial tiré en arrière sur $\widetilde{\X}_\infty$}\label{mvhzaf938}

On a construit dans le chapitre précédent un morphisme $\GL_n
(F)$-équivariant 
$$
\widetilde{\X}_\infty \ldrt \widehat{\Omega}
$$
à partir de l'application de Hodge-Tate de $H^\vee$ où $H$ est un
groupe $p$-divisible universel sur une cellule de
$\widetilde{\X}_\infty$.

Soit $G$ le $\O_D$-module formel spécial universel sur
$\widehat{\Omega}$ muni de sa rigidification, une quasi-isogénie de
hauteur $0$
$$
\rho_G : \mathbb{G}\times_{\Fqb} (\widehat{\Omega}\text{ mod }\pi)
\ldrt
G\times_{\widehat{\Omega}}(\widehat{\Omega}\text{ mod }\pi)
$$
Soit $a=[\La,M]$ un sommet de l'immeuble associé à $\X_\infty$ et 
$$
f_a:\widetilde{\DD}_{a,\infty} \ldrt \widehat{\Omega}
$$
le morphisme défini précédemment. Nous considérerons le $\O_D$-module
formel spécial $f_a^* G$ sur $\widetilde{\DD}_{a,\infty}$. 

Soit $H$ le $\O$-module $\pi$-divisible universel sur $\DD_a$ muni de
sa rigidification 
$$
\rho_H : \Hb\times_{\Fqb} (\DD_a \text{ mod }\pi) \ldrt
H\times_{\DD_a} (\DD_a \text{ mod }\pi)
$$
Le rigidification $\eta$ induit pour tout $k\geq 1$ 
$$
\eta : \pi^{-k}\La/\La \ldrt H\times_{\DD_a} \DD_{a,k}
$$
Elle fournit  donc un élément de  
$$
\Hom ( \La \otimes F/\La, H\times_{\DD_a}
\widetilde{\DD}_{a,\infty} )
$$
où $\La\otimes F/\La$ désigne le groupe $\pi$-divisible étale constant
$\underset{i}{\limi}\underline{\pi^{-i}\La/\La}$. 
Rappelons que $\La\subset F^n$ qui induit une quasi-isogénie 
$(F/\O_F)^n\ldrt \La \otimes F/\La$.
 Il y a donc une égalité
$$
\Hom ( \La [\frac{1}{\pi}]/\La, H\times_{\DD_a}
\widetilde{\DD}_{a,\infty} ) [\frac{1}{\pi}] = \Hom ( (F/\O_F)^n,  H\times_{\DD_a}
\widetilde{\DD}_{a,\infty} ) [\frac{1}{\pi}]
$$
et fournit donc $n$-éléments
$$
\zeta_1,\dots,\zeta_n\in \Hom ( F/\O_F,  H\times_{\DD_a}
\widetilde{\DD}_{a,\infty} ) [\frac{1}{\pi}]
$$
Considérons le composé 
$$
F/\O_F \xrig{(\zeta_1,\dots,\zeta_n)} H^n\times
(\widetilde{\DD}_{a,\infty}\text{ mod }\pi) \xrig{\; \rho_H^{-1}\;}
\Hb^n\times (\widetilde{\DD}_{a,\infty}\text{ mod }\pi )
\xrig{\;\Delta\;} \Gb\times (\widetilde{\DD}_{a,\infty}\text{ mod }\pi
)
\xrig{f_a^* \rho_G} f_a^* G\text{ mod }\pi
$$
qui définit un élément 
$$
\chi \in \Hom \left ( F/\O_F, f_a^* G \times_{\widetilde{\DD}_{a,\infty}}
({\widetilde{\DD}_{a,\infty}} \text{ mod }\pi ) \right  ) [\frac{1}{\pi}]
$$
L'application de réduction modulo $\pi$ induit une injection
$$
\Hom \left ( F/\O_F, f_a^* G   ) \right  )
[\frac{1}{\pi}]
\hookrightarrow
\Hom \left ( F/\O_F, f_a^* G \times_{\widetilde{\DD}_{a,\infty}}
({\widetilde{\DD}_{a,\infty}} \text{ mod }\pi ) \right  ) [\frac{1}{\pi}]
$$
dont l'image est caractérisée par la théorie de la déformation de
Messing relativement à l'idéal $(\pi)$ qui est muni de $\O$-puissances divisées
(cf. \cite{Messing1} ainsi que l'appendice B de \cite{Cellulaire}). 

\begin{theo}\label{relevisofd24}
Le morphisme $\chi$  se relève en caractéristique zéro via
l'application de réduction modulo $\pi$ précédente :
$$
\chi \in \Hom \left ( F/\O_F, f_a^* G   ) \right  )
[\frac{1}{\pi}]
$$
\end{theo}
\dem 
Rappelons le critère de relèvement de Messing.
Soit $\mathfrak{Z}$ un schéma formel $\pi$-adique sur $\spf
(\O)$. Soient $H_1$ et $H_2$ deux $\O$-modules $\pi$-divisibles sur
$\mathfrak{Z}$. On s'intéresse à l'image de l'injection
$$
\Hom (H_1,H_2)\hookrightarrow \Hom (H_1\text{ mod }\pi, H_2\text{ mod }\pi)
$$
Soit $H\longmapsto \Lie \, E(H)$ le foncteur algèbre de Lie de la
$\O$-extension vectorielle universelle (appendice B de \cite{Cellulaire}) muni de sa partie vectorielle $\Fil \, \Lie\, E(H)$. D'après la nature cristalline de la $\O$-extension
vectorielle universelle il y a un morphisme
$$
\Hom (H_1\text{ mod }\pi,H_2\text{ mod }\pi) \ldrt
\Hom_{\O_\mathfrak{Z}} ( \Lie\, E(H_1),\Lie\, E(H_2))
$$
Alors le critère de relèvement dit que le diagramme suivant est cartésien
$$
\xymatrix{
\Hom (H_1,H_2) \ar@{^(->}[r] \ar [d] & \Hom (H_1\text{ mod
}\pi,H_2\text{ mod }\pi) \ar[d] \\
  \Hom_{\O_\mathfrak{Z}-\text{ modules filtrés}} ( \Lie\, E(H_1),\Lie\, E(H_2))\ar@{^(->}[r] & \Hom_{\O_\mathfrak{Z}} ( \Lie\, E(H_1),\Lie\, E(H_2))
}
$$
Pour le groupe $F/\O_F$ on a $\Lie\, E(F/\O_F)=\Fil\, \Lie (F/\O_F) = \O_\mathfrak{Z}$. 

Revenons à la démonstration. Il s'agit de montrer que via la composée
$$
f_a^* \rho_G \circ \Delta \circ \rho_H^{-1}\circ (\zeta_1,\dots,\zeta_n)
 : F/\O_F \ldrt f_a^* G \text{ mod }\pi
$$
le morphisme induit au niveau de l'évaluation des cristaux sur l'épaississement 
$\widetilde{\DD}_\infty \text{ mod }\pi \hookrightarrow \widetilde{\DD}_\infty$
\begin{eqnarray*}
\O_{\widetilde{\DD}_\infty}\unpi &\ldrt & f_a^*\Lie \, E(G) \unpi \\
1 &\longmapsto & x
\end{eqnarray*}
vérifie $x\in f_a^*( \Fil\, \Lie E(G))\unpi$.

La quasi-isogénie $f_a^*\rho_G$ induit un isomorphisme
$$
\DD (\Gb)\otimes \O_{\widetilde{\DD}_\infty} \unpi \iso f_a^* \Lie \, E(G) \unpi
$$
Alors via cet isomorphisme 
$$
x\in \DD (\Gb)\otimes \O_{\widetilde{\DD}_\infty}\unpi = \bigoplus_{j\in \Z/n\Z} \DD (\Gb)_{\Q,j}\otimes  \O_{\widetilde{\DD}_\infty} \unpi
$$
La filtration définissant le morphisme $\widetilde{\DD}_\infty \ldrt \widehat{\Omega}$
(où $\widehat{\Omega}$ est l'espace de Drinfeld associé à l'espace vectoriel $\DD (\Gb)_{\Q,0}^{V=\Pi}$) 
 est une filtration localement facteur direct $\Fil \subset \DD (\Gb )_{\Q,0} \O_{\widetilde{\DD}_\infty} \unpi$ et via l'isomorphisme précédent
$$
f_a^* (\Fil\, \Lie \, E(G))\unpi = \bigoplus_{j\in \Z/n\Z} \Pi^j \Fil
$$
Rappelons qu'on a fixé un isomorphisme $\DD( \Hb)_{\Q,0}\simeq \breve{F}$ et que via l'isogénie  $\Delta$ cela nous a permis d'identifier
$$
\breve{F}^n \simeq \DD(\Hb)_{\Q,0}^n \underset{\simeq}{\xrig{\;\Delta_*\;}} \DD( \Gb)_{\Q,0}
$$
Donc le morphisme définissant $f_a:\widetilde{\DD}_\infty \ldrt \widehat{\Omega}$ 
(où maintenant il s'agit de l'espace $\widehat{\Omega}$ associé à l'espace vectoriel $F^n$) 
est donné  par une filtration
$$
\Fil' \subset \DD (\Hb)_{\Q,0}^n \otimes  \O_{\widetilde{\DD}_\infty} \unpi \simeq  \O_{\widetilde{\DD}_\infty} \unpi^n
$$
et il s'agit de voir que le morphisme composé 
$$
F/\O_F \xrig{\; (\zeta_1,\dots,\zeta_n)\;} H^n\times (\widetilde{\DD}_\infty \text{ mod }\pi) \ldrt \Hb^n \times (\widetilde{\DD}_\infty \text{ mod }\pi )
$$
induit au niveau de l'évaluation des cristaux un morphisme tel que 
$$
1\longmapsto X\in \bigoplus_{j\in \Z/n\Z} \Pi^j \Fil' \subset \DD (\Hb)^n\otimes   \O_{\widetilde{\DD}_\infty} \unpi
$$
Si $\a_H$ désigne l'application de Hodge-Tate de $H$ alors
$$
X=(\a_H (\zeta_i))_{1\leq i\leq n} \in \omega_{H^\vee} \otimes \O_{\widetilde{\DD}_\infty}[\frac{1}{\pi}]^n \subset \DD (\Hb)_\Q^n  \otimes \O_{\widetilde{\DD}_\infty}[\frac{1}{\pi}]^n
$$
et 
$$
X=(X_{ij})_{1\leq i\leq n,j\in \Z/n\Z}\in M_n ( \O_{\widetilde{\DD}_\infty}[\frac{1}{\pi}] )
$$
via 
$$
\DD (\Hb)_\Q^n \iso \bigoplus_{j\in \Z/n\Z} \DD (\Hb)_{\Q,j}^n \underset{\simeq}{\xrig{\; \sum \Pi^{-j}\;}} \bigoplus_{j\in \Z/n\Z} \DD( \Hb)_{\Q,0}^n \simeq \bigoplus_{j\in \Z/n\Z} \breve{F}^n
$$
Rappelons maintenant qu'on a identifié $F^n$ à $(F^n)^*$ et qu'alors $\Fil'\subset  \O_{\widetilde{\DD}_\infty}[\frac{1}{\pi}]^n$ est le noyau de $\a_{H^\vee}(-1)$. On vérifie alors que via l'isomorphisme 
$$
\DD (\Hb)_{\Q,0}^n\otimes  \O_{\widetilde{\DD}_\infty}[\frac{1}{\pi}] \simeq  \O_{\widetilde{\DD}_\infty}[\frac{1}{\pi}]^n \simeq \left (  \O_{\widetilde{\DD}_\infty}[\frac{1}{\pi}]^n \right )^*
$$
on a 
$$
\forall j\in \Z/n\Z\;\;\; (X_{ij})_{1\leq i\leq n} = \,^t \a_H (\ph_j )
$$
où $\ph_j \in \omega_{H^\vee}^*\otimes  \O_{\widetilde{\DD}_\infty}[\frac{1}{\pi}]$ est la forme linéaire composée
$$ \ph_j: 
\omega_{H^\vee}\otimes  \O_{\widetilde{\DD}_\infty}[\frac{1}{\pi}] \hookrightarrow \DD (\Hb)_\Q \otimes  \O_{\widetilde{\DD}_\infty}[\frac{1}{\pi}] \twoheadrightarrow \DD (\Hb)_{\Q,j}\otimes  \O_{\widetilde{\DD}_\infty}[\frac{1}{\pi}] \xrig{\;\Pi^{-j}\;}  \DD (\Hb)_{\Q,0}\otimes  \O_{\widetilde{\DD}_\infty}[\frac{1}{\pi}] \simeq \O_{\widetilde{\DD}_\infty}[\frac{1}{\pi}]
$$
Le théorème est donc une conséquence de la proposition \ref{exhodateg}.
\qed

\begin{prop}\label{countyrog}
Soit $f:\widetilde{\X}_\infty \ldrt \widehat{\Omega}$. Les morphismes $\chi$ définis sur les différentes cellules se recollent en un
$$
\chi\in \GG ( \widetilde{\X}_\infty, \underline{\Hom} (F/\O_F, f^* G)[\frac{1}{\pi}])
$$
au sens où $\forall \mathcal{U}$ ouvert quasicompact dans $\widetilde{\X}_\infty$ $\;
\chi\in \Hom (F/\O_F,f^* G_{|\mathcal{U}}) [\frac{1}{\pi}]$. 
Pour tout ouvert quasicompact $\mathcal{U} \subset \widetilde{\X}_\infty$ $\;\; \chi_{|\mathcal{U}}\neq 0$. Pour tout $g\in \GL_n (F)$ via l'égalité $g^*f^* G= f^* G$ où $g:\widetilde{\X}_\infty \ldrt \widetilde{\X}_\infty$ on a $g^* \chi =\chi$.
\end{prop}
\dem
Les assertions concernant le recollement et l'équivariance ``$g^*\chi=\chi$'' ne posent pas de problème; il suffit de les vérifier dans $\underline{\Hom} (F/\O_F, f^* G\text{ mod }\pi)[\frac{1}{\pi}]$. De même l'assertion $\chi_{|\mathcal{U}}\neq 0$ se vérifie modulo $\pi$ où, avec les notations de la démonstration précédente, puisque $\rho_H, \rho_G$ et $\Delta$ sont inversibles, elle est équivalente à
$$
\forall a=[\La,M]\;\forall \mathcal{U}\subset \widetilde{\DD}_{a,\infty}\;\;
\forall \mathcal{U} \subset \widetilde{\DD}_{a,\infty}\;\; (\zeta_1,\dots,\zeta_n): F/\O_F \ldrt 
H^n\times_{\widetilde{\DD}_{a,\infty}} (\mathcal{U} \text{ mod }\pi) \;\text{ est non-nul}
$$
Cette assertion est claire puisque sur $\mathcal{U}$ $\;\forall k\geq
1\; (\zeta_1,\dots,\zeta_n)$ définissent une structure de niveau de
Drinfeld $\pi^{-k}\La/\La \ldrt H[\pi^k]\times \mathcal{U}$, mais
$\mathcal{U}$ étant sans $\pi$-torsion une telle structure de niveau
de Drinfeld est non-triviale après inversion de $\pi$ et est donc
non-triviale. 
\qed

\subsection{Construction du morphisme}

\subsubsection{Un canular}

La remarque qui suit justifie que nous devrons procéder à quelques
vérifications concernant les composantes connexes des revêtements de
la tour de Lubin-Tate. 

Plaçons nous dans le cas $n=1$ et $F=\Qp$ et, afin de simplifier les
notations, travaillons sur $F$ et non $\breve{F}$. 
Soit $R$ une $\Zp$-algèbre $p$-adique sans $p$-torsion telle que
$\spec (R\unp)$ soit connexe. Soit $x\in \Hom ( \Qp/\Zp,
\mu_{p^\infty})\unp$ tel que $x\neq 0$.
On a vu dans la proposition
\ref{crutozp} qu'un tel $x$ induit un morphisme 
$$
\spf (R) \ldrt \mathcal{Y}_\infty = \coprod_{\Z} \spf (\widehat{\Zp^{\text{ab}}})
$$
Rappelons en effet qu'il existe 
 $ n\in\Z$ tel que 
$$
p^n x : \Qp/\Zp \iso \mu_{p^\infty/R\unp}
$$
et alors $\forall k\geq 1\;\; p^n x: p^{-k}\Z/\Z\iso \mu_{p^k/R\unp}$
d'où 
$$
\xymatrix@R=6mm{
 &  \mathcal{Y}_\infty \ar@{..>}[d] \\
 & \mathcal{Y}_{k+1}   \ar[d] \\
\spf (R) \ar[r] \ar[ru] \ar[ruu]^\a &  \mathcal{Y}_{k}
}
$$
et on définit alors le morphisme cherché par
$$
\spf (R) \xrig{\;\a\;} \mathcal{Y}_\infty \xrig{\;p^{-n}\;} \mathcal{Y}_\infty
$$
où l'action de $p^{\Z}$ su $\mathcal{Y}_\infty$ se fait par
translations des composantes de $\coprod_{\Z} \spf
(\widehat{\Zp^{\text{ab}}})$. 
\\

Considérons maintenant l'exemple suivant. Soit 
$$
\spf (R)=\underset{n}{\limp} \coprod_{\Z/p^n\Z} \spf ( \widehat{\Zp^{\text{ab}}})
$$
dans la catégorie des schémas formels $p$-adiques, 
c'est à dire 
$$
R = \left ( \underset{n}{\limi} \prod_{\Z/p^n \Z}
  \widehat{\Zp^{ab}}\right )^{\hspace{-2mm}\widehat{\;\;\;\;\;}} = \mathcal{C}^0 ( \Zp, \widehat{\Zp^{\text{ab}}})
$$
les fonctions continues  de $\Zp$ dans $\widehat{\Zp^{\text{ab}}}$.
Fixons $(\zeta_{p^k})_{k\geq 1} \in \mu_{p^\infty}
(\widehat{\Zp^{\text{ab}}})$ un générateur de $T_p (\mu_{p^\infty})$,
c'est à dire tel que $\zeta_p\neq 1$. 
Soit $x:\Qp/\Zp\ldrt \mu_{p^\infty /R}$ défini par $x=(x_k)_{k\geq 1},
x_k\in \mu_{p^k} (R)$  et $x_k$ est la fonction continue 
$$
\forall t\in \Zp\;\; x_k (t) = \zeta_{p^k}^t
$$
Alors, pour tout ouvert $\mathcal{U}$ de $\spf (R)$ $\;
x_{|\mathcal{U}}\neq 1$. Mais il n'existe pas de morphisme associé 
$\spf (R)\ldrt \mathcal{Y}_\infty$ car les composantes connexes de
$\spec (R)$ ne sont pas ouvertes (les ``diracs'' comme fonctions de
$\Zp$ à valeurs dans $\widehat{\Zp^{\text{ab}}}$ ne sont pas des
fonctions continues). On a $\pi_0 (\spec (R))\simeq \Zp$.

\subsection{Un remède au canular}\label{okk65tr}

\begin{prop}\label{canulartt}
Soit $\X$ l'espace de Lubin-Tate sans niveau sur $\spf
(\breve{\O})$ et pour $C\subset \GL_n (\O_F)$ un sous-groupe ouvert
$\X^{rig}_C$ l'espace de
Lubin-Tate rigide en niveau $C$. 
 Soit $H$ le groupe $p$-divisible universel sur
$\X$. Soit $U\subset \X^{rig}_C$ un ouvert admissible connexe. Soit
$s\in U(\overline{\breve{F}})$ un point géométrique et
$$
\rho : \pi_1 (U,s) \ldrt \GL_{\O_F} (T_p ( H^{rig}))
$$
la représentation de monodromie associée. Alors l'image de $\rho$ est
ouverte. En d'autres termes si $(U_{C'})_{C'\subset C}$ désigne l'image
réciproque de $U$ dans les revêtements de la tour de Lubin-Tate en
niveaux $C'\subset C$ alors 
$$
\underset{C'\subset C}{\limp} \pi_0 (U_{C'})
$$
est un ensemble fini (pour $C''\subset C'$  avec $C'\subset C$ suffisamment petit $\;\pi_0 (U_{C''})\iso \pi_0 (U_{C'})$). 
\end{prop}
\dem
Il suffit de montrer qu'il existe une extension de degré fini
$K|\breve{F}$ ainsi qu'un point $x\in U(K)$ tel que si $H_x$ désigne
le groupe de Lubin-Tate sur $\O_K$ spécialisé en $x$ alors l'image de 
$$
\rho_x: \Gal (\overline{K}|K)\ldrt \GL_{\O_F} ( T_p (H_x))
$$
est ouverte. Du point de vue des composantes connexes cela exprime que
si $x'\in U$ désigne le point ``du spectre maximal''
 associé à $x$ l'ensemble
$\underset{k}{\limp} \Pi_{C',C}^{-1} (x')$ est fini où $\Pi_{C',C}$ est
l'application de l'espace de Lubin-Tate en niveau $C'$ vers 
celui en niveau $C$. Il suffit de le faire pour $C=\GL_n (\O_F)$
puisque l'image de $U$ dans $\X^{rig}$ par l'application d'oubli du
niveau est un ouvert admissible. Supposons donc $C=\GL_n (\O_F)$. 

Si $K|\breve{F}$ est finie et $x\in \X^{rig} (K)$ soit $\Fil_x\subset
\DD ( \Hb)_\Q\otimes_{\breve{F}} K$, $\Fil_x \in \mathbb{P}^{n-1}
(K)$, la filtration de Hodge associée dans l'espace des périodes.
Alors 
$$
\End (H_x)_\Q \simeq \End_{\Gal (\overline{K}|K)} (V_p (H_x)) \simeq
\End ( \DD (\Hb)_\Q,\ph,\Fil_x) \simeq \text{Stab}_{D} ( \Fil_x)
$$
qui est un corps commutatif $E|F$ tel que $E\subset D$ (le fait que
$E$ est commutatif résulte de ce qu'étant donné que $H_x$ est un
groupe formel l'application qui à un endomorphisme
associe l'endomorphisme tangent sur l'algèbre de Lie est injective,
mais ici l'algèbre de Lie de $H_x$ est de dimension $1$). 

Supposons $E=F$. Alors  puisque l'isocristal de $H_x$ est simple la représentation cristalline $V_p (H_x)$ est
irréductible et le reste comme représentation $\Gal (\overline{K}
|K')$ pour toute extension de degré fini $K'$ de $K$. 
 De plus pour toute extension de degré
fini $K'|K\;$  
$\End_{\Gal (\overline{K}|K')} (V_p
(H_x))=F$. Cela implique par la théorie de Sen (cf. par exemple
le théorème 5 de \cite{Serre2} où le cas $F=\Qp$ est traité, le cas $F$
général étant identique)
que l'algèbre de Lie de l'image de $\Gal (\overline{K}|K)$  est $End (
V_p (H_x))$ et que donc l'image est ouverte.

Reste donc à voir que si $\breve{\pi}:\X^{rig}\ldrt \mathbb{P}^{n-1}$
désigne l'application des périodes alors
$$
\exists x\in U\;\; \text{Stab}_{D^\times} ( \breve{\pi} (x)) =F^\times 
$$
Mais le morphisme des périodes étant étale $\breve{\pi} (U)$ contient un
ouvert admissible de $\mathbb{P}^{n-1}$. Or les sous-variétés ``de
type Hodge'' associées aux corps $E$ tels que $F\subset E\subset D$
sont les points fixes de $E^\times$ agissant sur
$\mathbb{P}^{n-1}$ via $E^\times \subset D^\times$.
Lorsque $E$ varie 
elles  forment une union
finie de $D^\times$-orbites de
sous-variétés algébriques dans $\mathbb{P}^{n-1}$.
Ces $D^\times$-orbites sont en bijection avec les 
classes de conjugaison de tores non-triviaux
du groupe algébrique associé à $D^\times$.
 Donc d'après la
proposition qui suit
$\breve{\pi} (U)$ (ou plutôt ses points ``classiques'') 
n'est pas contenu dans l'union de ces
sous-variétés. 
\qed

\begin{prop}
Soit $L$ un corps valué complet pour une valuation discrète. Soit $X$
un $L$-espace analytique de Berkovich
lisse équidimensionnel. Soit $H$ un groupe
topologique compact agissant continûment sur $X$ au sens où
l'application $H\times |X|\ldrt |X|$ est continue. Soit $Z\subset X$
un sous-ensemble compact tel que tout $x\in X$ possède un voisinage
$U$ tel que $U\cap Z$ soit un sous-ensemble analytique Zariski fermé
dans $U$ de dimension $< \dim X$. Alors pour tout domaine analytique
$V\subset X$ non-vide $V(\overline{L})\setminus H.Z(\overline{L})\neq \emptyset$.
\end{prop}
\dem 
L'ensemble $H.Z$ étant compact dans $|X|$, $V\setminus H.Z$ est ouvert
dans $V$. Or l'image de $V(\overline{L})$ par l'application
$V(\overline{L})\ldrt |V|$ est dense dans $|V|$. Il suffit donc de montrer
que $V\setminus H.Z\neq \emptyset$. Soit $ \mathcal{M}(\mathcal{A})\subset V$ un
domaine affinoïde. Puisque la valuation de $L$ est discrète 
la $\O_L$-algèbre des éléments topologiquement
bornés $\mathcal{A}^0$ dans $\mathcal{A}$ est une $\O_L$-algèbre
  telle que si
$\widetilde{\mathcal{A}} = \mathcal{A}^0/\mathcal{A}^{00} =
(\mathcal{A}^0/\pi_L\mathcal{A}^0)_{\text{red}}$ désigne sa fibre
spéciale réduite on ait
$$
\dim \mathcal{A} = \dim \widetilde{\mathcal{A}}
$$
(cf. le lemme \ref{affghb}). 
Soit $\text{sp} : \mathcal{M} (\mathcal{A})\ldrt \spec
(\widetilde{\mathcal{A}})$ l'application de spécialisation. Soit $\xi$
un point générique d'une composante irréductible 
de $\spec (\widetilde{\mathcal{A}})$ de dimension $\dim X$.
D'après la proposition 2.4.4. page 36 de \cite{BerkSpectral} il existe $x\in\mathcal{M} (\mathcal{A})$ tel que
$\text{sp} (x)=\xi$. 
 De plus d'après le lemme \ref{affghb} qui suit 
 pour un tel $x$ le corps résiduel du corps
valué complet $\mathcal{K} (x)$ a pour degré de transcendance $\dim
X$ sur le corps résiduel de $L$.
 De cela on déduit que $x\notin H.Z$ puisque toujours d'après le lemme
 qui suit 
$\forall z\in Z\;\;
\text{deg.tr.} \widetilde{\mathcal{K}(z)}< \dim X$.
\qed

\begin{lemm}\label{affghb}
Soit $L$ un corps valué complet pour une valuation discrète de corps
résiduel $k$. Soit $\mathcal{A}$ une $L$-algèbre affinoïde. Soit
$\mathcal{A}^0 =\{ f\in\mathcal{A}\;|\; \| f\|_\infty\leq 1\}$,
$\mathcal{A}^{00} =\{ f\in\mathcal{A}\;|\; \| f\|_\infty < 1\}$ et
$\widetilde{\mathcal{A}}=\mathcal{A}/\mathcal{A}^0$. Alors 
 $\widetilde{\mathcal{A}}$ est une $k$-algèbre de type fini.  
De plus
$$
\dim\, \mathcal{A} = \dim\, \widetilde{\mathcal{A}}
$$
Soit $\mathcal{M}(\mathcal{A})$ l'espace de Berkovich associé à
$\mathcal{A}$ et $\forall x\in\mathcal{M} (\mathcal{A})$
$\,\mathcal{K} (x)$ le corps résiduel de $x$ (ou son complété). Alors
$$
\forall x\in\mathcal{M} (\mathcal{A}) \;\; \text{deg.tr.}_k
\widetilde{\mathcal{K}(x)} \leq \dim \,\mathcal{A} 
$$
Si de plus $sp (x)$ est un point générique d'une composante
irréductible de dimension $\dim \,\widetilde{\mathcal{A}}$ 
 alors
$$
\text{deg.tr.}_k
\widetilde{\mathcal{K}(x)} = \dim \,\mathcal{A}
$$
\end{lemm}
\dem
D'après \cite{BGR}, théorème 1 section 6.3.5, un morphisme d'algèbres affinoïdes
$\mathcal{B}\ldrt \mathcal{C}$ est fini ssi le morphisme induit 
$\widetilde{\mathcal{B}}\ldrt \widetilde{\mathcal{C}}$ l'est. On en
déduit que  $\widetilde{\mathcal{A}}$ est une $k$-algèbre de type
fini puisque c'est le cas de $\mathcal{A}=L<T_1,\dots,T_n>$ pour
tout $n$. 
\\
D'après le théorème de normalisation de Noether il existe un morphisme
injectif fini $\ph : L<T_1,\dots,T_n>\hookrightarrow \mathcal{A}$ où
$n=\dim \mathcal{A}$. Le morphisme $\ph$ étant fini injectif c'est une
isométrie (lemme 6 p.170 de \cite{BGR}). Donc,
$\widetilde{\ph}:k[T_1,\dots,T_n] \ldrt \widetilde{\mathcal{A}}$ est
injectif fini ce qui implique que $\dim
\,\widetilde{\mathcal{A}}=n=\dim\, \mathcal{A}$. 
\\
Maintenant si $\ph^* : \mathcal{M} ( \mathcal{A})\ldrt \mathbb{B}^n$
désigne le morphisme d'espaces de Berkovich associé à un $\ph$ comme
précédemment $\forall x\in \mathcal{M} ( \mathcal{A})$ si $y=\ph^*
(x)$ alors $\mathcal{K} ( x) |\mathcal{K} (y)$ est une extension de
degré fini. Donc l'extension de corps résiduels 
 $\widetilde{\mathcal{K}} (x) | \widetilde{\mathcal{K}} (y)$ est
 algébrique. On est alors ramené à montrer que $\forall y\in
 \mathbb{B}^n$ le degré de transcendance sur $k$ de
 $\widetilde{\mathcal{K}} (y)$ est plus petit que $n$. On procède
 par récurrence sur $n$. Le cas $n=1$ ne pose pas de problème car
 on a un description complète de $\mathcal{M} (L<T>)$ (section 1.4.4
 p.18 de \cite{BerkSpectral}). La récurrence se fait alors en utilisant
 le cas $n=1$ pour d'autres corps que $L$ et 
la projection $pr: \mathbb{B}^{n}\ldrt \mathbb{B}^{n-1}$ qui à
 $(x_1,\dots,x_n)$ associe $(x_1,\dots,x_{n-1})$.  
 En effet, pour $z\in \mathbb{B}^n$, la fibre au dessus de $pr (z)$ de
 $pr$ est isomorphe à $\mathbb{B}^1\hat{\otimes}_{L} \mathcal{K} ( pr
 (z))\ni z$.
\\
Supposons maintenant que $sp (x)$ est un point générique d'une
composante irréductible de $\widetilde{\mathcal{A}}$ de dimension
maximale. Alors étant donné que $\text{deg.tr.}_k
\widetilde{\mathcal{K}(x)} \leq \dim \,\mathcal{A}$ et que $
\widetilde{\mathcal{K}} (x) | k(sp (x))$ on a égalité des deux degrés
de transcendance.
\qed

\begin{rema}
Le principe des démonstrations précédentes consiste à vérifier
que le ``groupe de Mumford-Tate $p$-adique'' est génériquement le plus
gros possible où génériquement signifie que tout ouvert admissible
possède un point où ce groupe est maximal.   
\end{rema}

\subsection{Construction du morphisme}

\begin{lemm}
Soit $R_\infty = \underset{m\geq 1}{\limi} R_m$ où $\forall m\;  R_m$
est une $\breve{\O}$-algèbre $\pi$-adique sans $\pi$-torsion
intégralement fermée dans $R_m \unpi$. Alors 
$\widehat{R}_\infty$ est intégralement fermé dans $\widehat{R}_\infty
[\frac{1}{\pi}]$. 
\end{lemm}
\dem
On vérifie successivement que pour les anneaux sans $\pi$-torsion
la propriété d'être intégralement fermé
dans sa fibre générique (i.e. après inversion de $\pi$)
est stable par limite inductive et complétion $\pi$-adique.
\qed

\begin{lemm}
Tout ouvert quasicompact de $\widetilde{\X}_\infty$ possède un nombre
fini de composantes connexes et $\O_{\widetilde{\X}_\infty}$ est
intégralement fermé dans $\O_{\widetilde{\X}_\infty}\unpi$. 
\end{lemm}
\dem
Soit $k_0\geq 4$. 
Soit $\mathcal{U}_\infty \subset \widetilde{\X}_\infty$ un ouvert
affine de la forme 
$$
\mathcal{U}_\infty = \underset{k\geq k_0}{\limp} \mathcal{U}_k
$$
où $\mathcal{U}_k \subset \widetilde{\DD}_{a,k}$ pour un $a=[\La,M]$
et
$\forall k\;\mathcal{U}_{k+1}$ est l'image réciproque de
$\mathcal{U}_k$ via le morphisme de changement de niveau. 
Les $\mathcal{U}_k$ sont des schémas formels admissibles sur
$\breve{\O}$. De plus les $(\mathcal{U}_k^{rig})_k$ forment un
revêtement pro-galoisien de groupe 
 $K=Id +\pi^{k_0}\End (\La) \subset \GL_n (F)$. 
Rappelons que les $\mathcal{U}_k$ étant normaux $\forall k\; \pi_0
(\mathcal{U}_k)=\pi_0 (\mathcal{U}_k^{rig})$.
\'Ecrivons $\mathcal{U}_k=\spf (R_k)$, $\mathcal{U}_\infty= \spf
(R_\infty)$ avec
$R_\infty = (\underset{k}{\limi} R_k)^{\widehat{\;\;}}$. Alors,
d'après la proposition \ref{canulartt} 
$$\exists k_1\geq k_0 \; \forall k\geq
k_1\;\;\;
\pi_0 (\spec (R_\infty ))\iso \pi_0 (\spec (R_k))$$
qui est un ensemble fini et donc les composantes connexes de $\spec
(R_\infty)$ sont ouvertes. Maintenant si $f\in R_\infty/\pi R_\infty$
et $D(f)\subset \spf (R_\infty)$ est l'ouvert associé alors $\exists
k'$ tel que $f$ provienne d'un élément de $R_{k'}/\pi R_{k'}$ et donc
l'ouvert $D(f)$ également. De cela on déduit 
d'après l'analyse précédente appliquée à $D(f)$ 
 que $D(f)$ possède un nombre fini de
composantes connexes. Donc, $\widetilde{\X}_\infty$ possède une base
d'ouverts ayant un nombre fini de composantes connexes duquel on
déduit que tout ouvert quasicompact possède un nombre fini de
composantes connexes.
\\
D'après le lemme précédent on obtient également ainsi la seconde
assertion (cf. plus généralement la proposition \ref{dljbyryt} de
l'appendice). 
\qed
\\

L'élément $\chi$ de la proposition \ref{countyrog}
 fournit alors d'après la proposition
\ref{crutozp} un morphisme  
$\GL_n (F)\times D^\times$-équivariant 
$$
\widetilde{\X}_\infty \ldrt \mathcal{Y}_\infty
$$
où équivariant signifie que 
l'action de $(g,d)\in \GL_n (F)\times D^\times$ est transformée en
celle de $(\,^t g,d^{-1})$. 

\section{Construction du morphisme $\widetilde{\mathcal{Y}}_\infty \ldrt
  \widehat{\mathbb{P}}^{n-1}$}\label{kfjeg264}

Dans cette section nous entamons la construction du morphisme de la
tour de Drinfeld vers celle de Lubin-Tate en construisant un morphisme
d'un éclaté de $\mathcal{Y}_\infty$ vers l'espace des périodes associé
à l'espace de Lubin-Tate.
\\

\subsection{Applications de Hodge-Tate}

Soit comme précédemment $G$ le $\O_D$-module formel spécial universel
sur $\mathcal{Y}=\coprod_{\Z} \widehat{\Omega}$. Notons $\forall k\;
\mathcal{Y}_k = \mathcal{Y}_{1+\pi^k \O_D}$. Les schémas formels
$\mathcal{Y}_k$ étant normaux la structure de niveau en fibre
générique induit des morphismes
$$
\forall k\geq 1\;\;\;\; \underline{\pi^{-k}\O_D/\O_D}\ldrt
G[\pi^k]\times_{\mathcal{Y}} \mathcal{Y}_k
$$
qui induisent des isomorphismes sur $\mathcal{Y}_k^{rig}$. 
\\

{ \it Fait admis : Comme pour la tour de Lubin-Tate on admettra que
$\O_F/\pi^k \O_F (1)$ devient trivial sur $\mathcal{Y}_k$. Cela
résulte de l'existence d'une application déterminant que l'auteur
espère construire dans \cite{Periodes} }
\\

Comme dans la section \ref{chaverztu} on construit à partir de ces morphismes des
applications de Hodge-Tate 
$$
\forall k\geq 1\;\; \a_{G^\vee [\pi^k]}(-1) : \Hom_{\O_F} (\O_D, \O_F)
\otimes 
\O_{\mathcal{Y}_{k}} /\pi^k \O_{\mathcal{Y}_{kn}}  \ldrt
\omega_{G} \otimes \O_{\mathcal{Y}_{k}} /\pi^k \O_{\mathcal{Y}_{k}} (-1) 
$$
Ces morphismes satisfont une condition de compatibilité naturelle  et
induisent alors une application de Hodge-Tate en niveau infini 
$$ 
\a_{G^\vee}(-1) :\Hom_{\O_F} (\O_D,\O_F)\otimes \O_{\mathcal{Y}_\infty}  \ldrt
\omega_{G}\otimes \O_{\mathcal{Y}_\infty}(-1)
$$
L'action à gauche de $\O_D$ sur $G$ induit une action à droite de
$\O_D$ sur $G^\vee$, c'est à dire un morphisme $\O_D\ldrt \End
(G^\vee)^{\text{opp}}$. Le module de Tate de $G^\vee$ est donc un
$\O_D$-module à droite et via sa rigidification ce module est
 le $\O_D$-module  $\Hom (\O_D,\O_F)(1)$ où $\forall d\in \O_D\; \forall
 h\in\Hom (\O_D,\O_F)(1)  \;\; h.d
(\bullet) = h(d \bullet )$. De même $\omega_G$ est un $\O_D$-module à
droite. Les morphismes  de Hodge-Tate précédentes sont
$\O_D$-équivariantes 
pour les actions précédentes. 

Rappelons que pour un $\O_D\otimes_{\O_F} \breve{\O}$-module $M$ on
note  $M=\bigoplus_{j\in \Z/n\Z} M_j$ sa décomposition
en facteurs directs où $\O_{F_n}\subset  \O_D$ agit sur le facteur $M_j$ via
$\s^{-j}:F_n\hookrightarrow \breve{F}$.

On a donc des décompositions
$$
 \a_{G^\vee [\pi^k]}(-1) : 
\bigoplus_{j\in \Z/n\Z} \Hom_{\O_{F_n},\s^{-j}} (\O_D,
\breve{\O})\otimes_{\breve{\O}} \O_{\mathcal{Y}_k} /\pi^k
\O_{\mathcal{Y}_k}
\xrig{\; \oplus_j \a_{G^\vee [\pi^k],j} \;} \bigoplus_{j\in \Z/n\Z}
\omega_{G,j}\otimes \O_{\mathcal{Y}_k} /\pi^k
\O_{\mathcal{Y}_k} (-1)
$$
$$
 \a_{G^\vee}(-1) : 
\bigoplus_{j\in \Z/n\Z} \Hom_{\O_{F_n},\s^{-j}} (\O_D,
\breve{\O})\otimes_{\breve{\O}} \O_{\mathcal{Y}_\infty} 
\xrig{\; \oplus_j \a_{G^\vee ,j} \;} \bigoplus_{j\in \Z/n\Z}
\omega_{G,j}\otimes \O_{\mathcal{Y}_\infty} (-1)
$$
sur lesquelles $\Pi$ agit avec $\deg \Pi =-1$.

\subsection{\'Eclatements  et conoyau de l'application de Hodge-Tate}

Nous utiliserons le lemme clef suivant.

\begin{lemm}\label{krezacww}
Soient $M$ et $N$ deux $\O_F$-modules $\pi$-adiquement complets et 
$u,v:M\ldrt N$ deux morphismes tels que $u\equiv v\text{ mod }\pi^2$. Alors
$$
\text{coker } u\text{ est annulé par } \pi \lssi \text{coker } v\text{ est annulé par }\pi
$$
et si c'est le cas $Im \,u = Im \, v$. 
\end{lemm}

De plus nous utiliserons le résultat suivant 

\begin{theo}\label{koyzakox}
Soit $K|F$ un corps valué complet pour une valuation à valeurs dans $\R$. Soit $H$ un $\O$-module $\pi$-divisible sur $\O_K$. Alors le conoyau de l'application de Hodge-Tate
$$
\a_H : T_p (H)\otimes \O_{\widehat{\overline{K}}} \ldrt \omega_{H^\vee} \otimes \O_{\widehat{\overline{K}}} 
$$
est annulé par $\pi$.
\end{theo}

\begin{rema}
En fait nous n'utiliserons le résultat précédent que pour des
extensions $K$ de $\breve{F}$ de degré fini. 
\end{rema}

Considérons maintenant pour tout $j,\, 0\leq j\leq n-1$,
$$
\a_{G^\vee  [\pi^2],j}(-1) : \Hom_{\O_{F_n},\s^{-j}} (\O_D,\breve{\O}) \otimes \O_{\mathcal{Y}_2}/\pi^2
\O_{\mathcal{Y}_2} \ldrt \omega_{G,j}  \otimes \O_{\mathcal{Y}_2}/\pi^2
\O_{\mathcal{Y}_2}  (-1)
$$
Soit $\mathcal{I}_j$ l'idéal cohérent de $\O_{\mathcal{Y}_2}$ tel que 
$\pi^2 \O_{\mathcal{Y}_2}\subset \mathcal{I}_j$ et 
$$
Im ( \a_{G^\vee  [\pi^2],j} )=\omega_{G,j} \otimes \mathcal{I}_j /\pi^2 \O_{\mathcal{Y}_2}
$$

\begin{defi}
Pour tout $k\geq 2$ 
on note $\widetilde{\mathcal{Y}}_k$ le normalisé
de l'éclatement formel admissible des idéaux
$\O_{\mathcal{Y}_k}.\Pi_{k,2}^{-1} \mathcal{I}_j$,
 $0\leq j\leq n-1$,  où $\Pi_{k,2} : \mathcal{Y}_k\ldrt \mathcal{Y}_2$.
\end{defi}

Ainsi $\widetilde{\mathcal{Y}}_k$ est le normalisé du transformé
strict de $\mathcal{Y}_k\ldrt \mathcal{Y}_2$ relativement à
l'éclatement $\widetilde{\mathcal{Y}}_2 \ldrt \mathcal{Y}_2$
(cf. section \ref{sdguiec23} de l'appendice). 
En particulier les morphismes $\widetilde{\mathcal{Y}}_{k+1} \ldrt \widetilde{\mathcal{Y}}_{k}$ sont finis.

\begin{defi}
On note $\widetilde{\mathcal{Y}}_\infty =\underset{k}{\limp} \widetilde{\mathcal{Y}}_k$ dans la catégorie des schémas formels $\pi$-adiques sur $\spf (\breve{\O})$. 
\end{defi}

\begin{rema}
D'après le corollaire \ref{sdmlkgiob35} de l'appendice on peut construire
$\widetilde{\mathcal{Y}}_\infty$ directement en niveau infini. 
\end{rema}

\begin{prop}
Pour tout $k\geq 2$ et $j$, $0\leq j\leq n-1$, 
le morphisme de Hodge-Tate sur l'éclaté $\widetilde{\mathcal{Y}}_k$ 
$$
\a_{G^\vee[\pi^k],j} (-1) : \Hom_{\O_{F_n},\s^{-j}}(\O_D,\breve{\O})\otimes  \O_{\widetilde{\mathcal{Y}}_k}  /\pi^k \O_{\widetilde{\mathcal{Y}}_k}
 \ldrt \omega_{G,j} 
\otimes \O_{\widetilde{\mathcal{Y}}_k}  /\pi^k
\O_{\widetilde{\mathcal{Y}}_k} (-1)
$$
a un conoyau annulé par $\pi$ et son image est de la forme
$\omega_{G,j} 
\otimes \mathcal{J}_{j,k}/\pi^k  \O_{\widetilde{\mathcal{Y}}_k}(-1)$ où
$\mathcal{J}_{j,k}$ est un idéal localement libre de rang $1$ vérifiant
$\mathcal{J}_{j,k}=  \O_{\widetilde{\mathcal{Y}}_k} .\Pi_{k,2}^{-1}
\mathcal{J}_{j,2}$ avec $\Pi_{k,2} : \widetilde{\mathcal{Y}}_k\ldrt
\widetilde{\mathcal{Y}}_2$. 

En niveau infini l'image de
$$
\a_{G^\vee,j}(-1) :   \Hom_{\O_{F_n},\s^{-j}} (\O_D,\breve{\O}) \otimes \O_{\widetilde{\mathcal{Y}}_\infty}
\ldrt \omega_{G,j} 
\otimes \O_{\widetilde{\mathcal{Y}}_\infty} (-1)
$$
est localement libre de rang $1$ égale à $\omega_{G,0} 
\otimes
\O_{\widetilde{\mathcal{Y}}_\infty}.\Pi_{\infty,2}^{-1} \mathcal{J}_{j,2} (-1)$. 
\end{prop}
\dem
Oublions les torsion à la Tate dans cette démonstration. Fixons l'entier $j$. 
Par définition de $\widetilde{\mathcal{Y}}_2$ le morphisme
$$
\a_{G^\vee[\pi^2],j} :  \Hom_{\O_{F_n},\s{-j}}(\O_D,\breve{\O}) \otimes
\O_{\widetilde{\mathcal{Y}}_2}/\pi^2  \O_{\widetilde{\mathcal{Y}}_2} 
\ldrt \omega_{G,j} \otimes
\O_{\widetilde{\mathcal{Y}}_2}/\pi^2   \O_{\widetilde{\mathcal{Y}}_2}
$$
vérifie $Im\, \a_{G^\vee [\pi^2],j} = \omega_{G,j} \otimes\mathcal{J}_{j,2}/ \pi^2 \O_{\widetilde{\mathcal{Y}}_2}$ où $\mathcal{J}_{j,2}$ est localement libre de rang $1$. Montrons que $\pi \O_{\widetilde{\mathcal{Y}}_2} \subset \mathcal{J}_{j,2}$.
Il suffit de le vérifier localement sur $\widetilde{\mathcal{Y}}_2$. Soit donc, localement sur $\widetilde{\mathcal{Y}}_2$ $\; t\in \omega_{G,j}$ une section engendrant $\omega_{G,j}$ et écrivons 
$$
Im ( \a_{G^\vee [\pi^2],j}) = t\otimes \bar{f}
$$
où $f\in \mathcal{J}_{j,2}$. Pour tout $x\in \widetilde{\mathcal{Y}}_2^{rig}$, $x:\spf (\O_K)\ldrt \widetilde{\mathcal{Y}}_2$ pour $K|\breve{F}$ finie, on a $x^*\a_{G^\vee [\pi^2]} = \a_{(x^* G)^\vee [\pi^2]}$ où $x^* G$ vit sur $\spf (\O_K)$. Du théorème \ref{koyzakox} on déduit $|\frac{\pi}{f} (x)|\leq 1$. Donc la fonction rigide $\frac{\pi}{f}$ vérifie
$$
\| \frac{\pi}{f} \|_\infty \leq 1
$$
ce qui implique, $\widetilde{\mathcal{Y}}_2$ étant normal, que $\pi/f \in \O_{\widetilde{\mathcal{Y}}_2}$.

On déduit le reste de la proposition en utilisant le lemme \ref{krezacww}.
\qed
\\

\begin{rema}
Le lecteur attentif aura remarqué que contrairement au cas des espaces
de Lubin-Tate, c'est à dire  le cas de la section \ref{maishegndgj},
on doit d'abord éclater/normaliser avant d'avoir des renseignements
sur le conoyau de l'application de Hodge-Tate. La raison en est que le 
la proposition \ref{rugytok} est plus précise que le théorème 
\ref{koyzakox}. 
\end{rema}

\begin{lemm}
L'action de $\GL_n(F)\times D^\times$ sur
$\mathcal{Y}_\infty$ s'étend en une action de
$\widetilde{\mathcal{Y}}_\infty$. 
\end{lemm}
\dem
En ce qui concerne l'action de $\GL_n(F)\times \O_D^\times$ le
résultat est clair puisque ce groupe laisse invariant les idéaux
$\mathcal{I}_j$. 
Reste à voir que $\forall i, \,0\leq i\leq n-1$, l'action de $\Pi^i \in D^\times$  
s'étend (puisque $\pi \in D^\times$ agit comme $\pi\in \GL_n (F)$
l'action
de $\Pi^\Z$ s'étend ssi c'est le cas pour les $\Pi^i, \,0\leq i\leq
n-1$). 
 Le module de Tate de $(G/G[\Pi^i])^\vee$ est égal à $
\Hom_{\O_F} (\O_D,\O_F).\Pi^i \subset \Hom_{\O_F} (\O_D,\O_F)$.
Il suffit de voir que pour $k\geq 3$ l'image de ce sous-module par
l'application
$\a_{G^\vee [\pi^k]}$ est somme directe de modules localement libres.
Mais cela résulte de 
$$
\a_{G^\vee[\pi^k]} ( \Hom_{\O_F} (\O_D,\O_F).\Pi^i\otimes \O_{\mathcal{Y}_k}/\pi^k \O_{\mathcal{Y}_k} )
= \bigoplus_{j\in \Z/n\Z} Im (\a_{G^\vee [\pi^k],j+i}).\Pi^j
$$
et de ce qu'étant donné que $\widetilde{\mathcal{Y}}_\infty$ n'a pas
de $\pi$-torsion et que $\forall j\;\Pi^j$ possède un inverse après inversion de
$\pi$ alors l'image par  $\Pi^j$ d'un module localement libre de rang
1 est localement libre de rang 1.
\qed

\begin{rema}
Nous n'utiliserons que $\a_{G^\vee,n-1}$ pour définir le morphisme de
$\widetilde{\mathcal{Y}}_\infty$ vers l'espace des périodes de
Lubin-Tate. Mais la raison pour laquelle nous avons rendue localement
libre l'image de tous les $(\a_{G^\vee,j})_{0\leq j\leq n-1}$ provient
du lemme précédent, afin de pouvoir relever l'action de $\Pi^{\Z}$ 
à $\widetilde{\mathcal{Y}}_\infty$. On aurait également pu rendre
l'image de $\a_{G^\vee,n-1}$ localement libre de rang un après
éclatement/normalisation des idéaux
$(\O_{\mathcal{Y}_k}.\Pi_{k,2}^{-1} \mathcal{I}_0)_{k\geq 2}$ puis
remarquer qu'étant donné un idéal $\mathcal{K}$ $\O_D^\times\times \GL_n (F)$-invariant 
l'idéal $\prod_{0\leq j\leq n-1} \Pi^{j*}\mathcal{K}$ est
$D^\times\times \GL_n (F)$-invariant et qu'on peut donc
l'éclater/normaliser de nouveau afin d'obtenir un schéma formel
éclaté/normalisé  
$D^\times\times \GL_n (F)$-invariant. 
\end{rema}

\subsection{Construction du morphisme $\widetilde{\mathcal{Y}}_\infty
  \ldrt \widehat{\mathbb{P}}(\DD (\Hb))$}
\label{iqghuqd}
\subsubsection{Une identification}\label{guyphng}

Rappelons (cf. \cite{Points}) qu'avec les choix faits 
il y a un isomorphisme 
$$
\DD( \Hb)\simeq \O_D\otimes_{\O_{F_n}}\breve{\O}
$$
où si $V$ désigne le Verschiebung et $\iota : \O_D\iso \End (\DD (\Hb),V)$
\begin{eqnarray*}
\forall d\otimes \l\in \O_D\otimes_{\O_{F_n}} \breve{\O}\;\;\;\; V(d\otimes \l)& = & d\Pi\otimes \l^{\s^{-1}} \\
\forall d'\in \O_D\;\;\;\; \iota (d') (d\otimes \l) &= & d' d\otimes \l
\end{eqnarray*}
Avec ces notations
$$
\O_D\otimes_{\O_{F_n}} \breve{\O} =\bigoplus_{0\leq j\leq n-1}\underbrace{ \Pi^j\otimes 1.\breve{\O}}_{\DD ( \Hb)_j}
$$
et donc $\forall j, \, 0\leq j\leq n-1$, $\iota(\Pi^j) : \DD(\Hb)_0\iso \DD (\Hb)_j$ est un isomorphisme. 
De cela on déduit le lemme suivant 

\begin{lemm}
L'application
$$
\Hom_{\O_D} (\O_D,\DD (\Hb))\ldrt \Hom_{\O_{F_n}} (\O_D,\DD (\Hb)_{n-1})
$$
 où
$\DD (\Hb)_{n-1}$ est muni de l'action de $\O_{F_n}$ dérivant de celle
de $\O_D$ sur $\DD (\Hb)$,
qui à $h:\O_D\ldrt \DD(\Hb )$ associe $h$ composé avec la projection
$\DD (\Hb)\twoheadrightarrow \DD (\Hb)_{n-1}$ est une bijection.
\end{lemm}

Il y a donc une bijection
$$
\DD(\Hb)\iso \Hom_{O_{F_n}} (\O_D,\DD(\Hb)_{n-1})
$$
qui à $x=\sum_{j=0}^{n-1} \Pi^{j-n+1}. x_j\in \DD ( \Hb)$ où $\forall
j\; x_j\in \DD(\Hb)_{n-1}$ associe le morphisme de $\O_D$ vers
$\DD(\Hb)_{n-1}$ qui à $\Pi^{n-1-j}$ associe $x_j$. Réciproquement, à
$h\in \Hom_{O_{F_n}} (\O_D,\DD(\Hb)_{n-1})$ on associe
$\sum_{j=0}^{n-1} \Pi^{j-n+1}. h ( \Pi^{n-1-j}) \in \DD (\Hb)$. 

\begin{rema}\label{kro4512}
Via l'isomorphisme précédent $\DD (\Hb ) \iso \Hom_{\O_{F_n}} ( \O_D,
\DD (\Hb)_{n-1})$ l'action de $\O_D$ sur $\DD (\Hb)$ se fait via 
$$
\forall d\in \O_D\; \forall h\in  \Hom_{\O_{F_n}} ( \O_D,
\DD (\Hb)_{n-1}) \;\; (d.h)(\bullet)= h(\bullet d)
$$
\end{rema}

Au final on a donc une identification
\begin{eqnarray*}
\Hom_{\O_{F}} (\O_D, \breve{\O})_{n-1}\otimes_{\breve{\O}} \DD (\Hb)_{n-1}
 &=& \Hom_{\O_{F_n},\s } (
\O_D, \breve{\O})\otimes_{\breve{\O}} \DD (\Hb)_{n-1}  \\
&=& \Hom_{\O_{F_n}} ( \O_D, \DD (\Hb)_{n-1}) \\ 
&=& \DD (\Hb)
\end{eqnarray*}

\begin{rema}\label{qspfihb794}
Plus généralement, considérons le foncteur qui à un
$\O_D\otimes_{\O_F} \breve{\O}$-module $M$ associe le
$\breve{\O}$-module $\Hom_{\O_D\otimes \breve{\O}} ( M, \DD
(\Hb))$. Il y a alors un isomorphisme naturel en $M$ 
$$
\Hom_{\O_D\otimes \breve{\O}} ( M, \DD
(\Hb)) \iso M^*_{n-1}\otimes_{\breve{\O}} \DD (\Hb)_{n-1}
$$
où $M^* = \Hom_{\breve{\O}} (M,\breve{\O})$ et donc $M\mapsto
M^*_{n-1}$ est le foncteur $M\mapsto \Hom_{(\O_{F_n},\s),\breve{\O}} (
M,\breve{\O})$. 
\end{rema}

\subsubsection{Construction du morphisme}\label{qsqfugv354}

Considérons maintenant le morphisme
 $$
\a_{G^\vee,n-1}(-1)\otimes_{\breve{\O}}\DD(\Hb )_{n-1} : \Hom_{\O_{F_n},\s} ( \O_D, \DD (\Hb)_{n-1})\otimes_{\breve{\O}}\O_{\widetilde{\mathcal{Y}}_\infty} \ldrt \omega_{G,n-1}\otimes \DD (\Hb)_{n-1}\otimes \O_{\widetilde{\mathcal{Y}}_\infty} (-1)
$$
D'après l'identification précédente il définit un morphisme 
$$
\DD(\Hb)\otimes_{\breve{\O}}\O_{\widetilde{\mathcal{Y}}_\infty} \ldrt \omega_{G,n-1}\otimes \DD (\Hb)_{n-1}\otimes \O_{\widetilde{\mathcal{Y}}_\infty} (-1)
$$
qui fournit donc un morphisme 
$$
\widetilde{\mathcal{Y}}_\infty \ldrt \widehat{\mathbb{P}}(\DD (\Hb))
$$

\begin{defi}
On muni $\widehat{\mathbb{P}}( \DD (\Hb))$ de l'action à gauche de
$\O_D^\times$ définie par l'action l'action à droite de $\O_D^\times$
sur $\DD (\Hb)$ 
\begin{eqnarray*}
\O_D^\times &\ldrt & \O_D^\times \ldrt \Aut ( \DD (\Hb)) \\
d  &\longmapsto &  d^{-1}
\end{eqnarray*}
\end{defi}

\begin{lemm}
Le morphisme précédent de $\widetilde{\mathcal{Y}}_\infty$ vers
$\widehat{\mathbb{P}}( \DD (\Hb))$ 
est $\O_D^\times$-équivariant au sens où l'action à droite de $d\in \O_D^\times$ sur $\widetilde{\mathcal{Y}}_\infty$ est transformée en celle à gauche de $d^{-1}$ sur $\widehat{\mathbb{P}}(\DD (\Hb))$. 
\end{lemm}
\dem
C'est une conséquence de la remarque \ref{kro4512} couplée à la
définition de l'action de $\O_D^\times$ sur la rigidification du
module de Tate du $\O_D$-module formel spécial universel donnée dans
la définition \ref{zer3476}.
\qed

\section{Relèvement du morphisme $\widetilde{\mathcal{Y}}_\infty\ldrt
  \widehat{\mathbb{P}}^{n-1}$ vers une cellule de l'espace de
  Lubin-Tate} \label{lqvout159}

\subsection{\'Eclatement équivariant de l'espace projectif formel} \label{qdfkobji1247}

 Il est aisé de voir qu'après un éclatement
formel admissible l'action de $\O_D^\times$ sur
$\widehat{\mathbb{P}}(\DD (\Hb))$ 
 s'étend en une action de $D^\times$.
En effet, il suffit d'éclater les idéaux
$$
u(\Pi^j\DD (\Hb)\otimes \O_{\widehat{\mathbb{P}}(\DD (\Hb))} ) (-1)\;\; 1\leq j\leq n-1
$$
où 
$$
u : \DD (\Hb)\otimes \O_{\widehat{\mathbb{P}}(\DD (\Hb))} \twoheadrightarrow \O_{\O_{\widehat{\mathbb{P}}(\DD (\Hb))}} (1)
$$
est l'application universelle. 

Fixons maintenant une cellule $a=[\La,M]$ dans l'immeuble paramétrant les cellules 
de $\X_\infty$. La cellule $\DD_{a}$ ne dépend que de $M$ dans le
couple $[\La,M]$ et non de $\O_F^n$.
Rappelons  que d'après le théorème de Gross-Hopkins
(\cite{HopkinsGross},\cite{Cellulaire}
théorème 1 section 2.3.3) le morphisme des périodes de Hodge De-Rham
$$
\breve{\pi} : \DD^{rig}_a\ldrt \mathbb{P}(\DD (\Hb)_\Q )^{rig}
$$
est un isomorphisme sur son image et que l'ouvert admissible $\breve{\pi} (  \DD_a^{rig})$ vérifie 
$$
\mathbb{P}(\DD (\Hb)_\Q )^{rig} = \bigcup_{0\leq j\leq n-1} \Pi^j. \breve{\pi} (  \DD_a^{rig})
$$
où $\Pi^j.\breve{\pi} ( \DD_a^{rig})= \breve{\pi} (
\DD_{\Pi^j.a}^{rig})$. 
Il existe donc un éclatement formel admissible 
$$
\widetilde{\widehat{\mathbb{P}}} (\DD (\Hb)) \ldrt \widehat{\mathbb{P}} (\DD (\Hb))
$$
muni d'un ouvert $\mathcal{W}\subset \widetilde{\widehat{\mathbb{P}}} (\DD (\Hb))$ tels que l'action de $D^\times$ sur $\mathbb{P}(\DD (\Hb)_\Q )^{rig}$ se prolonge à $\widetilde{\widehat{\mathbb{P}}} (\DD (\Hb))$, 
$$
\widetilde{\widehat{\mathbb{P}}} (\DD (\Hb)) = \bigcup_{0\leq j\leq n-1} \Pi^j.\mathcal{W}
$$
et il y a un isomorphisme 
$$
\psi_a : \DD_a \iso \Pi^j.\mathcal{W}
$$
où si $a=[\La,M]$ $\;M=\Pi^k \O_D$, $k\equiv j\text{ mod }n$,  et
induisant en fibre générique le morphisme des périodes (utiliser
également le
fait que $\DD_a$ est normal pour voir qu'il existe un $\mathcal{W}$;
un isomorphisme $\mathcal{A}\iso \mathcal{B}$ entre algèbres
affinoïdes induit un isomorphisme $\mathcal{A}^0 \iso \mathcal{B}^0$
entre leurs boules unité pour la norme infini). 
 Cette dernière assertion signifie que 
si $\Fil \subset \DD (\Hb)\otimes \O_{\widetilde{\widehat{\mathbb{P}}}(\DD (\Hb))}$ désigne la filtration localement facteur direct de rang $n-1$ définissant l'espace projectif alors
$$
\psi_a^*\Fil_{|\Pi^j.\mathcal{W}} = \Fil\, \DD(\Hb)\unpi
$$
où si $(H_a,\rho_{H_a})$ désigne la déformation universelle sur $\DD_a$,
$$
\rho_{H_a} : \Hb\times_{\Fqb} (\DD_a
\text{ mod }\pi) \ldrt H_a\times_{\DD_a} (\DD_a
\text{ mod }\pi)
$$
$\Fil\, \DD(\Hb)\unpi = (\rho_{H_a*})^{-1} V(H_a)\unpi$ où $\rho_{H_a*}$ est l'isomorphisme
$$
\DD(\Hb)\otimes\O_{\DD_a}\unpi \iso \Lie \,E(H_a)\unpi
$$
induit par $\rho_{H_a}$ via l'évaluation des cristaux sur l'épaississement $(\DD_a
\text{ mod }\pi) \hookrightarrow  \DD_a$ et 
où $E(H_a)$ désigne la $\O$-extension vectorielle universelle de $H_a$ et $V(H_a)$ sa partie vectorielle.

\subsection{Tiré en arrière de l'éclatement de l'espace projectif vers $\widetilde{\mathcal{Y}}_\infty$}

Soit $\mathcal{K}\subset \O_{\widehat{\mathbb{P}}(\DD(\Hb))}$ l'idéal admissible définissant l'éclatement précédent. Soit $h:\widetilde{\mathcal{Y}}_\infty \ldrt \widehat{\mathbb{P}}(\DD(\Hb))$ le morphisme défini précédemment. L'idéal 
$$
\mathcal{K}'= \O_{\widetilde{\mathcal{Y}}_\infty}. h^{-1} \mathcal{K}
$$
vérifie 
$$
\exists N\;\;\; \pi^N\O_{\widetilde{\mathcal{Y}}_\infty}\subset \mathcal{K}' \subset 
O_{\widetilde{\mathcal{Y}}_\infty} 
$$
Quitte à remplacer $\mathcal{K}'$ par $\prod_{j=0}^{n-1} \Pi^{j*}
\mathcal{K'}$ on peut de plus supposer que $\mathcal{K'}$ est tel que
$\Pi^* \mathcal{K'} =\mathcal{K'}$. En effet, $\pi\in D^\times$ agit
sur $\mathcal{Y}_\infty$ comme $\pi\in \GL_n (F)$, or $\mathcal{K}'$
est $\GL_n (F)$-invariant. 
 On fera cette hypothèse. 
De plus le morphisme $h$ étant $\GL_n (F)$-invariant cet idéal l'est. L'idéal $\mathcal{K}$ étant cohérent et ``$\GL_n (F)\bc \widetilde{\mathcal{Y}}_\infty $'' quasicompact au sens où
$\exists \mathcal{V}\subset \widetilde{\mathcal{Y}}_\infty $
 un ouvert quasicompact tel que $\widetilde{\mathcal{Y}}_\infty  =\GL_n (F).\mathcal{V}$ 
on en déduit 
$$
\exists k_0\geq 2\;\; \exists \mathcal{K}''\;\; \pi^N \O_{\widetilde{\mathcal{Y}}_{k_0}} \subset \mathcal{K}'' \subset \O_{\widetilde{\mathcal{Y}}_{k_0}} \text{ et } \mathcal{K}' = \O_{\widetilde{\mathcal{Y}}_\infty}.\Pi^{-1}_{\infty,k_0} \mathcal{K}
$$
où $\Pi_{\infty,k_0} : \widetilde{\mathcal{Y}}_\infty \ldrt
\widetilde{\mathcal{Y}}_{k_0}$ et de plus on peut supposer que $\mathcal{K}''$ est $\GL_n
(F)\times D^\times$-invariant. 

\begin{defi}
Pour tout $k\geq k_0$ on note $\widetilde{\widetilde{\mathcal{Y}}}_k$
le normalisé de l'éclatement formel admissible de l'idéal
$\O_{\widetilde{\mathcal{Y}}_k}.\Pi_{k,k_0}^{-1} \mathcal{K}''$. On
note 
$$
\widetilde{\widetilde{\mathcal{Y}}} =\underset{k\geq k_0}{\limp} \widetilde{\widetilde{\mathcal{Y}}}_k
$$
dans la catégorie des schémas formels $\pi$-adiques sur $\spf
(\breve{\O})$. 
\end{defi}

D'après les propriétés des idéaux données précédemment l'action de
$\GL_n (F)\times D^\times$ s'étend à
$\widetilde{\widetilde{\mathcal{Y}}}_\infty$. 
De plus il y a un morphisme $D^\times$-équivariant et $\GL_n
(F)$-invariant
$$
\xi : \widetilde{\widetilde{\mathcal{Y}}}_\infty \ldrt
\widetilde{\widehat{\mathbb{P}}}(\DD (\Hb))
$$
où $D^\times$-équivariant signifie que l'action à droite de $D^\times$ à la source est
transformée en l'action à gauche au but via $d\mapsto d^{-1}$. 

\begin{rema}
D'après les résultats de l'appendice on peut construire
$\widetilde{\widetilde{\mathcal{Y}}}_\infty$ directement comme le
normalisé dans sa fibre générique de l'éclatement formel admissible de
l'idéal $\mathcal{K}'$ dans $\widetilde{\mathcal{Y}}_\infty$. Cela
évite d'avoir à repasser en niveau fini, passage qui est donc en
quelque sorte artificiel d'après le corollaire \ref{sdmlkgiob35}. 
\end{rema}

\subsection{Relèvement vers la cellule}

Pour tout $j, \, 0\leq j\leq n-1$ si $a=[\La,M]$ est un sommet de l'immeuble
tel que $[\O_D :M]\equiv j \text{ mod n}$ 
on a donc un morphisme 
$$
\psi_a^{-1} \circ \xi :  \xi^{-1} (  \mathcal{W}) . \Pi^{-j} \ldrt \DD_a
$$
où
$$
 \widetilde{\widetilde{\mathcal{Y}}}_\infty = \bigcup_{0\leq j\leq
   n-1} \xi^{-1} ( \mathcal{W}) .\Pi^{-j}
$$

\section{Construction du morphisme
  $\widetilde{\widetilde{\mathcal{Y}}}_\infty \ldrt \X_\infty$}
\label{dsjdmi184}

\subsection{Caractérisation modulaire de $\X_\infty$}

Comme dans la section \ref{rrkcd34} on montre la proposition suivante.

\begin{prop}
Soit $a=[\La,M]$ et $\DD_a$ la cellule associée. Soit $(H_a,\rho_a)$
la déformation universelle sur $\DD_a$. Soit $\DD_{a}=\spf (A)$. On
note encore $H_a$ pour le groupe $p$-divisible sur $\spec (A)$ associé
au groupe $p$-divisible $H_a$ sur $\spf (A)$. 
Soit $Y_k$ le $\spec (A)$-schéma représentant
le faisceau $\underline{\Hom} (\underline{\pi^{-k}\La/\La}, H[\pi^k])$
(après choix d'une base de $\La$ ce schéma est isomorphe à $H[\pi^k]^n$).  
 Soit $U_k\subset Y_{k,\eta}$ l'ouvert/fermé du schéma étale
fibre générique de $Y_k$ défini par
$$
\forall v\in \pi^{-k}\La/\La\setminus \{0\} \;\; \eta (v)\neq 0
$$
où $\eta: \underline{\pi^{-k}\La/\La} \ldrt H[\pi^k]\times_{spec (A)} Y_k$ désigne la
 section universelle.  Alors si $\spec (A_k)$ désigne l'adhérence schématique de $U_k$ dans
$Y_k^{\text{normalisé}}$, c'est à dire l'unique ouvert/fermé de $Y_k^{\text{normalisé}}$ induisant
$U_k$ en fibre générique, on a 
$$
\spf (A_k) \simeq \DD_{a,Id +\pi^k \End (\La)}
$$
\end{prop}

\begin{prop}\label{ghotopolmk}
Soit $\mathfrak{Z}$ un schéma formel $\pi$-adique sur $\spf
(\breve{\O})$ tel que
\begin{itemize}
\item tout ouvert quasicompact de $\mathfrak{Z}$ possède un nombre
  fini de composantes connexes
\item $\O_{\mathfrak{Z}}$ est intégralement fermé dans
  $\O_{\mathfrak{Z}}\unpi $
\end{itemize}
Soit $b$ un sommet de l'immeuble paramétrant les cellules de
$\X_\infty$  
et $f:\mathfrak{Z}\ldrt \DD_b$ un morphisme. Soit $(H,\rho)$ la
déformation universelle sur $\DD_b$. Soit 
$$
(x_1,\dots,x_n)\in \GG( \mathfrak{Z}, \underline{\Hom} (F/\O_F,f^*H )\unpi^n)
$$
tel que $\forall \mathcal{U}\subset \mathfrak{Z}$ un ouvert
quasicompact la famille $(x_1,\dots,x_n)$ soit linéairement
indépendante sur $F$ dans le $F$-espace vectoriel 
$$
 \GG( \mathfrak{Z}, \underline{\Hom} (F/\O_F,f^*H )\unpi)
$$
On peut alors construire naturellement un morphisme 
$$
\mathfrak{Z}\ldrt \coprod_{a} \DD_{a,\infty}
$$
\end{prop}
\dem 
Soit $\spf (R)\subset \mathfrak{Z}$ un ouvert affine connexe. Soit $H'$ le groupe $p$-divisible sur $\spec (R)$ associé au groupe $p$-divisible $f^* H\times_{\mathfrak{Z}} \spf (R)$. 
D'après le lemme \ref{purifitazm} $\spec (R\unpi)$ est connexe. 
Le $n$-uplet 
$(x_1,\dots,x_n)$ induit donc une quasi-iosgénie 
$$
(F/\O_F)^n \ldrt H'\times_{\spec (R)}\spec (R\unpi)
$$
Il existe donc un unique réseau $\La\subset F^n$ 
(le ``module de Tate'' à l'intérieur du ``module de Tate rationnel'') 
tel que la quasi-isogénie composée 
$$
\delta : \La\otimes F/\La \ldrt (F/\O_F)^n \ldrt H'\times_{\spec (R)}\spec (R\unpi)
$$
soit un isomorphisme. Alors $\delta$ induit un système compatible d'isomorphismes 
$$
\underline{\pi^{-k}\La/\La} \iso H'[\pi^k]\times_{\spec (R)} \spec (R\unpi)
$$
Si $b=[\La_0,M_0]$ posons $a=[\La,M_0]$
 Il y a donc 
d'après la proposition précédente 
un système compatible de morphismes
$$
\xymatrix@C=14mm{
 & \ar@{..>}[d] \\
 & \DD_{a,Id+ \pi^{k+1} \End (\La)} \ar@{..>}[d] \\
\spf (R)  \ar[r] \ar[ru]\ar@{..>}[ruu] &  \DD_{a,Id+\pi^{k} \End (\La)}
}
$$
D'où, puisque $\DD_{a,\infty} =\underset{k}{\limp} \DD_{a,Id +\pi^k \End (\La)}$ dans la catégorie des schémas formels $\pi$-adiques, un morphisme
$$
\spf (R) \ldrt \DD_{a,\infty}
$$
Il est facile de voir que ces différents morphismes se recollent en un
morphisme de $\mathfrak{Z}$ vers $\coprod_{a} \DD_{a,\infty}$.  
\qed

\subsection{Sur la suite de Hodge-Tate en niveau infini}

Comme précédemment pour $G^\vee$ on définit une application de
Hodge-Tate pour $G$ 
$$
\a_G: \O_D \otimes\O_{\mathcal{Y}_\infty} \ldrt \omega_{G^\vee}
\otimes \O_{\mathcal{Y}_\infty}
$$

\begin{prop}\label{krutyopml59}
Dans la suite de Hodge-Tate de $G^\vee$ sur
$\widetilde{\mathcal{Y}}_\infty$
$$
\omega_{G^\vee}^*\otimes \O_{\widetilde{\mathcal{Y}}_\infty} \xrig{\;
  \,^t \a_G \;}
\Hom_{\O_F} (\O_D,\O_F) \otimes \widetilde{\mathcal{Y}}_\infty
\xrig{\; \a_{G^\vee} (-1)\;} \omega_{G} \otimes
\O_{\widetilde{\mathcal{Y}}_\infty }(-1)
$$
on a $\a_{G^\vee} (-1)\circ \,^t\a_{G} =0$. 
\end{prop}
\dem
De la même façon que l'on a une décomposition 
$\a_{G^\vee}=\oplus_{j\in \Z/n\Z} \a_{G^\vee,j}$ on a une
décomposition
$\a_G= \oplus_{j\in \Z/n\Z} \a_{G,j}$. Il suffit alors de montrer que
$\forall j\;\; \a_{G^\vee,j} (-1)\circ \,^t\a_{G,j} =0$. La
démonstration est alors identique à celle de la proposition
\ref{exhodateg}.
\qed

\subsection{Construction d'éléments dans le module de Tate du groupe
  de Lubin-Tate universel tiré en arrière sur
  $\widetilde{\widetilde{\mathcal{Y}}}_\infty$} \label{kqiqeou184}

Notons pour abréger $\widetilde{\widehat{\mathbb{P}}} :=
\widetilde{\widehat{\mathbb{P}}} (\DD (\Hb))$. 
 Rappelons que pour $a=[\La,M]$ $\;\psi_a$
désigne le morphisme
$$ \psi_a :
\DD_a \iso \Pi^j.\mathcal{W}
$$
où $\mathcal{W}\subset \widetilde{\widehat{\mathbb{P}}}$ et $j\equiv
[\O_D:M]\text{ mod } n$. Rappelons que l'on note 
$$
\xi : \widetilde{\widetilde{\mathcal{Y}}}_\infty \ldrt \widetilde{\widehat{\mathbb{P}}}
$$
le morphisme construit dans la section \ref{iqghuqd}. Fixons $j,\, 0\leq
j\leq n-1$ et $a$ comme précédemment associé à $j$.  
On note 
$$
\widetilde{\widetilde{\mathcal{Y}}}_{\infty,j} = \xi^{-1} (\mathcal{W}).\Pi^{-j}
$$
et
$$
\xi_j :  \widetilde{\widetilde{\mathcal{Y}}}_{\infty,j} \xrig{\; \xi
  \;} \Pi^{j}.\mathcal{W} 
$$ 
On note $(H_a,\rho_{H_a})$ la déformation universelle sur $\DD_a$
 ainsi que
$(G,\rho_G)$ celle sur $\mathcal{Y}$. 

Les structures de niveau sur $(\mathcal{Y}_k)_{k\geq 1}$ définissent un
morphisme associé à $1\in \O_D$ 
$$
F/\O_F\ldrt G\times_{\mathcal{Y}} \mathcal{Y}_\infty
$$
 Notons 
$$
\mu : F/\O_F \ldrt G\times_{\mathcal{Y}}  \widetilde{\widetilde{\mathcal{Y}}}_\infty
$$
le morphisme associé sur $
\widetilde{\widetilde{\mathcal{Y}}}_\infty$. Soit maintenant $\kappa
\in \Hom_{\O_D} ( \Gb,\Hb^n)\unpi$. On peut alors considérer le morphisme
composé 
$$
\xymatrix@C=14mm@R=14mm{
F/\O_F\ar[r]^(.37){\mu\text{ mod }\pi}\ar[rrrd]_{\nu} & G\times
(\widetilde{\widetilde{\mathcal{Y}}}_{\infty,j}\text{ mod }\pi)
\ar[r]^(.52){\rho_G^{-1}\times Id } & 
\Gb\times (\widetilde{\widetilde{\mathcal{Y}}}_{\infty,j}\text{ mod }\pi)
\ar[r]^{\kappa\times Id} & \Hb^n \times
(\widetilde{\widetilde{\mathcal{Y}}}_{\infty,j}\text{ mod }\pi) 
\ar[d]^{( \psi^{-1}_a\circ\xi_j)^* \rho_{H_a}^n}\\
&&& (\psi^{-1}_a \circ\xi_j)^* H^n_a \text{ mod }\pi 
}
$$
où $\nu\in \GG  (\widetilde{\widetilde{ \mathcal{Y}}}_{\infty,j}, \underline{\Hom} (F/\O_F,
( \psi^{-1}_a\circ\xi_j)^* H^n_a \text{ mod }\pi )\unpi  )$. 

\begin{theo}
Le morphisme $\nu$ se relève en caractéristique zéro :
$$
\nu \in \GG  (\widetilde{\widetilde{\mathcal{Y}}}_{\infty,j} ,
\underline{\Hom} (F/\O_F, ( \psi^{-1}_a\circ
\xi_j)^* H^n_a)\unpi  )
$$
\end{theo}
\dem
Comme dans la démonstration du théorème \ref{relevisofd24} on applique le critère de
relèvement de Messing.
Considérons le morphisme  composé de $\O$-modules $\pi$-divisibles 
$$
\xymatrix@C=14mm@R=14mm{
F/\O_F\ar[r]^(.37){\mu\text{ mod }\pi} & G\times
(\widetilde{\widetilde{\mathcal{Y}}}_{\infty}\text{ mod }\pi)
\ar[r]^(.52){\rho_G^{-1}\times Id } & 
\Gb\times (\widetilde{\widetilde{\mathcal{Y}}}_{\infty}\text{ mod }\pi)
\ar[r]^{\kappa\times Id} & \Hb^n \times
(\widetilde{\widetilde{\mathcal{Y}}}_{\infty}\text{ mod }\pi) 
}
$$
Le morphisme induit au niveau de l'évaluation des cristaux sur
l'épaississement $(\widetilde{\widetilde{\mathcal{Y}}}_{\infty} \text{
  mod }\pi\hookrightarrow
\widetilde{\widetilde{\mathcal{Y}}}_{\infty})$ induit un élément
$$
1 \longmapsto Y\in \DD(\Hb)^n\otimes \O_{\widetilde{\widetilde{\mathcal{Y}}}_{\infty}}\unpi
$$
Soit $\Fil\subset \DD(\Hb)\otimes
\O_{\widetilde{\widehat{\mathbb{P}}}}$ la filtration définissant
l'espace projectif.   \'Etant donné que $\psi_a$ est un modèle entier
de l'application des périodes on doit vérifier que 
$$
Y\in (\xi^*\Fil\unpi )^n\subset  \DD(\Hb)^n\otimes \O_{\widetilde{\widetilde{\mathcal{Y}}}_{\infty}}\unpi
$$
Rappelons que l'application de Hodge-Tate 
$$
\a_G: \O_D\otimes\O_{\widetilde{\widetilde{\mathcal{Y}}}_\infty} \ldrt
\omega_{G^\vee} \otimes \O_{\widetilde{\widetilde{\mathcal{Y}}}_\infty}
$$
est telle que 
$$
\xymatrix@C=12mm@R=6mm{
 \O_D\otimes\O_{\widetilde{\widetilde{\mathcal{Y}}}_\infty} \unpi
 \ar[r]^{\a_G} & \omega_{G^\vee} \otimes
 \O_{\widetilde{\widetilde{\mathcal{Y}}}_\infty}\unpi \ar@{^(->}[r] &
\DD (\Gb)\otimes\O_{\widetilde{\widetilde{\mathcal{Y}}}_\infty} \unpi
\ar[r]^{\kappa_*\otimes Id} & \DD (\Hb)^n
\otimes\O_{\widetilde{\widetilde{\mathcal{Y}}}_\infty} \unpi \\
1 \otimes 1\ar@{|->}[rrr] &&& Y
}
$$
Écrivons $Y=(Y_i)_{1\leq i\leq n}$, $Y_i\in \DD
(\Hb)\otimes\O_{\widetilde{\widetilde{\mathcal{Y}}}_\infty}
\unpi$. Rappelons qu'on identifie $\DD (\Hb)$ à
$\Hom_{\O_D} ( \O_D,\DD(\Hb))$ (cf. section \ref{guyphng}) et qu'alors d'après
la suite précédente 
$$
Y_i = \,^t \a_G (\gamma_i)
$$ 
où $\,^t\a_G$ désigne $\Hom_{\O_D} ( \a_G,\DD(\Hb))$, le transposé
relativement à la dualité $\Hom_{\O_D} (-,\DD(\Hb))$, et 
$$
\xymatrix{
\gamma_i : \omega_{G^\vee} \otimes
 \O_{\widetilde{\widetilde{\mathcal{Y}}}_\infty}\unpi \ar@{^(->}[r] &
\DD (\Gb)\otimes\O_{\widetilde{\widetilde{\mathcal{Y}}}_\infty} \unpi
\ar[r]^{\kappa_*\otimes Id} & \DD (\Hb)^n
\otimes\O_{\widetilde{\widetilde{\mathcal{Y}}}_\infty} \unpi
\ar@{->>}[r]^{pr_i}  & \DD (\Hb)
\otimes\O_{\widetilde{\widetilde{\mathcal{Y}}}_\infty} \unpi
}
$$
Donc d'après la remarque \ref{qspfihb794}
si $\delta_i\in \omega_{G^\vee,n-1}^*\otimes \DD (\Hb)_{n-1}=
\Hom_{\breve{\O}} ( \omega_{G^\vee},\DD (\Hb)_{n-1})\unpi$ est la forme
linéaire
$$
\xymatrix{
\delta_i : \omega_{G^\vee,n-1} \otimes
 \O_{\widetilde{\widetilde{\mathcal{Y}}}_\infty}\unpi \ar@{^(->}[r] &
\DD (\Gb)_{n-1}\otimes\O_{\widetilde{\widetilde{\mathcal{Y}}}_\infty} \unpi
\ar[r]^{\kappa_*\otimes Id} & \DD (\Hb)^n_{n-1}
\otimes\O_{\widetilde{\widetilde{\mathcal{Y}}}_\infty} \unpi
\ar[r]^{pr_i}  & \DD (\Hb)_{n-1}
\otimes\O_{\widetilde{\widetilde{\mathcal{Y}}}_\infty} \unpi
}
$$
on a avec les notations de la section \ref{qsqfugv354}
$$
Y_i = ((\,^t \a_{G,n-1})\otimes \DD ( \Hb)_{n-1}) ( \delta_i)
$$
Le théorème est donc une conséquence de la proposition \ref{krutyopml59}.
\qed

\begin{prop}\label{kafeterui}
La construction précédente définit un morphisme 
$$
\Hom_{\O_D} ( \Gb,\Hb^n)\unpi \ldrt \GG  (\widetilde{\widetilde{\mathcal{Y}}}_{\infty,j} ,
\underline{\Hom} (F/\O_F, ( \psi^{-1}_a\circ
\xi_j)^* H_a)\unpi  )^n
$$
Si $\nu = (\nu_1,\dots,\nu_n)$ correspond à $\Delta^{-1}\in
\Hom_{\O_D} ( \Gb,\Hb^n)\unpi$ alors pour tout ouvert quasicompact
$\mathcal{U}$ de $\widetilde{\widetilde{\mathcal{Y}}}_{\infty,j}$ 
les éléments $(\nu_{1|\mathcal{U}},\dots, \nu_{n|\mathcal{U}})$ sont
linéairement indépendants dans le $F$-espace vectoriel 
$\Hom (F/\O_F, (\psi_a^{-1}\circ \xi_{j|\mathcal{U}})^*H_a)\unpi$. 
\end{prop}
\dem
Pour démontrer l'indépendance linéaire on procède par spécialisation
afin de se ramener au cas d'un point. 
Soit donc $\mathcal{U}$ comme dans l'énoncé. On peut supposer $\mathcal{U}
=\spf (R)$. La $\breve{\O}$-algèbre $R$ étant $\pi$-adique sans
$\pi$-torsion 
il existe un corps valué complet $K$ pour une valuation de
rang $1$ étendant celle de $\breve{F}$ et un morphisme continu $R\ldrt
\O_K$
(utiliser le théorème 1.2.1 p.13 de \cite{BerkSpectral} appliqué à
l'algèbre de Banach $R\unpi$). 
 Soit $\theta : \spf (\O_K)\ldrt \spf (R) \hookrightarrow
\widetilde{\widetilde{\mathcal{Y}}}_{\infty,j}$. 
Il suffit alors de montrer l'indépendance linéaire de
$(\theta^*\nu_1,\dots,\theta^*\nu_n)$. Mais cela est démontré à la fin
de la section 9.1 de \cite{Points}.
\qed

\subsection{Construction du morphisme de
  $\widetilde{\widetilde{\mathcal{Y}}}_\infty$ vers $\X_\infty$}

Comme dans la section \ref{okk65tr} on montre la proposition suivante.

\begin{prop}
Soit $G$ le $\O_D$-module formel universel sur
$\widehat{\Omega}$. Soit $C\subset \O_D^\times$  un sous-groupe
ouvert, 
$U\subset \Omega_C$ un ouvert admissible non-vide et $\bar{s}\in U
(\overline{\breve{F}})$. Alors le sous-groupe de monodromie
arithmétique image de 
$$
\pi_1 (U,\bar{s}) \ldrt \Aut_{\O_D} ( T_p (G)) \simeq {\O_D^{opp}}^\times
$$
est un sous-groupe ouvert de $ {\O_D^{opp}}^\times$. 
 Donc, si $\Pi_{C',C}:\Omega_{C'}\ldrt \Omega_C$ l'ensemble 
$\underset{C'\subset C}{\limp} \pi_0 ( \Pi_{C',C}^{-1} (U))$ est fini.
\end{prop}

\begin{coro}
Tout ouvert quasicompact de
$\widetilde{\widetilde{\mathcal{Y}}}_\infty$ possède un nombre fini de
composantes connexes et de plus $\O_{
  \widetilde{\widetilde{\mathcal{Y}}}_\infty}$ est intégralement fermé
dans $\O_{
  \widetilde{\widetilde{\mathcal{Y}}}_\infty} \unpi$. 
\end{coro}

Soit $j, \, 0\leq j\leq n-1$ et $b$ un sommet de l'immeuble associé à $j$ c'est à dire tel que 
si $b=[\La,\Pi^k\O_D]$ alors $k\equiv j\text{ mod }n$. D'après le corollaire précédent, 
la proposition \ref{kafeterui} et la proposition \ref{ghotopolmk} on peut construire un morphisme
$$
\omega_j : \widetilde{\widetilde{\mathcal{Y}}}_{\infty,j}\ldrt \coprod_a \DD_{\infty,a}
$$
dont on vérifie aussitôt qu'il ne dépend pas du choix fait de $b$. 
 Il est également aisé de vérifier que ce morphisme est $\GL_n (F)$-équivariant où l'action à gauche de $\GL_n (F)$ sur $ \widetilde{\widetilde{\mathcal{Y}}}_{\infty,j}$ est transformée en celle à droite sur $\coprod_a \DD_{\infty,a}$ via $g\mapsto \,^t g$.

\begin{prop}
Pour $j_1,j_2$
le diagramme suivant est commutatif
$$
\xymatrix{
 & \coprod_a \DD_{\infty,a} \ar[rd] \\
 \widetilde{\widetilde{\mathcal{Y}}}_{\infty,j_1}
\cap 
\widetilde{\widetilde{\mathcal{Y}}}_{\infty,j_2} \ar[ru]^{\omega_{j_1}}
\ar[rd]_{\omega_{j_2}} & & \X_\infty \\
& \coprod_a \DD_{\infty,a} \ar[ru]
}
$$
\end{prop}
\dem
Soient $0\leq j_1<j_2\leq n-1$
Avec les notations de \cite{Cellulaire} on a pour $b_1$ associé à $j_1$ et $b_2$ associé à $j_2$
$$
\xymatrix@R=6mm{
\partial_{j_2-j_1}\DD_{b_1} \ar[rd]^{\psi_{b_1}}_\simeq \\
&\Pi^{j_1}.\mathcal{W} \cap \Pi^{j_2}.\mathcal{W} \\
\partial_{j_2-j_1}\DD_{b_2} \ar[ru]_{\psi_{b_2}}^\simeq
}
$$
il est aisé d'en déduire qu'il y a des factorisations
$$
\xymatrix@C=12mm{
 \widetilde{\widetilde{\mathcal{Y}}}_{\infty,j_1}\cap \widetilde{\widetilde{\mathcal{Y}}}_{\infty,j_2}\ar[r] \ar@(ru,lu)[rr]^{\omega_{j_1}} \ar@(rd,ld)[rr]_{\omega_{j_2}}
 & \dpt{\coprod_{a\drt a'}} \DD_{\infty,a\drt a'} \ar@<.8ex>[r] \ar@<-.8ex>[r] & \dpt{\coprod_{a}} \DD_{\infty,a} 
}
$$
\qed
\\

Il y a donc un morphisme  $\GL_n (F)\times D^\times$-équivariant
$$
\widetilde{\widetilde{\mathcal{Y}}}_\infty \ldrt \X_\infty
$$
 où l'action de $\GL_n ( F)$ est transformée via $g\mapsto\,^t g$ et celle de $D^\times$ via $d\mapsto d^{-1}$. 

\section{Construction de l'isomorphisme}\label{consteug27}

Nous utiliserons maintenant les résultats de l'appendice afin de ne
pas avoir à repasser en niveau fini à chaque fois que l'on veut
éclater un idéal, ce qui deviendrait rapidement inextricable comme
l'auteur a pu le vérifier en tentant de rédiger cette section sans les
résultats de l'appendice. 

\subsection{De nouveaux éclatements}

On a construit dans les section précédentes un diagramme $\GL_n (F)\times D^\times$-équivariant  de morphismes de schémas formels $\pi$-adiques sur $\spf (\breve{\O})$ 
$$
\xymatrix{
\widetilde{\X}_\infty  \ar[rd]\ar[d] & \widetilde{\widetilde{\mathcal{Y}}}_\infty \ar[d]\ar[ld] \\
\X_\infty & \mathcal{Y}_\infty
}
$$
Rappelons que $\widetilde{\X}_\infty$ n'est pas vraiment un éclatement
formel admissible d'un idéal de $\O_{\X_\infty}$ mais plutôt d'une
famille d'idéaux. Plus précisément 
$$
\X_\infty = \bigcup_a \DD_{a,\infty}
$$
et pour tout $a$ il y a un idéal admissible $\mathfrak{I}_a \subset
\O_{\DD_{a,\infty}}$ tel que 
\begin{itemize}
\item $\forall (g,d)\in \GL_n (F)\times D^\times\;\; (g,d)^*
  \mathfrak{I}_{(g,d).a} = \mathfrak{I}_a$
\item $\forall a,a'\;\; \mathfrak{I}_{a'|\DD_{a',\infty}\cap
    \DD_{a,\infty}}$ devient localement libre de rang un sur l'éclaté
  formel admissible de $\mathfrak{I}_a$ dans $\DD_{a',\infty}\cap
  \DD_{a,\infty}$ (mais en général
  $\mathfrak{I}_{a'|\DD_{a,\infty}\cap \DD_{a,\infty}}\neq
  \mathfrak{I}_{a|\DD_{a,\infty}\cap \DD_{a,\infty}}$)
\item $\widetilde{\DD}_{a,\infty}$ est le normalisé dans sa fibre
  générique de l'éclatement formel admissible de $\mathfrak{I}_a$
\end{itemize}

Notons 
$$
f : \widetilde{\widetilde{\mathcal{Y}}}_\infty \ldrt \X_\infty
$$
Notons pour tout $a$
$$
\mathcal{V}_a = f^{-1} ( \DD_{a,\infty}) \;\;\text{ et }\;\;
\mathfrak{J}_a = \O_{\mathcal{V}_a}. f^{-1} \mathfrak{I}_a 
$$
un idéal admissible de $\O_{\mathcal{V}_a}$. Soit
$\widetilde{\mathcal{V}}_a$ le normalisé dans sa fibre générique de l'éclatement formel
admissible de $\mathfrak{J}_a$. Pour tout couple $a,a'$ l'idéal 
$\mathfrak{J}_{a'|\mathcal{V}_{a'}\cap \mathcal{V}_a}$ devient
localement libre de rang un sur $\widetilde{\mathcal{V}}_a$. Les
$\widetilde{\mathcal{V}}_a$ se recollent donc (utiliser la propriété
universelle des normalisations dans la fibre générique/éclatements
formels admissibles, cf. appendice) en un schéma formel 
$$
\widetilde{\mathcal{Y}}^{(3)}_\infty = \bigcup_{a} \widetilde{\mathcal{V}}_a
$$
qui est $\GL_n (F)\times D^\times$-équivariant. De plus il y a un
morphisme équivariant
$$
\widetilde{\mathcal{Y}}^{(3)}_\infty \ldrt \widetilde{\X}_\infty
$$
Dans l'autre sens soit 
$$
h: \widetilde{\X}_\infty \ldrt \mathcal{Y}_\infty
$$
Il existe un idéal admissible $\GL_n (F)\times D^\times$-invariant
 $\mathfrak{J}\subset \O_{\mathcal{Y}_\infty}$ tel que $\widetilde{\widetilde{\mathcal{Y}}}_\infty$
soit le normalisé dans sa fibre générique de l'éclatement formel
admissible de $\mathfrak{J}$ dans $ \mathcal{Y}_\infty$. Soit
$\widetilde{\X}^{(2)}_\infty$ le normalisé dans sa fibre générique de
l'idéal admissible $\O_{\widetilde{\X}_\infty}.f^{-1}
\mathfrak{J}$. Il est muni d'une action de $\GL_n (F)\times D^\times$
et il y a de plus un morphisme équivariant
$$
\widetilde{\X}_\infty^{(2)} \ldrt \widetilde{\widetilde{\mathcal{Y}}}_\infty
$$
On a donc un nouveau diagramme au dessus du précédent
$$
\xymatrix{
\widetilde{\X}_\infty^{(2)} \ar[rd]\ar[d] &
\widetilde{\mathcal{Y}}_\infty^{(3)} \ar[d] \ar[ld] \\
\widetilde{\X}_\infty & \widetilde{\widetilde{\mathcal{Y}}}_\infty
}
$$
Mais si $h': \widetilde{\X}_\infty^{(2)} \ldrt
\widetilde{\mathcal{Y}}_\infty$ on a un diagramme pour tout $a$
$$
\xymatrix{
h'^{-1} ( \mathcal{V}_a) \ar[rd]^{h'} \ar[d] \\
\widetilde{\DD}_{a,\infty} \ar[d] & \mathcal{V}_a \ar[ld]^{f} \\
\DD_{a,\infty}
}
$$
et donc l'image réciproque à $h'{-1} ( \mathcal{V}_a)$ par $h$' de
l'idéal $\mathfrak{I}_a$, $\O_{h'^{-1} (\mathcal{V}_a)}. h'^{-1}
\mathfrak{I}_a$ 
 est localement libre de rang un sur $h'^{-1}
(\mathcal{V}_a)$ puisque c'est le cas sur $\widetilde{\DD}_{a,\infty}$
(on utilise toujours le fait que sur un schéma formel $\pi$-adique sans $\pi$-torsion
un idéal admissible est localement
libre de rang un ssi il est localement monogène). Donc d'après la
propriété universelle du ``normalisé dans la fibre générique/éclaté'' le
morphisme $h'$ se relève à $\widetilde{\mathcal{Y}}_\infty^{(3)}$. 
De la même façon on relève le morphisme
$\widetilde{\mathcal{Y}}_\infty^{(3)}\ldrt 
\widetilde{\X}_\infty$ à $\widetilde{\X}_\infty^{(2)}$. 

On obtient donc au final un diagramme
$$
\xymatrix{
\widetilde{\X}_\infty^{(2)} \ar[rd]\ar[d] \ar@<.6ex>[r] & \ar@<.6ex>[l]
\widetilde{\mathcal{Y}}_\infty^{(3)} \ar[d] \ar[ld] \\
\widetilde{\X}_\infty  \ar[rd]\ar[d] & \widetilde{\widetilde{\mathcal{Y}}}_\infty \ar[d]\ar[ld] \\
\X_\infty & \mathcal{Y}_\infty
}
$$

Le but des sections qui suivent est maintenant de démontrer le
théorème suivant.

\begin{theo}
Les deux morphismes $\xymatrix{
\widetilde{\X}_\infty^{(2)} \ar@<.6ex>[r] & \ar@<.6ex>[l]
\widetilde{\mathcal{Y}}_\infty^{(3)}}$ sont inverses l'un de l'autre
et fournissent donc un isomorphisme équivariant entre
$\widetilde{\X}_\infty^{(2)}$ et $\widetilde{\mathcal{Y}}_\infty^{(3)}$.
\end{theo}

\subsection{Retour aux suites de Hodge-Tate en niveau infini}\label{ivhjubg1358}

Nous allons améliorer les résultats précédents d'exactitude des suites
de Hodge-Tate en niveau infini.

\begin{lemm}\label{mfcclat75}
Soit $(X,\O_X)$ un espace annelé et $u:\E_1\ldrt \E_2$ un morphisme de
$\O_X$-modules localement libres de rang $r$. Soit $t\in \GG
(X,\O_X)$.
\begin{itemize}
\item Si $\text{coker} \;u$ est annulé par $t$ alors $\text{coker} \;(\det u )$ est
  annulé par $t^r$
\item Si $\text{coker} \;(\det u)$ est annulé par $t$ alors
  $\text{coker} \;u$ est
  annulé par $t$ et localement sur $X$ il existe un morphisme $v$
$$
\xymatrix{
\E_1 \ar[r]_u & \E_2 \ar@{-->}@(ul,dl)[l]_(.22){v} 
}
$$
tel que $u\circ v = t $. 
\end{itemize}
\end{lemm}

\begin{theo}\label{lvjozr48}
Soit $K|F$ un corps valué complet pour une valuation de rang $1$
étendant celle de $F$. Soit $H_0$ un $\O$-module  $\pi$-divisible sur
$\O_K$. Dans la suite de Hodge-Tate
$$
\omega_{H_0}^*\otimes \O_{\widehat{\overline{K}}} (1) \xrig{\; \,^t
  \a_{H_0^\vee} (1)\;}
T_p (H)\otimes  \O_{\widehat{\overline{K}}} \xrig{\; \a_{H_0}\;} \omega_{H_0^\vee}
\otimes \O_{\widehat{\overline{K}}} 
$$
on a $\ker \a_{H_0}/\text{Im}  \,^t
  \a_{H_0^\vee} (1)$ est annulé par $\pi^{n-1}$.
\end{theo}
\dem
On renvoie à \cite{Points}.
\qed

On reprend maintenant les notations de la section \ref{decosdza}.

\begin{prop}
Soit $a=[\La,M]$ un sommet de l'immeuble paramétrant les cellules de
$\X_\infty$ et $\widetilde{\DD}_\infty := \widetilde{\DD}_{a,\infty}$
la cellule éclatée en niveau infini. Soit la suite de Hodge-Tate sur
$\widetilde{\DD}_\infty$
$$
\omega_{H^\vee}^*\otimes\O_{\widetilde{\DD}_\infty} \xrig{\; \,^t
  \a_H\;} \La^*\otimes \O_{\widetilde{\DD}_\infty} \xrig{\;\a_{H^\vee}
  (-1)\;}  \omega_H\otimes \O_{\widetilde{\DD}_\infty} (-1)
$$
Alors $\ker \, ( \,^t
  \a_H )/\text{Im}\, ( \a_{H^\vee}
  (-1))$ est annulé par $\pi^{n-1}$. 
\\
De plus, localement
  sur $\widetilde{\DD}_\infty$ il y a un morphisme $\beta$
$$
\xymatrix@C=18mm{
\omega_{H^\vee}^*\otimes\O_{\widetilde{\DD}_\infty} 
 \ar[r]_{\,^t
  \a_{H}} & \ar@{..>}@(ul,ur)[l]_{\beta}  \La^*\otimes \O_{\widetilde{\DD}_\infty} 
}
$$
tel que $\,^t \a_H\circ \beta = \pi^{n-1}$ et $  \beta \circ \,^t \a_H
=\pi^{n-1}$. 
\end{prop}
\dem
On sait, par construction de $\widetilde{\DD}_\infty$ que $\text{Im}\, ( \a_{H^\vee}
  (-1))$ est localement libre de rang un. On en déduit aisément que
  $\ker \, (\a_{H^\vee}
  (-1))$ est localement libre de rang $n-1$. D'après  la proposition
  \ref{exhodateg} il y a un morphisme
$$
u: \omega_{H^\vee}^*\otimes\O_{\widetilde{\DD}_\infty} \xrig{\; \,^t
  \a_H\;} \ker \, (\a_{H^\vee}
  (-1))
$$
Intéressons-nous à $\det\, u =\wedge^{n-1} u$. Localement sur un
ouvert quasicompact $\mathcal{U}$ de $\widetilde{\DD}_\infty$ $\; \det \, u$ est
donné par une fonction $f\in
\GG(\mathcal{U},\O_{\widetilde{\DD}_\infty})$. Soit  $x\in
\mathcal{U} ( \O_{\C_p}),\, x: \spf (\O_{\C_p})\ldrt \mathcal{U}$.
 \'Etant donné que la suite exacte 
$$
0 \ldrt  \ker \, (\a_{H^\vee} )\ldrt  \La^*\otimes \O_{\widetilde{\DD}_\infty} \xrig{\;\a_{H^\vee}\;}
 \text{Im } ( \a_{H^\vee} (-1) )
  (-1)) 
$$
est localement scindée on en déduit que $\ker ( \a_{x^* H^\vee} (-1))
= x^* \ker \, (\a_{H^\vee} )$ et donc d'après le théorème
\ref{lvjozr48} couplé au premier point du lemme \ref{mfcclat75}
$$
|f(x)|\geq |\pi^{n-1}|
$$
Du lemme qui suit on en déduit que $\pi^{n-1} \in \O_\mathcal{U} .f$ et
on conclut grâce au second point du lemme \ref{mfcclat75}.
\qed

\begin{lemm}
Soit $\mathcal{U}$ un ouvert quasicompact de $\widetilde{\DD}_\infty$, $\a\in\N$ et $f\in \GG
(\mathcal{U},\O_{\widetilde{\DD}_\infty})$ tel que $\forall x\in
\mathcal{U} ( \O_{\C_p})$ $\; |f(x)|\geq |\pi^\a|$. Alors $\pi^\a \in
\O_\mathcal{U}. f$. 
\end{lemm}
\dem
L'ouvert $\mathcal{U}$ étant quasicompact il existe $k\geq 3$ ainsi
qu'un ouvert $\mathcal{V}\subset \widetilde{\DD}_k$ et $g\in \GG (
\mathcal{V}, \O_{ \widetilde{\DD}_k})$ tels que l'image réciproque de
$\mathcal{V}$ en niveau infini soit $\mathcal{U}$ et $f\equiv g \text{
  mod } \pi^{\a+1}$. L'application $\mathcal{U} ( \O_{\C_p}) \ldrt
\mathcal{V} (\O_{\C_p})$ est surjective (on a $ \mathcal{V}
(\O_{\C_p}) = \mathcal{V}^{rig} ( \Cp)$ et, revenant à
l'interprétation modulaire de $\DD_{\infty}^{rig}$, on peut toujours
prolonger une structure de niveau finie en une infinie sur $\C_p$). 
Par hypothèse  on a donc
$$
\forall x\in \mathcal{V} (\O_{\C_p}) \; |g(x)|\geq |\pi^{\a}|
$$
Mais $\mathcal{V}$ étant normal on en déduit que $\pi^\a \in (g)$ et
donc $\pi^{\a} \in (f,\pi^{\a+1})$ ce qui implique $\pi^\a \in (f)$
 puisque
$\widetilde{\DD}_\infty$ est $\pi$-adique.
\qed
\\

\begin{coro}\label{dgjigy124}
La suite 
$$ 0\ldrt 
\omega_{H^\vee}^*\otimes\O_{\widetilde{\DD}_\infty}\unpi \xrig{\; \,^t
  \a_H\;} \La^*\otimes \O_{\widetilde{\DD}_\infty} \unpi \xrig{\;\a_{H^\vee}
  (-1)\;}  \omega_H\otimes \O_{\widetilde{\DD}_\infty} \unpi (-1)
\ldrt 0
$$
est exacte localement scindée. De plus le conoyau du morphisme 
$$
\a_H : \La\otimes\O_{\widetilde{\DD}_\infty} \ldrt \omega_{H^\vee}
\otimes \O_{\widetilde{\DD}_\infty}
$$
est annulé par $\pi^{n-1}$. 
\end{coro}

De la même façon on montre avec les notations de la proposition
\ref{krutyopml59} la proposition qui suit.

\begin{prop}\label{dvjlee138}
Dans la suite de Hodge-Tate de $G^\vee$ sur
$\widetilde{\mathcal{Y}}_\infty$
$$
\omega_{G^\vee}^*\otimes \O_{\widetilde{\mathcal{Y}}_\infty} \xrig{\;
  \,^t \a_G \;}
\Hom_{\O_F} (\O_D,\O_F) \otimes \O_{\widetilde{\mathcal{Y}}_\infty}
\xrig{\; \a_{G^\vee} (-1)\;} \omega_{G} \otimes
\O_{\widetilde{\mathcal{Y}}_\infty }(-1)
$$
on a $\ker (\a_{G^\vee} (-1)) / \text{Im} (\,^t\a_{G})$ est annulé par
$\pi^{n-1}$ et il existe un morphisme $\beta :  \Hom_{\O_F}
(\O_D,\O_F) \otimes \widetilde{\mathcal{Y}}_\infty \ldrt
\omega_{G^\vee}^*\otimes \O_{\widetilde{\mathcal{Y}}_\infty}$ tel que 
$\,^t \a_G \circ \beta = \pi^{n-1}$ et $\beta \circ \,^t \a_G
=\pi^{n-1}$. La suite
$$
\omega_{G^\vee}^*\otimes \O_{\widetilde{\mathcal{Y}}_\infty} \unpi \xrig{\;
  \,^t \a_G \;}
\Hom_{\O_F} (\O_D,\O_F) \otimes \O_{\widetilde{\mathcal{Y}}_\infty} \unpi
\xrig{\; \a_{G^\vee} (-1)\;} \omega_{G} \otimes
\O_{\widetilde{\mathcal{Y}}_\infty }\unpi (-1)
$$ 
est exacte localement scindée. 
\end{prop}

\subsection{Démonstration du théorème principal}
\subsubsection{La composée $\widetilde{\X}_\infty^{(2)} \ldrt \widetilde{\mathcal{Y}}_\infty^{(3)} \ldrt \widetilde{\X}_\infty^{(2)}$ est l'identité}

Nous utiliserons constamment les deux faits suivants :
\begin{itemize}
\item Soient $\mathfrak{Z}$ et $\X$ deux schémas formels $\pi$-adiques sans $\pi$-torsion. Soit $\widetilde{\mathfrak{Z}}\xrig{\; p\;} \mathfrak{Z}$ un éclatement formel admissible. Soient $\xymatrix{\X\ar@<.6ex>[r]^{f_1} \ar@<-.6ex>[r]_{f_2} & \widetilde{\mathfrak{Z}}}$ deux morphismes  tels que $p\circ f_1 = p\circ f_2$. Alors $f_1=f_2$. 
\item Soit $\X$ $\pi$-adique sans $\pi$-torsion. Deux morphismes $\xymatrix{
\X\ar@<.6ex>[r]  \ar@<-.6ex>[r] & \widehat{\mathbb{P}}^n}$ sont égaux
ssi les morphismes associés le sont après inversion de $\pi$ le sont. Plus précisément si 
$$
\xymatrix@R=4mm{
 & \L_1 \\
\O_\X^n \ar@{->>}[rd]\ar@{->>}[ru] \\
& \L_2
}
$$
correspond à nos deux morphismes alors ceux-ci sont égaux ssi $\exists u$
$$
\xymatrix@R=4mm{
 & \L_1 \unpi \ar[dd]^{\simeq}_{u} \\
(\O_\X\unpi )^n \ar@{->>}[rd]\ar@{->>}[ru] \\
& \L_2\unpi
}
$$

\end{itemize}

Soit $a=[\La,M]$ un sommet de l'immeuble paramétrant les cellules de $\X_\infty$. Soit $\widetilde{\DD}_\infty^{(2)} : = \DD_{a,\infty}^{(2)}$ l'image réciproque de 
$\widetilde{\DD}_{a,\infty}$ dans $\widetilde{\X}_\infty^{(2)}$. Commençons par montrer que le diagramme suivant est commutatif.
$$
\xymatrix{
\widetilde{\DD}_\infty^{(2)} \ar@{^(->}[r]\ar[d] & \widetilde{\X}^{(2)}_\infty \ar[r] & \widetilde{\mathcal{Y}}^{(3)} \ar[d] \\
\DD_\infty \ar[d] & & \widetilde{\widetilde{\mathcal{Y}}}_\infty \ar[lld] \\
\widetilde{\widehat{\mathbb{P}}} (\DD (\Hb))
}
$$
où le morphisme $\DD_\infty \ldrt \widetilde{\widehat{\mathbb{P}}} (\DD (\Hb))$ est celui noté $\psi_a$ dans la section \ref{qdfkobji1247}.
Il suffit de montrer que le diagramme qui suit est commutatif.
$$ 
\xymatrix{
\widetilde{\DD}_\infty^{(2)} \ar[r]\ar[dd] & \widetilde{\widetilde{\mathcal{Y}}}_\infty \ar[d] \\
& \widetilde{\mathcal{Y}}_\infty \ar[ld] \\
\widehat{\mathbb{P}} ( \DD (\Hb))
}
$$
Avec les notations des sections précédentes le morphisme composé 
$$
\widetilde{\DD}_\infty^{(2)} \ldrt
\widetilde{\widetilde{\mathcal{Y}}}_\infty \ldrt
\widetilde{\mathcal{Y}}_\infty \ldrt 
\widehat{\mathbb{P}} ( \DD (\Hb))
$$
est obtenu à partir de l'application de Hodge-Tate de $G$ tiré en
arrière sur $\widetilde{\DD}_\infty^{(2)}$ 
$$
\a_{G^\vee,n-1}\otimes Id :\underbrace{ \Hom_{\O_F} ( \O_D,
\breve{\O})_{n-1} \otimes
\DD(\Hb)_{n-1}}_{\DD(\Hb)}   \otimes \O_{\widetilde{\DD}_\infty^{(2)} }\unpi
\ldrt \omega_{G,n-1}\otimes \DD(\Hb)_{n-1} \otimes \O_{\widetilde{\DD}_\infty^{(2)} }\unpi
$$
Mais d'après la remarque \ref{qspfihb794} et la proposition
\ref{dvjlee138} le noyau de l'application précédente est obtenu en
appliquant $\Hom_{\O_D} ( -,\DD(\Hb))$ à 
$$
\a_G : \O_D\otimes \O_{\widetilde{\DD}_\infty^{(2)}} \unpi \ldrt
\omega_{G^\vee}\otimes \O_{\widetilde{\DD}_\infty^{(2)}}\unpi 
$$
via $\Hom_{\O_D} ( \O_D ,\DD(\Hb))\simeq \DD(\Hb)$. Mais si $H$
désigne le groupe de Lubin-Tate universel sur
$\widetilde{\DD}_\infty^{(2)}$ l'image de la composée
$$
\xymatrix@R=3mm{
\O_D\otimes \O_{\widetilde{\DD}_\infty^{(2)}}\unpi \ar[r]^{\a_G} &
\omega_{G^\vee} \otimes \O_{\widetilde{\DD}_\infty^{(2)}}\unpi
\ar@{^(->}[r]^{\rho_{G*}^{-1}} & \DD(\Gb)\otimes
\O_{\widetilde{\DD}_\infty^{(2)}}\unpi
\ar[r]_(.5){\simeq}^(.5){\Delta_*^{-1}} &
\DD(\Hb)^n \otimes \O_{\widetilde{\DD}_\infty^{(2)}}\unpi  \\
1\otimes 1 \ar[rrr] &&& (x_1,\dots,x_n)
}
$$
est telle que si 
$$
\a_H : \La \otimes \O_{\widetilde{\DD}_\infty^{(2)}}\unpi =
\O_{\widetilde{\DD}_\infty^{(2)}}\unpi^n \ldrt \omega_{H^\vee} \otimes
\O_{\widetilde{\DD}_\infty^{(2)}}\unpi
\overset{\rho_{H*}^{-1}}{\hookrightarrow}  \DD(\Hb) \otimes \O_{\widetilde{\DD}_\infty^{(2)}}\unpi
$$
alors $\text{Im } \a_H = \O_{\widetilde{\DD}_\infty^{(2)}}\unpi\, x_1+
\dots + \O_{\widetilde{\DD}_\infty^{(2)}}\unpi \, x_n$. Cela découle
en effet de ce que modulo $\pi$ les éléments du module de Tate de $H$
et $G$ sont reliés par (en notations abrégées)
$$
F/\O_F \ldrt H^n \text{ mod }\pi \xrig{\;\rho_H^{-1}\; } \Hb^n
\xrig{\; \Delta\;} \Gb \xrig{\;\rho_G\;} G \text{ mod }\pi
$$
où $F/O_F \ldrt  H^n \text{ mod }\pi$ fournit l'identification du
module de Tate de $H$ avec $F^n$ et la composée $F/\O_F\ldrt G\text{
  mod }\pi$ fournit l'élément $1\in \O_D$ associé à l'identification
du module de Tate de $G$ avec $D$. 
\\
Appliquant maintenant $\Hom_{\O_D} ( -,\DD(\Hb))$ à 
$$
\O_D\otimes  \O_{\widetilde{\DD}_\infty^{(2)}}\unpi \ldrt
\DD(\Hb)^n\otimes  \O_{\widetilde{\DD}_\infty^{(2)}}\unpi
$$
on en déduit aisément que la filtration de $\DD (\Hb)$ associé à la
composée 
$$
\widetilde{\DD}_\infty^{(2)} \ldrt
\widetilde{\widetilde{\mathcal{Y}}}_\infty \ldrt
\widetilde{\mathcal{Y}}_\infty \ldrt 
\widehat{\mathbb{P}} ( \DD (\Hb))
$$
est donnée par $\text{Im } (\a_H)$. Mais d'après le corollaire
\ref{dgjigy124} $\text{Im } (\a_H)$ coïncide avec $\omega_{G^\vee}
\otimes  \O_{\widetilde{\DD}_\infty^{(2)}}\unpi$ qui définit la
filtration associée au morphisme $\widetilde{\DD}_\infty^{(2)}\ldrt
\DD \ldrt \widehat{\mathbb{P}} ( \DD (\Hb))$. On a donc bien montré que
le diagramme voulu de morphismes vers $\widetilde{\widehat{\mathbb{P}}} ( \DD
(\Hb))$ est commutatif.

Reste à voir que le morphisme $\DD_\infty^{(2)} \ldrt \widetilde{\widetilde{\mathcal{Y}}}_\infty \ldrt
\widetilde{\widehat{\mathbb{P}}} ( \DD (\Hb))$ relevé à $\X_\infty$
grâce au éléments du module de Tate de $H$ construits dans la section
\ref{lqvout159} coïncide avec la projection
$\widetilde{\DD}_\infty^{(2)} \ldrt \DD_\infty \subset
\X_\infty$. Mais c'est une conséquence immédiate de ce que ces
éléments du module de Tate coïncident modulo $\pi$ puisque modulo
$\pi$ 
les
constructions des sections \ref{mvhzaf938} et \ref{kqiqeou184} sont
inverses l'une de l'autre.

\subsubsection{La composée $\widetilde{\mathcal{Y}}^{(3)} \ldrt
  \widetilde{\X}_\infty^{(2)} \ldrt \widetilde{\mathcal{Y}}^{(3)}$ est
  l'identité}

La méthode de démonstration est identique à celle de la section
précédente. On démontre grâce aux résultats de la section \ref{ivhjubg1358}
 que deux morphismes
$\xymatrix{\widetilde{\mathcal{Y}}^{(3)}_\infty \ar@<.6ex>[r] \ar@<-.6ex>[r] &
  \widehat{\Omega}}$ coïncident. Puis comme précédemment on constate
que les relèvements associés vers $\mathcal{Y}_\infty$
coïncident.

\appendix

\section{Compléments sur les schémas formels $\pi$-adiques}

On fixe $\breve{\O}$ un anneau de valuation de hauteur $1$ et $\pi$ un élément de $\breve{\O}$ de valuation strictement positive. On note $k$ le corps résiduel de $\breve{\O}$. 
Tous les schémas formels considérés seront supposés quasi-séparés.

\subsection{Quelques lemmes d'algèbre $\pi$-adique} 

Si $R$ est une $\breve{\O}$-algèbre on note $\widehat{R}$ le complété $\pi$-adique de $R$.
Si $R$ est une $\breve{\O}$-algèbre sans $\pi$-torsion on note $\overline{R}$ la fermeture
intégrale de $R$ dans $R\unpi$.

Nous utiliserons constamment le lemme suivant dont la démonstration ne pose pas de problème.

\begin{lemm}\label{carotidegh}
Soit $R$ une $\breve{\O}$-algèbre sans $\pi$-torsion. Alors,
\begin{itemize}
\item $\widehat{R}$ est sans $\pi$-torsion
\item Soit $\a: R\unpi \ldrt \widehat{R}\unpi$. Alors $\overline{R} =
  \a^{-1} \left (\, \overline{\widehat{R}} \,\right)$. 
\item si $R$ est intégralement fermé dans $R\unpi$ alors $\widehat{R}$ est intégralement fermé dans $\widehat{R}\unpi$.
\item si $R$ s'écrit comme une limite inductive filtrante 
$R=\underset{i}{\limi} R_i$ où les $(R_i)_i$ sont des $\breve{ \O}$-algèbre sans $\pi$-torsion telle que $\forall i\; R_i$ soit intégralement fermé dans $R_i\unpi$ alors $R$ est intégralement fermé dans $R\unpi$.
\end{itemize}
\end{lemm}

\begin{coro}\label{mutk256}
Soit $R$ une $\breve{\O}$-algèbre sans $\pi$-torsion. Alors, canoniquement, 
$$
\widehat{\overline{R}}  = \widehat{\overline{\widehat{R}}}
$$
\end{coro}
\dem
Le morphisme $R\ldrt \widehat{R}$ induit par application de $\widehat{\overline{(-)}}$ un morphisme $\widehat{\overline{R}}  \ldrt \widehat{\overline{\widehat{R}}}$. Dans l'autre sens, le
morphisme $R\ldrt \overline{R}$ induit par complétion $\widehat{R} \ldrt \widehat{\overline{R}}$. Mais d'après le lemme précédent $\widehat{\overline{R}}$ est intégralement fermé dans $
\widehat{\overline{R}}\unpi$, donc le morphisme précédent se prolonge en un morphisme $\overline{\widehat{R}} \ldrt \widehat{\overline{R}}$. Par complétion $\pi$-adique ce morphisme induit un morphisme $\widehat{\overline{\widehat{R}}} \ldrt \widehat{\overline{R}}$. 
\\
On vérifie que les deux morphismes précédents sont inverses l'un de l'autre.
\qed

\begin{rema}
Dans le lemme et corollaire précédent on ne suppose pas que $R$ est séparé pour la topologie $\pi$-adique. De plus même si $R$ est $\pi$-adique en général $\overline{R}$ n'est pas forcément séparé. 
\end{rema}

\subsection{Rappels sur les schémas formels $\pi$-adiques}

\begin{defi}
On appelle schéma formel $\pi$-adique (sous-entendu sur $\spf (\breve{\O})$) un schéma formel $\mathfrak{Z}$ sur $\spf (\breve{\O})$ tel que $\pi\O_\mathfrak{Z}$ soit un idéal de définition de $\mathfrak{Z}$. La catégorie des schémas formels $\pi$-adiques est donc équivalente à la 2-limite projective de la catégorie fibrée des schémas sur $(\spec (\breve{O}/\pi^k\breve{O}))_{k\geq 1}$. Cela signifie que se donner un schéma formel $\pi$-adique est équivalent à se donner une famille $(Z_k)_{k\geq 1}$ où $Z_k$ est un $\spec (\breve{\O}/\pi^k\breve{\O})$-schéma munie d'isomorphismes $Z_{k+1}\otimes \breve{\O}/\pi^k \breve{\O}\iso Z_k$ satisfaisant une condition de cocyle évidente (un schéma formel $\pi$-adique n'est rien d'autre qu'un cas particulier d'ind-schéma). 
\end{defi}

Pour $\mathfrak{Z}$ un schéma formel $\pi$-adique on note $\O_\mathfrak{Z}\unpi$ le faisceau associé au préfaisceau $\mathfrak{U}\mapsto \O_\mathfrak{Z} (\mathfrak{U})\unpi$. Cela signifie que pour $\mathfrak{U}$ un ouvert quasicompact $\GG (\mathfrak{U}, \O_{\mathfrak{Z}} \unpi) = \GG (\mathfrak{U},\O_\mathfrak{Z})\unpi$. 

\begin{exem}
Si $\mathfrak{Z} = \coprod_\N \spf (\breve{\O})$ alors 
$\GG ( \mathfrak{Z}, \O_\mathfrak{Z})\unpi = ( \breve{\O}^\N )\unpi \subsetneq 
\breve{\O}\unpi^\N = \GG ( \mathfrak{Z}, \O_\mathfrak{Z}\unpi )$.
\end{exem}

\begin{defi}
 On dit que $\mathfrak{Z}$ est sans $\pi$-torsion si le faisceau
 $\O_{\mathfrak{Z}}$ l'est c'est à dire $\O_\mathfrak{Z}\xrig{\;\times
   \pi\;} \O_\mathfrak{Z}$ est un monomorphisme. 
\end{defi} 

\begin{lemm}
La schéma formel $\pi$-adique $\mathfrak{Z}$ est sans $\pi$-torsion ssi
il possède un
recouvrement affine $(\spf (R_i))_{i\in I}$ tel que $\forall i\; R_i$
soit sans $\pi$-torsion.
\end{lemm}
\dem
Utiliser le fait que pour $R$ $\pi$-adique
sans $\pi$-torsion et $f\in R$ $\; R<\frac{1}{f}>$ est sans
$\pi$-torsion puisque $R[1/f]$ l'est 
(lemme \ref{carotidegh}). 
\qed

\subsection{Morphismes affines}

\begin{defi}
Un morphisme de schémas formels $\pi$-adiques $\X\ldrt \mathfrak{Y}$ est dit affine si le morphisme de schémas induit entre les fibres spéciales $\X\otimes k \ldrt \mathfrak{Y}\otimes k$ l'est ou encore de façon équivalente si $\forall k$ le morphisme de schémas induit  $\X\otimes \breve{\O}/\pi^k \breve{\O} \ldrt \mathfrak{Y}\otimes \breve{\O}/\pi^k \breve{\O}$ l'est. 
\end{defi} 

Ainsi si $\X$ est un schéma formel $\pi$-adique la catégorie des $\X$-schémas formels $\pi$-adiques finis est équivalente à la catégories des faisceaux de $\O_\X$-algèbres $\mathcal{A}$ tels que l'application canonique $\mathcal{A}\ldrt \underset{k}{\limp} \mathcal{A}/\pi^k \mathcal{A}$ est un isomorphisme et $\forall k\; \mathcal{A}/\pi^k \mathcal{A}$ est une $\O_{\X\otimes \breve{\O}/\pi^k\breve{\O}}$-algèbre quasi-cohérente. 

\subsection{Limite projective dans la catégorie des schémas formels $\pi$-adiques}

\begin{prop}
Soit $(I, \geq )$ un ensemble ordonné filtrant décroissant et $((\mathfrak{Z}_i)_{i\in I}, (\ph_{ij})_{i\geq j})$ un système projectif de schémas formels $\pi$-adiques tel que les morphismes de transition $\ph_{ij}$ soient affines. Alors $\underset{i\in I}{\limp} \mathfrak{Z}_i$ existe dans la catégorie des $spf( \breve{\O})$-schémas formels 
et c'est un schéma formel $\pi$-adiques égal à 
$$
\underset{k\in \N}{\limi} \underset{i\in I}{\limp} (\mathfrak{Z}_i\otimes \breve{\O}/\pi^k\breve{\O})
$$
\end{prop}
\dem
La démonstration ne pose pas de problème. On renvoie au chapitre 8 de EGA IV pour les limites projectives de schémas à morphismes de transition affines. 
\qed

\begin{exem}
Si $\mathfrak{Z}_i = \spf (R_i)$ alors $\underset{i}{\limp} \spf (R_i) =\spf (R_\infty)$ où 
$$
R_\infty = \widehat{\underset{i}{\limi} R_i }
$$
Par exemple si on considère la limite projective de revêtements de Kümmer $\DD \leftarrow \dots \leftarrow \DD \leftarrow \dots$ où $\DD =\spf (\breve{\O}<T,T^{-1}>)$ et les morphismes de transition sont tous $t\mapsto t^p$ alors  
$$
R_\infty =\{ \sum_{\a\in \Z \unp} a_\a T^\a\;|\; a_\a\in \breve{\O}\;\;
 a_\a \underset{|\a| \drt +\infty}{\ldrt} 0 \;\;\;  a_\a \underset{v_p(\a)\drt -\infty}{\ldrt} 0 \;\}
$$
où $a_\a\drt 0$ signifie tendre vers $0$ pour la topologie $\pi$-adique.
\end{exem}

\subsection{Normalisation dans la fibre générique}\label{krutofj248}

\begin{defi}
Pour un schéma formel $\pi$-adique sans $\pi$-torsion $\mathfrak{Z}$ on dit que $\O_\mathfrak{Z}$ est intégralement fermé dans $\O_\mathfrak{Z} \unpi$ si pour tout ouvert $\mathfrak{U}$ 
$\GG (\mathfrak{U},\O_\mathfrak{Z})$ est intégralement fermé dans $\GG (\mathfrak{U},\O_\mathfrak{Z}\unpi )$. Cela est équivalent à dire que $\forall z\in \mathfrak{Z}$ $\O_{\mathfrak{Z},z}$ est intégralement fermé dans $\O_{\mathfrak{Z},z} \unpi$ ou bien encore que pour tout ouvert quasicompact $\mathfrak{U}$ $\GG( \mathfrak{U},\O_\mathfrak{Z})$ est intégralement fermé dans $\GG( \mathfrak{U},\O_\mathfrak{Z})\unpi$.
\end{defi}

\begin{lemm}
Soit $\mathfrak{Z}$ un schéma formel $\pi$-adique sans $\pi$-torsion. Alors $\O_\mathfrak{Z}$ est intégralement fermé dans $\O_\mathfrak{Z}\unpi$ ssi il existe un recouvrement affine $(\spf (R_i))_{i}$ de $\mathfrak{Z}$ tel que $\forall i\; R_i$ est intégralement clos dans $R_i\unpi$.
\end{lemm}
\dem
Soit $R$ une $\breve{\O}$-algèbre $\pi$-adique sans $\pi$-torsion. Il suffit
d'utiliser le fait que $R$ intégralement fermé dans $R\unpi$ implique
 que pour $f\in R$ il en est de même pour $R [\frac{1}{f}]$ et donc d'après le lemme \ref{carotidegh} il en est de même pour $R<\frac{1}{f}>$.
\qed

\begin{lemm}
Soit $\mathfrak{Z}$ un schéma formel $\pi$-adique sans
$\pi$-torsion. Soit $\mathcal{A}$ la $\O_{\mathfrak{Z}}$-algèbre
qui est le faisceau associé au préfaisceau qui
à $\mathfrak{U}$ quasicompact associe la fermeture intégrale de $\GG
(\mathfrak{U},\O_\mathfrak{Z})$ dans $\GG
(\mathfrak{U},\O_\mathfrak{Z})\unpi $. Alors, pour tout ouvert affine
$\spf (R)\subset \mathfrak{Z}$ $\; \GG(\spf (R),\mathcal{A})$ est la
fermeture intégrale de $R$ dans $R\unpi$. De plus $\forall k\geq 1\;
\mathcal{A}/\pi^k\mathcal{A}$ est une $\O_{\mathfrak{Z}\otimes
  \breve{\O}/\pi^k\breve{\O}}$-algèbre quasi-cohérente sur le schéma $\mathfrak{Z}\otimes
  \breve{\O}/\pi^k\breve{\O}$. 
\end{lemm}
\dem
Soit $\spf (R)\subset \mathfrak{Z}$ un ouvert affine. \'Etant donné que $\mathcal{A} \subset
\O_\mathfrak{Z}\unpi$ et que $\spf (R)$ est quasicompact $\GG (\spf
(R),\mathcal{A})\subset R\unpi$. Soit donc $s \in\GG (\spf
(R),\mathcal{A})$.  Il existe un recouvrement affine fini $\spf (R)=
\cup_{i\in I} D(f_i)$ où $\forall i\; f_i\in R$ tel que $\forall i\;
s_{|D(f_i)}$ vérifie
$$
s_{|D(f_i)} \in \overline{R<1/f_i>}
$$
Mais d'après le lemme \ref{carotidegh} pour tout $i\in I$ 
l'image réciproque de $ \overline{R<\frac{1}{f_i}>}$ dans
$R[\frac{1}{f_i}]\unpi$ par l'application $R[\frac{1}{f_i}]\unpi \ldrt R<\frac{1}{f_i}>\unpi$
est égale à $\overline{R[\frac{1}{f_i}]}$. 
On en déduit grâce au recouvrement affine 
$\spec (R) =\cup_{i\in I} D'(f_i)$ où cette fois-ci $D'(f_i)=\spec
(R[\frac{1}{f_i}])$ que l'élément de $R\unpi$ associé à $s$ est 
localement sur $\spec (R)$ entier sur $\O_{\spec (R)}$, donc entier
sur $R$ (cf. par exemple la démonstration de la proposition 6.1.4 de
EGA II p.110). On a donc démontré que $\GG ( \spf (R),
\mathcal{A})=\overline{R}$. 
\\
Montrons maintenant que $\forall k\geq 1 \;
\mathcal{A}/\pi^k\mathcal{A}$ est quasi-cohérente. Soit donc $spf (R)
\subset \mathfrak{Z}$ un ouvert affine et pour un $f\in R$ l'ouvert
$D(f) =\spf ( R<\frac{1}{f}>)$. Soit $A=R[1/f]$. On a $\overline{A} =
\overline{R} [1/f]$. Donc par application du corollaire \ref{mutk256}
à l'anneau $A$
$$
\widehat{\overline{R}}<\frac{1}{f}> = \widehat{ \overline{ R<1/f>}}
$$ 
Donc $\forall k\geq 1$ par tensorisation de l'égalité précédente par
$\breve{\O}/\pi^k \breve{\O}$ 
$$
\overline{R}/\pi^k \overline{R} [\frac{1}{f}] = \overline{R<1/f>}
/\pi^k  \overline{R<1/f>}
$$
Soit $\mathcal{B}$ le préfaisceau $\mathfrak{U} \mapsto \mathcal{A} (
\mathfrak{U}) /\pi^k \mathcal{A} ( \mathfrak{U})$. L'égalité
précédente s'écrit 
$$ 
\GG ( \spf (R), \mathcal{B}) [1/f] =
\GG ( D(f), \mathcal{B})
$$
Donc le faisceau associé $\mathcal{A}$ est quasi-cohérent. 
\qed

\begin{prop}
Soit $\mathfrak{Z}$ un schéma formel $\pi$-adique sans
$\pi$-torsion. Soit $\mathcal{A}$ la $\O_{\mathfrak{Z}}$-algèbre
définie dans le lemme précédent par normalisation de
$\O_{\mathfrak{Z}}$ dans $\O_\mathfrak{Z}\unpi$. Soit 
$$
\mathfrak{Z}^{norm} = \underset{k}{\limi} \spec_{\mathfrak{Z}\otimes \breve{\O}/\pi^k\breve{\O}} ( \mathcal{A}/ \pi^k \mathcal{A})
$$
qui est un $\mathfrak{Z}$-schéma formel affine.
Il vérifie $\O_{\mathfrak{Z}^{norm}}$ est intégralement fermé dans
$\O_{\mathfrak{Z}^{norm}}\unpi$. 
 De plus pour tout schéma
formel $\pi$-adique sans $\pi$-torsion $\mathfrak{X}$ tel que
$\O_\X$ soit intégralement fermé dans $\O_\X\unpi$ pour tout
morphisme $\X\ldrt \mathfrak{Z}$ il existe une unique extension 
$$
\xymatrix{
&\mathfrak{Z}^{norm} \ar[d] \\
\X\ar[ru]^\exists \ar[r] & \mathfrak{Z}
}
$$
\end{prop}

\begin{defi}
Le schéma formel précédent $\mathfrak{Z}^{norm}$ est appelé le normalisé de $\mathfrak{Z}$
dans sa fibre générique. 
\end{defi}

\begin{exem}
Soit $\mathfrak{Z}$ un schéma formel localement de type fini sur $\spf
(\breve{\O})$ sans $\pi$-torsion tel que $\mathfrak{Z}^{rig}$ soit
normal. Alors on a vu dans l'appendice A de \cite{Cellulaire} que 
$\mathfrak{Z}^{norm} = \spf (\mathcal{A})$ où $\mathcal{A} = sp_*
\O_{\mathfrak{Z}^{rig}}^0$ est une $\O_{\mathfrak{Z}}$-algèbre
cohérente. La proposition précédente généralise donc cette situation
sans condition de finitude. 
\end{exem}

\subsection{Commutation de la normalisation dans la fibre générique et
  du passage à la limite projective}

\begin{prop}\label{dljbyryt}
Soit $(I,\geq)$ un ensemble ordonné filtrant décroissant et
$((\mathfrak{Z}_i)_{i\in I},(\ph_{ij})_{i\geq j})$ un système
projectif de schémas formels $\pi$-adiques sans $\pi$-torsion dont les
morphismes de transition sont affines. Il y a alors un isomorphisme
canonique
$$\underset{i}{\limp} \mathfrak{Z}_i^{norm}\iso 
\left ( \underset{i}{\limp} \mathfrak{Z}_i \right )^{norm} 
$$
\end{prop}
\dem
Il y a un morphisme de systèmes projectifs 
$(\mathfrak{Z}_i^{norm})_i \ldrt (\mathfrak{Z}_i)_i$ qui induit un
morphisme affine
$$
\underset{i}{\limp} \mathfrak{Z}_i^{norm} \ldrt \underset{i}{\limp}
\mathfrak{Z}_i 
$$
D'après le lemme \ref{carotidegh} le schéma formel
$\underset{i}{\limp} \mathfrak{Z}_i^{norm}$ est ``intégralement fermé
dans sa fibre générique'' (au sens où si $\O$ est son faisceau
structural $\O$ est intégralement fermé dans $\O\unpi$).  
D'après la propriété universelle du normalisé le morphisme précédent
s'étend donc en
un morphisme 
$$
\underset{i}{\limp} \mathfrak{Z}_i^{norm} \ldrt \left ( \underset{i}{\limp}
\mathfrak{Z}_i \right )^{norm}
$$
Construisons un inverse à ce morphisme. Pour tout $j\in I$ il y a un
morphisme composé 
$$
\left ( \underset{i}{\limp}
\mathfrak{Z}_i \right )^{norm} \ldrt \underset{i}{\limp}
\mathfrak{Z}_i \xrig{\text{ projection }} \mathfrak{Z}_j
$$
qui s'étend par propriété universelle du normalisé en un morphisme
$$
\left ( \underset{i}{\limp}
\mathfrak{Z}_i \right )^{norm} \ldrt \mathfrak{Z}_j^{norm}
$$
Ces morphismes sont compatibles lorsque $j$ varie et fournissent donc
par la propriété universelle de la limite projective un morphisme 
$$
\left ( \underset{i}{\limp}
\mathfrak{Z}_i \right )^{norm} \ldrt  \underset{i}{\limp}  \mathfrak{Z}_i^{norm}
$$
On vérifie facilement que ces deux morphismes sont inverses l'un de
l'autre.
\qed

\subsection{\'Eclatements formels admissibles}

\subsubsection{Définition et premières propriétés}

\begin{defi}\label{sdghidgpp}
Soit $\mathfrak{Z}$ un schéma formel $\pi$-adique sans $\pi$-torsion
et $\mathcal{I}\subset \O_{\mathfrak{Z}}$ un idéal tel que localement
sur $\mathfrak{Z}$ $\;\exists N\in \N\; \pi^N\O_{\mathfrak{Z}}\subset
\mathcal{I}$ et $\mathcal{I}/\pi^N \O_{\mathfrak{Z}}$ est
quasi-cohérent de type fini. Un tel idéal est dit admissible.
 On appelle éclatement formel admissible
de $\mathcal{I}$ le $\mathfrak{Z}$-schéma formel $\pi$-adique 
$$
\widetilde{\mathfrak{Z}}= \underset{k}{\limi} \text{Proj} \left (
  \bigoplus_{i\geq 0} \mathcal{I}^i/\pi^k \mathcal{I}^i \right )
$$
\end{defi}

\begin{prop}
Avec les notations de la définition précédente 
\begin{itemize}
\item si $\spf (R)\subset \mathfrak{Z}$ est un ouvert affine et 
  $I=\GG(\spf (R),\mathcal{I})$ alors $\widetilde{\mathfrak{Z}}_{|\spf
    (R)}$ s'identifie au complété $\pi$-adique de l'éclatement de
  l'idéal $\widetilde{I}$ de $\spec (R)$ 
\item $\widetilde{\mathfrak{Z}}$ est sans $\pi$-torsion
\item si $\ph : \widetilde{\mathfrak{Z}}\ldrt \mathfrak{Z}$ alors
  $\O_{\widetilde{\mathfrak{Z}}}. \ph^{-1} \mathcal{I}$ est localement
  libre de rang un 
\item $\widetilde{\mathfrak{Z}}$ satisfait à la propriété universelle suivante :
  pout tout $\mathfrak{Z}$-schéma formel $\pi$-adique sans
  $\pi$-torsion $\mathfrak{Y}\xrig{\; \psi\;} \mathfrak{Z}$ tel que
  $\O_{\mathfrak{Y}}.\psi^{-1} \mathcal{I}$ soit localement libre de
  rang un il existe un unique $\mathfrak{Z}$-morphisme
  $\mathfrak{Y}\ldrt \widetilde{\mathfrak{Z}}$
\end{itemize}
\end{prop}
\dem
La première assertion découle de la définition de l'éclatement formel
admissible.
\\
La seconde résulte de la première car $\spec (R)$ étant sans
$\pi$-torsion l'éclatement de l'idéal $\widetilde{I}$ l'est aussi
(l'image réciproque d'un ouvert schématiquement dense par un
éclatement reste schématiquement dense) et donc d'après le lemme
\ref{carotidegh} son complété $\pi$-adique est encore sans
$\pi$-torsion.
\\
La troisième résulte également de la première. En effet, sur l'éclaté 
de $\widetilde{I}$ dans le schéma $\spec (R)$ l'idéal $\widetilde{I}$
est localement libre de rang un. Il est donc localement monogène sur
le complété $\pi$-adique de ce schéma. Donc $\mathcal{I}$ devient
localement monogène sur $\widetilde{\mathfrak{Z}}$. Mais étant donné
que $\mathcal{I}$ contient localement une puissance de $\pi$ et que
$\widetilde{\mathfrak{Z}}$ est sans $\pi$-torsion $\mathcal{I}$
est localement libre de rang un sur $\widetilde{\mathcal{Z}}$.
\\
La dernière assertion résulte aisément de son homologue pour les
schémas (et de la première assertion).
\qed

\begin{rema}
Soit $\mathfrak{Z} = \spf (R)$ $\pi$-adique sans $\pi$-torsion 
 et $\mathcal{I}$ un idéal admissible de $\O_\mathfrak{Z}$.
Il y a alors un idéal 
 $I = (f_1,\dots,f_n)$  de $R$ contenant une puissance de $\pi$ tel
 que  
 $\mathcal{I}$ soit l'image réciproque  de l'idéal quasi-cohérent
$\widetilde{ (I/\pi^N \O_\mathfrak{Z})}$ dans $\O_{\mathfrak{Z}}$ 
pour $N>>0$. 
La description donnée dans le lemme 2.2 de \cite{BLI} de
$\widetilde{\mathfrak{Z}}$ lorsque $\mathfrak{Z}$ est topologiquement
de type fini sur $\breve{\O}$ est en général fausse. On a en fait la
description suivante :
$$
\widetilde{\mathfrak{Z}} = \bigcup_{i=1}^n \mathcal{U}_i
$$
où $\mathcal{U}_i = \spf (A_i)$ est un ouvert affine tel que si
$$
B_i = R<T_1,\dots, \widehat{T_i},\dots, T_n>/(T_j f_i-f_j)_{1\leq j\leq
  n,j\neq i}
$$
et $B'_i =B_i /J_i$ avec $J_i =\{b\in B_i\;|\; \exists k\; \pi^k
b=0\}$ alors $A_i$ est le séparé de $B'_i$, $A_i= B'_i /\cap_{k\geq 0}
\pi^k B'_i$. 
\end{rema}

\begin{lemm}
Soit $\mathfrak{Z}$ comme précédemment.
Soient $\mathcal{I}_1, \mathcal{I}_2$ deux idéaux admissibles.
 Soit $\widetilde{\mathfrak{Z}}$ l'éclatement
formel de $\mathcal{I}_1$. 
Alors l'éclatement formel admissible de
$\mathcal{I}_1.\mathcal{I}_2$ s'identifie à l'éclatement formel
admissible de l'image réciproque de $\mathcal{I}_2$ à
$\widetilde{\mathfrak{Z}}$.  
\end{lemm}

\begin{rema}
On utilisera souvent la propriété suivante. Soit $\mathfrak{Z}$
$\pi$-adique sans $\pi$-torsion et $\mathfrak{I}\subset
\O_\mathfrak{Z}$ un idéal admissible. Alors $\mathfrak{I}$ est
localement libre de rang un ssi il est localement monogène.  
\end{rema}

\subsubsection{Adhérence ``schématique'' de la fibre générique}

Soit $\mathfrak{Z}$ un schéma formel $\pi$-adique. Notons pour tout
$k\geq 1$
$$
\mathcal{I}_k = \underset{i\geq 1}{\limi} \ker \left (
  \O_\mathfrak{Z}/\pi^k\O_\mathfrak{Z} \xrig{\; \times \pi^i\;}
  \O_\mathfrak{Z}/\pi^{k+i}\O_\mathfrak{Z}\right ) \subset
\O_{\mathfrak{Z}}/\pi^{k} \O_\mathfrak{Z}
$$
un faisceau d'idéaux quasicohérent sur $\mathfrak{Z}\otimes
\breve{\O}/\pi^k$. 
Si $\spf (R)\subset  \mathfrak{Z}$ et $I=\{x\in R\;|\; \exists i\geq
1\; \pi^i x=0 \}$ alors 
$$
\mathcal{I}_k = (I+\pi^k R/\pi^kR)^{\widetilde{\;\;\;\;\;}}
$$
où le tilda signifie ``le faisceau quasicohérent associé''. 
Notons alors 
$$
Z_k = V(\mathcal{I}_k) \subset  \mathfrak{Z}\otimes\breve{\O}/\pi^k\breve{\O}
$$
On a donc
$$
Z_{k+1} \otimes \breve{\O}/\pi^k\breve{\O} = Z_k
$$
Notons alors
$$
\mathfrak{Z}' = \underset{k}{\limi} Z_k
$$
un schéma formel $\pi$-adique. Si $spf (R)$ est un ouvert affine de
$\mathfrak{Z}$ et $I$ est l'idéal des éléments de $\pi^\infty$-torsion
comme précédemment alors l'ouvert correspondant de $\mathfrak{Z}'$ est
$\spf ( R/\overline{I})$ où $\overline{I}$ désigne l'adhérence de $I$
pour la topologie $\pi$-adique.

\begin{lemm}
Le schéma formel $\pi$-adique 
$\mathfrak{Z}'$ est sans $\pi$-torsion. De plus ``l'immersion fermée''
$\mathfrak{Z}' \hookrightarrow \mathfrak{Z}$ est telle que pour tout
schéma formel $\pi$-adique sans $\pi$-torsion $\mathfrak{Y}$  tout
morphisme 
$\mathfrak{Y} \ldrt \mathfrak{Z}$ se factorise via
$\mathfrak{Z}'\hookrightarrow \mathfrak{Z}$. 
\end{lemm}
\dem 
Elle ne pose pas de problème particulier.
\qed

\begin{defi}
Par abus de terminologie 
on appellera $\mathfrak{Z}'$ l'adhérence schématique de la fibre
générique de $\mathfrak{Z}$.
\end{defi}

\subsubsection{Transformée stricte} \label{sdguiec23}

Soit $\ph : \mathfrak{Y}\ldrt \mathfrak{Z}$ un morphisme de schémas
formels $\pi$-adiques sans $\pi$-torsion. Soit $\mathcal{I}\subset
\O_\mathfrak{Z}$ un faisceau d'idéaux satisfaisant aux hypothèses de
la définition \ref{sdghidgpp}. Notons $\widetilde{\mathfrak{Z}}\ldrt
\mathfrak{Z}$ l'éclatement formel admissible associé. 

\begin{defi}
On appelle transformé strict de $\mathfrak{Y}$ relativement à
l'éclatement $\widetilde{\mathfrak{Z}}\ldrt
\mathfrak{Z}$ l'adhérence schématique de la fibre générique de 
$\mathfrak{Y}\times_\mathfrak{Z} \widetilde{\mathfrak{Z}}$.
\end{defi}

\begin{prop}\label{subbermpzq}
Le transformé strict de $\mathfrak{Y}$ s'identifie à l'éclatement
formel admissible de l'idéal $\O_\mathfrak{Y}.\ph^{-1}\mathcal{I}$. 
\end{prop}
\dem
Notons $\X$ le transformé strict et $\widetilde{\mathfrak{Y}}$
l'éclatement formel de $\O_\mathfrak{Y}.\ph^{-1}\mathcal{I}$.
Puisqu'il y a une factorisation $ \mathfrak{Y}\times_{\mathfrak{Z}}
\widetilde{\mathfrak{Z}} \ldrt \widetilde{\mathfrak{Z}}\ldrt
\mathfrak{Z}$ 
l'image réciproque à $\mathfrak{Y}\times_{\mathfrak{Z}}
\widetilde{\mathfrak{Z}}$ de l'idéal $\mathcal{I}$ est localement
monogène. Donc, puisque cet idéal contient localement une puissance de
$\pi$ et puisque $\X\hookrightarrow \mathfrak{Y}\times_{\mathfrak{Z}}
\widetilde{\mathfrak{Z}}$
 est sans $\pi$-torsion son image réciproque à $\X$ est localement
 libre de rang un. Donc via le morphisme composé 
$\X \ldrt  \mathfrak{Y}\times_{\mathfrak{Z}}
\widetilde{\mathfrak{Z}} \ldrt \mathfrak{Y}$ l'image réciproque de 
$\O_\mathfrak{Y}.\ph^{-1} \mathcal{I}$ est localement libre de rang
un.
D'après la propriété universelle de $\widetilde{\mathfrak{Y}}$ il y a
donc un morphisme
$$
\X\ldrt \widetilde{\mathfrak{Y}}
$$
Construisons un morphisme dans l'autre sens. D'après la propriété
universelle de $\widetilde{\mathfrak{Z}}$ le morphisme
$\widetilde{\mathfrak{Y}}\ldrt \mathfrak{Z}$ s'étend en un morphisme
$\widetilde{\mathfrak{Y}}\ldrt \widetilde{\mathfrak{Z}}$. Il fournit
donc un morphisme
$$
\widetilde{\mathfrak{Y}}\ldrt  \mathfrak{Y}\times_{\mathfrak{Z}}
\widetilde{\mathfrak{Z}}
$$
Mais grâce à la propriété caractérisant l'adhérence schématique de la
fibre générique se morphisme se factorise en un morphisme
$$
\widetilde{\mathfrak{Y}}\ldrt \X
$$
On vérifie alors facilement que les deux morphismes précédents sont
inverses l'un de l'autre.
\qed

\subsubsection{Commutation à la limite projective}

Soit $(I,\geq)$ un ensemble ordonné filtrant décroissant et et
$((\mathfrak{Z}_i)_{i\in I}, (\ph_{ij})_{i\geq j})$ un système
projectif de schémas formels $\pi$-adiques sans $\pi$-torsion tel que les morphismes de
transition $\ph_{ij}$ soient affines. Soit $i_0\in I$ fixé et
$\mathcal{I}$ un idéal admissible de $\O_{\mathfrak{Z}_{i_0}}$. 
Notons pour $i\geq i_0$ $\widetilde{\mathfrak{Z}}_i$ l'éclatement
formel de l'image réciproque de $\mathcal{I}$. On a donc un système
projectif $(\widetilde{\mathfrak{Z}}_i)_{i\in I}$. D'après la
proposition \ref{subbermpzq} les morphismes de transition sont
affines.

\begin{prop}
La limite projective $\underset{i\geq i_0}{\limp}
\widetilde{\mathfrak{Z}}_i$ coïncide avec l'éclatement formel de
l'image réciproque à $\underset{i}{\limp} \mathfrak{Z}_i$ de $\mathcal{I}$.
\end{prop}

\begin{coro}\label{sdmlkgiob35}
Soit pour tout $i$ $\widetilde{\mathfrak{Z}}_i^{norm}$ le normalisé
dans sa fibre générique de l'éclaté $\widetilde{\mathfrak{Z}}_i$. Il y
a alors une identification entre $\underset{i\geq i_0}{\limp}
\widetilde{\mathfrak{Z}}_i^{norm}$ et le normalisé 
dans sa fibre générique de l'éclatement
formel de $\underset{i}{\limp} \mathfrak{Z}_i$. 
\end{coro}
\dem
Appliquer la proposition précédente couplée à la proposition
\ref{dljbyryt}.
\qed

\bibliographystyle{plain}
\bibliography{biblio}
\end{document}

%% file: desol.pstex_t
\begin{picture}(0,0)%
\includegraphics{desol.pstex}%
\end{picture}%
\setlength{\unitlength}{1381sp}%
\begingroup\makeatletter\ifx\SetFigFont\undefined%
\gdef\SetFigFont#1#2#3#4#5{%
  \reset@font\fontsize{#1}{#2pt}%
  \fontfamily{#3}\fontseries{#4}\fontshape{#5}%
  \selectfont}%
\fi\endgroup%
\begin{picture}(7484,5502)(301,-7381)
\put(601,-2161){\makebox(0,0)[lb]{\smash{{\SetFigFont{5}{6.0}{\rmdefault}{\mddefault}{\updefault}{\color[rgb]{0,0,0}$\widetilde{\mathfrak{X}}_\infty$}%
}}}}
\put(6226,-7111){\makebox(0,0)[lb]{\smash{{\SetFigFont{5}{6.0}{\rmdefault}{\mddefault}{\updefault}{\color[rgb]{0,0,0}$\widehat{\Omega}$}%
}}}}
\put(7051,-6211){\makebox(0,0)[lb]{\smash{{\SetFigFont{5}{6.0}{\rmdefault}{\mddefault}{\updefault}{\color[rgb]{0,0,0}$\displaystyle{\coprod_{\mathbb{Z}}\widehat{\Omega}}$}%
}}}}
\put(1801,-7336){\makebox(0,0)[lb]{\smash{{\SetFigFont{5}{6.0}{\rmdefault}{\mddefault}{\updefault}{\color[rgb]{0,0,0}$\widehat{\mathbb{P}}^{n-1}$}%
}}}}
\put(1051,-5536){\makebox(0,0)[lb]{\smash{{\SetFigFont{5}{6.0}{\rmdefault}{\mddefault}{\updefault}{\color[rgb]{0,0,0}$\mathbb{D}_a$}%
}}}}
\put(7051,-2536){\makebox(0,0)[lb]{\smash{{\SetFigFont{5}{6.0}{\rmdefault}{\mddefault}{\updefault}{\color[rgb]{0,0,0}$\widetilde{\mathcal{Y}}_\infty$}%
}}}}
\put(7078,-3243){\makebox(0,0)[lb]{\smash{{\SetFigFont{5}{6.0}{\rmdefault}{\mddefault}{\updefault}{\color[rgb]{0,0,0}$\mathcal{Y}_\infty$}%
}}}}
\put(1426,-6586){\makebox(0,0)[lb]{\smash{{\SetFigFont{5}{6.0}{\rmdefault}{\mddefault}{\updefault}{\color[rgb]{0,0,0}$\widetilde{\widehat{\mathbb{P}}^{n-1}}$}%
}}}}
\put(7051,-2011){\makebox(0,0)[lb]{\smash{{\SetFigFont{5}{6.0}{\rmdefault}{\mddefault}{\updefault}{\color[rgb]{0,0,0}$\widetilde{\widetilde{\mathcal{Y}}}_\infty$}%
}}}}
\put(1051,-2911){\makebox(0,0)[lb]{\smash{{\SetFigFont{5}{6.0}{\rmdefault}{\mddefault}{\updefault}{\color[rgb]{0,0,0}$\scriptstyle{\supset}\mathbb{D}_{a,\infty}$}%
}}}}
\put(451,-2911){\makebox(0,0)[lb]{\smash{{\SetFigFont{5}{6.0}{\rmdefault}{\mddefault}{\updefault}{\color[rgb]{0,0,0}$\mathfrak{X}_{\scriptscriptstyle{\infty}}$}%
}}}}
\end{picture}%

%% file: diagramme_final.pstex_t
\begin{picture}(0,0)%
\includegraphics{diagramme_final.pstex}%
\end{picture}%
\setlength{\unitlength}{3947sp}%
\begingroup\makeatletter\ifx\SetFigFont\undefined%
\gdef\SetFigFont#1#2#3#4#5{%
  \reset@font\fontsize{#1}{#2pt}%
  \fontfamily{#3}\fontseries{#4}\fontshape{#5}%
  \selectfont}%
\fi\endgroup%
\begin{picture}(7484,5502)(301,-7381)
\put(601,-2911){\makebox(0,0)[lb]{\smash{{\SetFigFont{10}{12.0}{\rmdefault}{\mddefault}{\updefault}{\color[rgb]{0,0,0}$\mathfrak{X}_\infty$}%
}}}}
\put(601,-2161){\makebox(0,0)[lb]{\smash{{\SetFigFont{10}{12.0}{\rmdefault}{\mddefault}{\updefault}{\color[rgb]{0,0,0}$\widetilde{\mathfrak{X}}_\infty$}%
}}}}
\put(6226,-7111){\makebox(0,0)[lb]{\smash{{\SetFigFont{10}{12.0}{\rmdefault}{\mddefault}{\updefault}{\color[rgb]{0,0,0}$\widehat{\Omega}$}%
}}}}
\put(7051,-6211){\makebox(0,0)[lb]{\smash{{\SetFigFont{10}{12.0}{\rmdefault}{\mddefault}{\updefault}{\color[rgb]{0,0,0}$\displaystyle{\coprod_{\mathbb{Z}}\widehat{\Omega}}$}%
}}}}
\put(866,-2911){\makebox(0,0)[lb]{\smash{{\SetFigFont{10}{12.0}{\rmdefault}{\mddefault}{\updefault}{\color[rgb]{0,0,0}$\supset\mathbb{D}_{a,\infty}$}%
}}}}
\put(1801,-7336){\makebox(0,0)[lb]{\smash{{\SetFigFont{10}{12.0}{\rmdefault}{\mddefault}{\updefault}{\color[rgb]{0,0,0}$\widehat{\mathbb{P}}^{n-1}$}%
}}}}
\put(1051,-5536){\makebox(0,0)[lb]{\smash{{\SetFigFont{10}{12.0}{\rmdefault}{\mddefault}{\updefault}{\color[rgb]{0,0,0}$\mathbb{D}_a$}%
}}}}
\put(7051,-2536){\makebox(0,0)[lb]{\smash{{\SetFigFont{10}{12.0}{\rmdefault}{\mddefault}{\updefault}{\color[rgb]{0,0,0}$\widetilde{\mathcal{Y}}_\infty$}%
}}}}
\put(7078,-3243){\makebox(0,0)[lb]{\smash{{\SetFigFont{10}{12.0}{\rmdefault}{\mddefault}{\updefault}{\color[rgb]{0,0,0}$\mathcal{Y}_\infty$}%
}}}}
\put(1426,-6586){\makebox(0,0)[lb]{\smash{{\SetFigFont{10}{12.0}{\rmdefault}{\mddefault}{\updefault}{\color[rgb]{0,0,0}$\widetilde{\widehat{\mathbb{P}}^{n-1}}$}%
}}}}
\put(7051,-2011){\makebox(0,0)[lb]{\smash{{\SetFigFont{10}{12.0}{\rmdefault}{\mddefault}{\updefault}{\color[rgb]{0,0,0}$\widetilde{\widetilde{\mathcal{Y}}}_\infty$}%
}}}}
\end{picture}%